\documentclass{article}
\usepackage{amsmath,amsthm,amssymb}

\title{Elementary Equivalence\\
of Endomorphism Rings of Abelian $\boldsymbol{p}$-Groups}

\author{E.~I.~Bunina, A.~V.~Mikhalev}

\newtheorem{theorem}{Theorem}
\newtheorem{corollary}{Corollary}
\newtheorem*{corol}{Corollary}
\newtheorem{lemma}{Lemma}
\newtheorem{proposition}{Proposition}
\newtheorem{statement}{Statement}

\makeatletter
\@addtoreset{theorem}{section}
\@addtoreset{lemma}{section}
\@addtoreset{proposition}{section}
\@addtoreset{corollary}{section}
\@addtoreset{statement}{section}
\makeatother

\newcommand{\prisv}{\mathbin{{:}\!=}}

\newcommand{\Dp}{\mathop{\mathrm{Dp}}\nolimits}
\newcommand{\rEndom}{\mathrm{Endom}}
\newcommand{\Base}{\mathrm{Base}}
\newcommand{\Hom}{\mathop{\mathrm{Hom}}\nolimits}
\newcommand{\Idem}{\mathrm{Idem}}
\newcommand{\tEnd}{\widetilde{\mathrm{End}}}
\newcommand{\GL}{\mathrm{GL}}
\newcommand{\PGL}{\mathrm{PGL}}
\newcommand{\PSL}{\mathrm{PSL}}
\newcommand{\SL}{\mathrm{SL}}

\newcommand{\Cn}{\mathrm{Cn}}
\newcommand{\Th}{\mathrm{Th}}
\newcommand{\Endom}{\mathop{\mathrm{End}}\nolimits}

\newcommand{\Dom}{\mathrm{Dom}}
\newcommand{\Un}{\mathrm{Un}}
\newcommand{\On}{\mathrm{On}}
\newcommand{\card}{\mathop{\mathrm{card}}\nolimits}
\newcommand{\cf}{\mathrm{cf}}
\newcommand{\rng}{\mathop{\mathrm{rng}}\nolimits}
\newcommand{\ucup}{\mathop{{\cup}}}

\newcommand{\Fin}{\mathop{\mathrm{Fin}}\nolimits}

\newcommand{\Gr}{\mathrm{Gr}}
\newcommand{\Rng}{\mathop{\mathrm{Rng}}\nolimits}
\newcommand{\Cl}{\mathop{\mathrm{Cl}}\nolimits}

\newcommand{\Func}{\mathrm{Func}}

\newcommand{\Ker}{\mathop{\mathrm{Ker}}\nolimits}
\newcommand{\Image}{\mathop{\mathrm{Im}}\nolimits}

{\begin{list}{}{%
\setlength{\itemsep}{-\parsep}%
\settowidth{\labelwidth}{#1}%
\setlength{\leftmargin}{\labelwidth}%
\addtolength{\leftmargin}{\labelsep}}}%
{\end{list}}

\newcommand{\logic}[1]{\mathrel{#1}}
\newcommand{\losp}{\,}
\newcommand{\ttm}{\mskip-1mu}

\newcounter{lastla}
\newcounter{lastar}
\newcounter{slast}

\allowdisplaybreaks

\begin{document}



\maketitle

\tableofcontents

\section*{Introduction}
\addcontentsline{toc}{section}{\hspace*{2.3em}Introduction}

In this paper, we consider elementary properties
(i.e., properties which are definable in the first order language)
of endomorphism rings of Abelian $p$-groups.

The first result on relationship of elementary properties of some
models with elementary properties of derivative models was proved by
A.~I.~Maltsev in 1961 in~\cite{Maltsev}.
He proved that the groups $G_n(K)$ and $G_m(L)$
(where $G=\GL,\SL,\PGL,\PSL$, $n,m\ge 3$, $K$,~$L$ are fields
of characteristic~$0$)
are elementarily equivalent if and only if $m=n$ and the fields
$K$~and~$L$ are elementarily equivalent.

This theory was continued in 1992 when with the help of the construction
of ultraproduct and the isomorphism theorem~\cite{Keisler} C.~I.~Beidar
and A.~V.~Mikhalev in~\cite{BeiMikh} formulated a~general approach
to problems of elementary equivalence
of various algebraic structures, and generalized Maltsev's theorem
for the case where $K$~and~$L$ are skewfields and associative rings.

In 1998--2001 E.~I.~Bunina continued to study some problems of this type
(see \cite{Bun1,Bun2,Bun3,Bun4}). She generalized the results of
A.~I.~Maltsev for unitary linear groups over skewfields
and associative rings with involution, and also for
Chevalley groups over fields.

In 2000 V.~Tolstykh in~\cite{Tolstyh} studied a~relationship
between second order properties of skewfields
and first order properties of automorphism groups of infinite-dimensional
linear spaces over them. In 2003 (see~\cite{categories})
the authors studied a~relationship between second order properties
of associative rings and first order properties of categories of modules,
endomorphism rings, automorphism groups, and projective spaces
of infinite rank over these rings.

In this paper, we study a~relationship between second order properties
of Abelian $p$-groups and first order
properties of their endomorphism rings.

The first section includes some basic notions from the set theory
and model theory:
definitions of first order and second order languages, models of
a~language, deducibility, interpretability,
basic notions of set theory,
which will be needed in next sections.

The second section contains all notions and statements about Abelian groups
which will be needed for our future constructions.
We have taken them mainly from~\cite{Fuks}.

In the third section, we show how to extend the results of
S.~Shelah from~\cite{Shelah} on interpreting the set theory in
a~category for the case of the endomorphism ring
of some special Abelian $p$-group, which is
a~direct sum of cyclic groups of the same order.

In Sec.~4, we describe the second order group language~$\mathcal L_2$,
and also its restriction $\mathcal L_2^\varkappa$
by some cardinal number~$\varkappa$, and then in Sec.~4.2
we introduce the \emph{expressible rank} $r_{\mathrm{exp}}$
of an Abelian group~$A$,
represented as the direct sum $D\oplus G$ of its divisible and reduced
components, as the maximum of the powers of the group~$D$ and some basic
subgroup~$B$ of~$A$, i.e., $r_{\mathrm{exp}}=\max (r(D), r(B))$.
In Sec.~4.2, we also formulate the main theorem of this work.

{\bf Theorem 1.}
\emph{For any infinite $p$-groups $A_1$~and~$A_2$
elementary equivalence of endomorphism rings $\Endom(A_1)$
and $\Endom(A_2)$ implies coincidence of the second order theories
$\Th_2^{r_{\mathrm{exp}}(A_1)}(A_1)$ and
$\Th_2^{r_{\mathrm{exp}}(A_2)}(A_2)$ of the groups $A_1$~and~$A_2$,
bounded by the cardinal numbers
$r_{\mathrm{exp}}(A_1)$ and $r_{\mathrm{exp}}(A_2)$, respectively.}

\smallskip 

Note that $r_{\mathrm{exp}}(A)=|A|$
in all cases except the case where $|D|<|G|$,
any basic subgroup of~$A$ is countable, and the group~$G$ itself is
uncountable. In this case, $r_{\mathrm{exp}}(A)=\omega$.

In Sec.~4.3, we prove two ``inverse implications''
of the main theorem.

{\bf Theorem 2.}
\emph{For any Abelian groups $A_1$~and~$A_2$, if the groups $A_1$~and~$A_2$
are equivalent in the second order logic~$\mathcal L_2$,
then the rings $\Endom(A_1)$ and $\Endom(A_2)$
are elementarily equivalent.}

\smallskip 

{\bf Theorem 3.}
\emph{If Abelian groups $A_1$~and~$A_2$ are reduced
and their basic subgroups are countable,
then $\Th_2^\omega(A_1)=\Th_2^\omega(A_2)$
implies elementary equivalence of
the rings $\Endom(A_1)$ and $\Endom(A_2)$.}

\smallskip

Therefore for all Abelian groups,
except the case where
$A=D\oplus G$, $D\ne 0$, $|D|< |G|$, and $|G|> \omega$,
a~basic subgroup in~$A$ is countable, and elementary
equivalence of the rings $\Endom(A_1)$ and $\Endom(A_2)$ is equivalent to
$$
\Th_2^{r_{\mathrm{exp}}(A_1)}(A_1) = \Th_2^{r_{\mathrm{exp}}(A_2)}(A_2).
$$

In Sec.~4.4, we divide the proof of the main theorem into three cases:
\begin{enumerate}
\item
$A_1$ and $A_2$ are bounded;
\item
$A_1=D_1\oplus G_1$, $A_2=D_2\oplus G_2$, $D_1$~and~$D_2$
are divisible, $G_1$~and~$G_2$ are bounded;
\item
$A_1$~and~$A_2$ have unbounded basic subgroups.
\end{enumerate}

In Secs.~5--7, these three cases are under consideration.

In Sec.~8, we prove the main theorem,
combining all three cases in one proof.

\section{Basic Notions from Model Theory}

\subsection{First Order Languages}\label{ss1.1}
A~\emph{first order language} $\mathcal L$ is a~collection of symbols.
It consists of
(1)~parentheses ${(}$,~${)}$; (2)~connectives $\land$~(``and'') and
$\neg$~(``not'');
(3)~the quantifier~$\forall$ (for all);
(4)~the binary relation symbol~$=$ (identity);
(5)~a countable set of variables~$x_i$;
(6)~a~finite or countable set of relation symbols~$Q_i^n$ ($n\ge 1$);
(7)~a~finite or countable set of function symbols~$F_i^n$ ($n\ge 1$);
(8)~a~finite or countable set of constant symbols~$c_i$.

Now we introduce
the most important for us examples of first order languages:
the group language $\mathcal L_G$ and the ring language~$\mathcal L_R$.

We assume that in the group language there are neither function nor
constant symbols, and there is a~unique 3-place relation symbol~$Q^3$,
which corresponds to multiplication.
Instead of $Q^3(x_1,x_2,x_3)$
we shall write $x_1=x_2\cdot x_3$, or $x_1=x_2 x_3$.

For the ring language we shall also suppose that there are
neither function nor constant symbols,
and there are two relation symbols:
a~3-place symbol of multiplication~$Q_1^3$
(instead of $Q_1^3(x_1,x_2,x_3)$
we shall write $x_1=x_2\cdot x_3$, or $x_1=x_2 x_3$) and a~3-place
symbol of addition~$Q_2^3$ (instead of $Q_2^3(x_1,x_2,x_3)$ we shall write
$x_1=x_2+x_3$).

A~\emph{symbol-string} is defined as follows:
(1)~every symbol~$\alpha$ of the language~$\mathcal L$ is
a~\emph{symbol-string}; (2)~if $\sigma$~and~$\rho$ are symbol-strings,
than $\sigma\rho$ is a~symbol-string.
A~\emph{designating symbol-string~$\sigma$ for a~symbol-string~$\rho$}
is the symbol-string $\sigma \prisv \rho$, or $\rho \prisv \sigma$
(\emph{$\sigma$~is a~designation for~$\rho$}).
If a~symbol-string~$\rho$ is a~part of a~symbol-string~$\sigma$,
staying in one of the three following positions:
$\ldots \rho$, $\rho \ldots$, $\ldots \rho \ldots$,
then $\rho$ is an \emph{occurrence in~$\sigma$}.

Some symbol-strings constructed from the
symbols of the language~$\mathcal L$ are called
\emph{terms} and \emph{formulas} of this language.

\emph{Terms} are defined as follows:
\begin{enumerate}
\item
a~variable is a~term;
\item
a~constant symbol is a~term;
\item
if $F^n$ is an $n$-place function symbol
and
$t_1,\dots,t_{n}$ are terms, then
$F^n(t_1,\dots,t_{n})$ is a~term;
\item
a~symbol-string is a~term only if it can be shown to be a~term
by a~finite number of applications of (1)--(3).
\end{enumerate}

In the cases of the languages
$\mathcal L_G$~and~$\mathcal L_R$,
terms have the form~$x_i$.

The \emph{elementary formulas} of the
language~$\mathcal L$ are symbol-strings of the form
given below:
\begin{enumerate}
\item
if $t_1$ and $t_2$ are terms of the language~$\mathcal L$, then $t_1=t_2$
is an elementary formula;
\item
if $Q^n$ is an $n$-place relation symbol
and $t_1,\dots,t_{n}$ are terms,
then the symbol-string $Q^n(t_1,\dots,t_{n})$ is an elementary formula.
\end{enumerate}

For the language $\mathcal L_G$ the elementary formulas have the form
$x_i=x_j$ and $x_i=x_j\cdot x_k$,
and for the language~$\mathcal L_R$ they have the form
$x_i=x_j$, $x_i=x_j\cdot x_k$, and $x_i=x_j+x_k$.

Finally, the \emph{formulas}
of the language~$\mathcal L$ are defined as follows:
\begin{enumerate}
\item
an elementary formula is a~formula;
\item
if $\varphi$~and~$\psi$ are formulas and $x$~is a~variable,
then $(\neg \varphi)$, $(\varphi\land \psi)$, and
$(\forall x \losp \varphi)$ are formulas;
\item
a~symbol-string is a~formula
only if it can be shown to be a~formula by a~finite number
of applications of (1)--(2).
\end{enumerate}

Let us introduce the following abbreviations:
\begin{itemize}
\item[]
$(\varphi\lor \psi)$ stands for $(\neg ((\neg \varphi)\land(\neg \psi)))$;
\item[]
$(\varphi\Rightarrow \psi)$ stands for $((\neg \varphi)\lor \psi)$;
\item[]
$(\varphi\Leftrightarrow \psi)$ stands for
$((\varphi\Rightarrow \psi)\land (\psi\Rightarrow \varphi))$;
\item[]
$(\exists x\losp \varphi)$ is an abbreviation for
$(\neg (\forall x\losp (\neg \varphi)))$;
\item[]
$\varphi_1\vee \varphi_2\vee\dots \vee\varphi_n$ stands for
$(\varphi_1\vee (\varphi_2\vee \dots \vee \varphi_n))$;
\item[]
$\varphi_1\wedge \varphi_2\wedge \dots \wedge \varphi_n$ stands for
$(\varphi_1\wedge (\varphi_2\wedge \dots \wedge \varphi_n))$;
\item[]
$\forall x_1\dots \forall x_n \varphi$ stands for
$(\forall x_1)\dots (\forall x_n)\varphi$;
\item[]
$\exists x_1\dots \exists x_n \varphi$ stands for
$(\exists x_1)\dots (\exists x_n)\varphi$.
\end{itemize}

Let us introduce the notions of \emph{free} and \emph{bound} occurrences
of a~variable in a~formula.

\begin{enumerate}
\renewcommand{\labelenumi}{\theenumi.}
\item
All occurrences of all variables in elementary formulas
are free occurrences.
\item
Every free (bound) occurrence of a~variable~$x$ in a~formula~$\varphi$
is a~free (bound) occurrence of the variable~$x$ in the formulas
$(\neg \varphi)$, $(\varphi\land \psi)$, and $(\psi \land \varphi)$.
\item
For any occurrence of a~variable~$x$ in a~formula~$\varphi$,
the occurrence of the variable~$x$ in the formula $\forall x\losp \varphi$
is bound.
If an occurrence of a~variable~$x$ in a~formula~$\varphi$ is free (bound),
then the occurrence of~$x$ in $\forall x'\losp \varphi$ is free (bound).
\end{enumerate}

Therefore, one variable can have free and bound occurrences in the
same formula. A~variable is called a~\emph{free}
(\emph{bound}) \emph{variable} in a~given formula if there exist
free (bound) occurrences of this variable in this formula.
Thus a~variable can be free and bound in the same time.

A~\emph{sentence} is a~formula with no free variables.

Let $\varphi$ be a~formula, $t$ be a~term, and $x$ be a~variable.
The \emph{substitution} of a~term~$t$
into the formula~$\varphi$
for the variable~$x$
is the formula $\varphi(t\mid x)$,
obtained by replacing every free occurrence of the variable~$x$
in~$\varphi$ by the term~$t$. The substitution $\varphi(t\mid x)$ is called
\emph{admissible} if for every variable~$x'$ occurring in the term~$t$
no free occurrence of~$x$ in~$\varphi$ is a~part of a~subformula
$\forall x'\losp \psi(x')$
or $\exists x'\losp \psi(x')$ of the formula~$\varphi$.

For example, in the case of the group language~$\mathcal L_G$
terms are variables.
If we have a~formula ${\forall x_1\losp (x_2=x_1)}$,
then the substitution of
the term~$x_1$ for~$x_2$ is not admissible, and the substitution
of the term~$x_3$ for~$x_2$ is admissible.

For the formula $\forall x_1\losp (x_2=x_1\cdot x_3)$
the substitution $x_1\mid x_2$
is not admissible, and the substitution $x_3\mid x_2$ is admissible.

Now let us introduce the following convention of notation: we use
$t(x_1,\dots,x_n)$ to denote
a~term~$t$ whose variables form a~subset of $\{ x_1,\dots,x_n\}$. Similarly,
we use $\varphi (x_1,\dots,x_n)$ to denote a~formula whose \emph{free}
variables form a~subset of $\{ x_1,\dots,x_n\}$.

We need
\emph{logical axioms} and \emph{rules of inference} to construct
a~formal system. Logical axioms are cited below.

Purely logical axioms.
\begin{enumerate}
\renewcommand{\labelenumi}{\theenumi.}
\item
$\varphi\Rightarrow (\psi\Rightarrow \varphi)$.
\item
$(\varphi\Rightarrow (\psi\Rightarrow \chi))\Rightarrow
((\varphi \Rightarrow \psi)\Rightarrow (\varphi\Rightarrow \chi))$.
\item
$(\neg \varphi\Rightarrow \neg \psi)\Rightarrow
((\neg \psi \Rightarrow \varphi)\Rightarrow \psi)$.
\item
$\forall x \varphi(x)\Rightarrow \varphi(t|x)$
if $t$ is a~term such that the substitution $t\mid x$ is admissible.
\item
$(\forall x (\psi\Rightarrow \varphi(x)))\Rightarrow
(\psi\Rightarrow (\forall x \varphi))$ if
$\psi$ does not contain any free occurrences of~$x$.%
\setcounter{lastla}{\value{enumi}}
\end{enumerate}

Identity axioms.
\begin{enumerate}
\renewcommand{\labelenumi}{\theenumi.}
\item
$x=x$.
\item
$y=z\Rightarrow t(x_1,\dots,x_{i-1},y,x_{i+1},\dots,x_n)=
t(x_1,\dots,x_{i-1},z,x_{i+1},\dots,x_n)$.
\item
$y=z\Rightarrow (\varphi(x_1,\dots,x_{i-1},y,x_{i+1},\dots,x_n)
\Leftrightarrow
\varphi (x_1,\dots,x_{i-1},z,x_{i+1},\dots,x_n))$,
where $x_1,\dots,x_n,y,z$ are variables, $t$~is a~term, and
$\varphi(x_1,\dots,x_n)$ is an elementary formula.%
\setcounter{lastar}{\value{enumi}}
\end{enumerate}

There are two inference rules.
\begin{enumerate}
\renewcommand{\labelenumi}{\theenumi.}
\item
The rule of detachment (modus ponens or MP):
from $\varphi$ and $\varphi\Rightarrow \psi$ infer $\psi$.
\item
The rule of generalization:
from $\varphi$ infer $\forall x\losp \varphi$.
\end{enumerate}

Let $\Sigma$ be a~collection of formulas and $\psi$ be a~formula of the
language~$\mathcal L$.
A~sequence $(\varphi_1,\dots, \varphi_n)$
of formulas of the language~$\mathcal L$ is called
a~\emph{deduction of the formula~$\psi$ from the collection~$\Sigma$}
if $\varphi_n=\psi$ and for any $1\le i\le n$
one of the following conditions is fulfilled:
\begin{enumerate}
\item
$\varphi_i$ belongs to~$\Sigma$ or is a~logical axiom;
\item
there exist $1\le k< j< i$ such that $\varphi_j$ is
$(\varphi_k \Rightarrow \varphi_i)$, i.e.,
$\varphi_i$~is obtained from $\varphi_k$ and
$\varphi_k\Rightarrow \varphi_i$ by the inference rule MP;
\item
there exists $1\le j< i$ such that $\varphi_i$ is
$\forall x\losp \varphi_j$,
where $x$ is not a~free variable of any formula from~$\Sigma$.
\end{enumerate}

Denote this deduction by
$ (\varphi_1,\dots,\varphi_n)\colon \Sigma\vdash \psi$.

If there exists a~deduction
$(\varphi_1,\dots,\varphi_n)\colon \Sigma\vdash \psi$,
then the formula~$\psi$
is called \emph{deducible in the language~$\mathcal L$
from the set~$ \Sigma$},
and the deduction $(\varphi_1,\dots,\varphi_n)$~is called
a~\emph{proof of~$\psi$}.

A~(\emph{first order}) \emph{theory}~$T$ in the language~$\mathcal L$
is some set of sentences of the language~$\mathcal L$.
A~\emph{set of axioms} of the theory~$T$ is any set of sentences which
has the same corollaries as~$T$.

\subsection{Theory of Classes and Sets NBG}\label{ss1.2}
The set theory of von Neumann, Bernays, and G\"odel NBG
(see~\cite{Mendelson}),
which will be a~base for all our constructions,
has one relation symbol~$P^2$, which denotes
a~2-place relation, no function, and no constant symbols.
We shall use Latin letters $X$, $Y$, and~$Z$
with subscripts as variables of this system.
We also introduce the abbreviations $X\in Y$ for $P(X,Y)$
and $X\notin Y$ for $\neg P(X,Y)$. The sign~${\in}$ can be interpreted
as the symbol of belonging.

The formula $X\subseteq Y$ is an abbreviation for the formula
$\forall Z\losp (Z\in X\Rightarrow Z\in Y)$
(\emph{inclusion}),
$X\subset Y$ is an abbreviation for
$X\subseteq Y\wedge X\ne Y$ (\emph{proper inclusion}).

Objects of the theory NBG are called \emph{classes}.
A class is called a~\emph{set} if it is an element of some class.
A class which is not a~set is called a~\emph{proper class}.
We introduce small Latin letters $x$, $y$, and~$z$ with subscripts
as special variables bounded by sets.
This means that the formula $\forall x\losp A(x)$ is an abbreviation
for $\forall X\losp (\text{$X$ is a~set}\Rightarrow A(X))$, and it has
the sense ``$A$~is true for all sets'',
and $\exists x\losp A(x)$ is an abbreviation for
$\exists X\losp (\text{$X$~is a~set} \logic\wedge A(X))$,
and it has the sense ``$A$~is true for some set.''

\begin{description}
\item[\hspace*{-\parindent}A1
{\mdseries (\emph{the extensionality axiom}).}]
$\forall X\losp \forall Y\losp
(X=Y\Leftrightarrow \forall Z\losp
(Z\in X\Leftrightarrow Z\in Y))$.
Intuitively, $X=Y$ if and only if $X$~and~$Y$ have the same elements.
\item[\hspace*{-\parindent}A2
{\mdseries (\emph{the pair axiom}).}]
$\forall x\losp \forall y\losp \exists z\losp \forall u\losp
(u\in z\Leftrightarrow u=x \logic\vee u=y)$,
i.e., for all sets $x$~and~$y$
there exists a~set~$z$ such that $x$~and~$y$ are the only
elements of~$z$.
\item[\hspace*{-\parindent}A3
{\mdseries (\emph{the empty set axiom}).}]
$\exists x\losp \forall y\losp \neg (y\in x)$,
i.e., there exists a~set which does not contain any elements.
\end{description}

Axioms \textbf{A1}~and~\textbf{A3} imply that this set is unique,
i.e., we can introduce a~constant symbol~$\varnothing$
(or~$0$), with the condition $\forall y\losp (y\notin \varnothing)$.

Also we can introduce a~new function symbol
$f(x,y)$ for the pair, and write it in the form
$\{ x,y\}$. Further, let $\{ x\}=\{ x,x\}$.
The set $\langle x,y\rangle \equiv\{ \{ x,\{ x,y\}\}$ is called
the \emph{ordered pair} of sets $x$~and~$y$.

\begin{proposition}\label{prop1}
$\vdash \forall x\losp\forall y\losp \forall u\losp \forall v\losp
(\langle x,y\rangle =\langle u,v\rangle \Rightarrow x=u\logic\wedge y=v)$.
\end{proposition}

In the same way we can introduce
\emph{ordered triplets of sets, ordered quadruplets of sets},
and so on.

\begin{description}
\item[\hspace*{-\parindent}AS4
{\mdseries (\emph{the axiom scheme of existence of classes}).}]
Let
$$
\varphi(X_1,\dots,X_n,Y_1,\dots,Y_m)
$$
be a~formula. We shall call this formula
\emph{predicative} if only variables for sets are bound in it
(i.e.,
if it can be transferred to this form with the help of abbreviations).
For every predicative formula $\varphi(X_1,\dots,X_n,Y_1,\dots,Y_m)$
$$
\exists Z\losp \forall x_1\dots \forall x_n\losp
(\langle x_1,\dots,x_n\rangle \in Z
\Leftrightarrow
\varphi(x_1,\dots,x_n,Y_1,\dots,Y_m)).
$$
\end{description}

The class~$Z$ which exists by the axiom scheme~\textbf{AS4}
will be denoted by
$$
\{ x_1,\dots,x_n\mid \varphi(x_1,\dots,x_n,Y_1,\dots,Y_m)\}.
$$

Now, by the axiom scheme \textbf{AS4}, we can define for arbitrary
classes $X$~and~$Y$ the following derivative classes:
\begin{itemize}
\item[]
$X\cap Y\equiv \{ u\mid u\in X \logic\land u\in Y\}$
(\emph{the intersection of classes $X$~and~$Y$});
\item[]
$X\cup Y\equiv \{ u\mid u\in X \logic\lor u\in Y\}$
(\emph{the union of classes $X$~and~$Y$});
\item[]
$\bar X\equiv \{ u\mid u\notin X\}$
(\emph{the complement of a~class~$X$});
\item[]
$V\equiv \{ u\mid u=u\}$ (\emph{the universal class});
\item[]
$X\setminus Y\equiv \{ u\mid u\in X \logic\land u\notin Y\}$
(\emph{the difference of classes $X$~and~$Y$});
\item[]
$\Dom(X)\equiv \{ u\mid \exists v\losp (\langle u,v\rangle \in X)\}$
(\emph{the domain of a~class~$X$});
\item[]
$\Rng(X)\equiv \{ u\mid \exists v (\langle v,u\rangle \in X)\}$
(\emph{the image of a~class~$X$});
\item[]
$X\times Y\equiv \{ u\mid \exists x\losp \exists y\losp
(u=\langle x,y\rangle \logic\land x\in X \logic\land y\in Y)\}$
(\emph{the Cartesian product of classes $X$~and~$Y$});
\item[]
$\mathcal P(X)\equiv \{ u\mid u\subseteq X\}$
(\emph{the class of all subsets of a~class~$X$});
\item[]
$\ucup X\equiv \{ u\mid \exists v\losp (u\in v \logic\land v\in X)\}$
(\emph{the union of all elements of a~class~$X$}).
\end{itemize}

Introduce now other axioms.

\begin{description}
\item[\hspace*{-\parindent}A5 {\mdseries (\emph{the union axiom}).}]
$\forall x\losp \exists y\losp \forall u\losp
(u\in y\Leftrightarrow \exists v\losp (u\in v \logic\land v\in x))$.
\end{description}

This axiom states that the union $\ucup x$ of all elements
of a~set~$x$ is also a~set.

\begin{description}
\item[\hspace*{-\parindent}A6 {\mdseries (\emph{the power set axiom}).}]
$\forall x\losp \exists y\losp \forall u\losp
(u\in y\Leftrightarrow u\subseteq x)$.
\end{description}

This axiom states that the class of all subsets of a~set~$x$
is a~set, which will be called the \emph{power set of~$x$}.

\begin{description}
\item[\hspace*{-\parindent}A7 {\mdseries (\emph{the separation axiom}).}]
$\forall x\losp \forall Y\losp \exists z\losp \forall u\losp
(u\in z\Leftrightarrow u\in x \logic\wedge u\in Y)$.
\end{description}

This axiom states that the intersection of a~class and a~set is a~set.

Denote the class $X\times X$ by $X^2$,
the class $X\times X\times X$ by~$X^3$,
and so on. Let the formula $\mathop{\mathrm{Rel}}(X)$ be an abbreviation
for the formula $X\subseteq V^2$ (\emph{$X$~is a~relation}),
$\Un(X)$ be an abbreviation
for the formula $\forall x\losp \forall y\losp \forall z\losp
(\langle x,y\rangle \in X \logic\wedge
\langle x,z\rangle \in X\Rightarrow y=z)$ (\emph{$X$~is functional}),
and
$\mathop{\mathrm{Fnc}}(X)$ be an abbreviation
for $X\subseteq V^2 \logic\wedge \Un(X)$
(\emph{$X$~is a~function}).

\begin{description}
\item[\hspace*{-\parindent}A8 {\mdseries (\emph{the replacement axiom}).}]
$\forall X\losp \forall x\losp
(\Un(X)\Rightarrow \exists y\losp \forall u\losp
(u\in y \Leftrightarrow \exists v\losp
(\langle v,u\rangle \in X \logic\wedge v\in x)))$.
\end{description}

This axiom states that if the class~$X$ is functional,
then the class of second components of pairs from~$X$ such that the first
component belongs to~$x$ is a~set.
The following axiom postulates existence of an infinite set.

\begin{description}
\item[\hspace*{-\parindent}A9 {\mdseries (\emph{the infinity axiom}).}]
$\exists x\losp (0\in x \logic\wedge
\forall u\losp (u\in x\Rightarrow u\cup \{ u\}\in x))$.
It is clear that for such a~set~$x$
we have
$\{ 0\}\in x$, $\{ 0,\{ 0\}\}\in x$,
$\{ 0,\{ 0\},\{ 0,\{ 0\}\}\}\in x$,\ldots{}
If we now set $1\prisv \{0\}$, $2\prisv \{ 0,1\}$,\ldots,
$n\prisv \{ 0,1,\dots,n-1\}$,
then for every integer $n\ge 0$ the condition $n\in x$ is fulfilled
and $0\ne 1$, $0\ne 2$, $1\ne 2$,\ldots
\item[\hspace*{-\parindent}A10 {\mdseries (\emph{the regularity axiom}).}]
$\forall X\losp (X\ne \varnothing\Rightarrow
\exists x\in X\losp (x\cap X = \varnothing))$.
\end{description}

This axiom states that every nonempty set is disjoint
from one of its elements.

\begin{description}
\item[\hspace*{-\parindent}A11 {\mdseries (\emph{the axiom of choice AC}).}]
For every set~$x$ there exists
a~mapping~$f$ such that for every nonempty subset
$y\subseteq x$ we have $f(y)\in y$
(this mapping is called a~\emph{choice} mapping for~$x$).
\end{description}

The list of axioms of the theory NBG is finished.

A~class~$P$ is called \emph{ordered by a~binary relation~$\le$ on~$P$},
if the following conditions hold
\begin{enumerate}
\item
$\forall p\in P\losp (p\le p) $;
\item
$\forall p,q\in P\losp (p\le q \logic\land q\le p\Rightarrow p=q)$;
\item
$\forall p,q,r\in P\losp (p\le q \logic\land q\le r\Rightarrow p\le r)$.
\end{enumerate}

If, in addition,\nopagebreak
\begin{enumerate}
\setcounter{enumi}{3}
\item
$\forall p,q\in P\losp (p\le q \logic\lor q\le p)$,
\end{enumerate}
then the relation~$\le$ is called a~\emph{linear order}
on the class~$P$.

An ordered class~$P$ is called \emph{well-ordered} if
\begin{enumerate}
\setcounter{enumi}{4}
\item
$\forall q\losp (\varnothing \ne q\subseteq P\Rightarrow
\exists x\in q\losp (\forall y\in q\losp (x\le y)))$,
i.e., every nonempty subset of the class~$P$
has the smallest element.
\end{enumerate}

A~class~$S$ is called \emph{transitive} if
$\forall x\losp (x\in S\Rightarrow x\subseteq S)$.

A~class (a~set)~$S$ is called an \emph{ordinal}
(an \emph{ordinal number}) if $S$
is transitive and well-ordered by the relation ${\in}\cup {=}$ on~$S$.

Ordinal numbers are usually denoted by Greek letters
$\alpha$, $\beta$, $\gamma$, and so on. The class of all
ordinal numbers is denoted by~$\On$.
The natural ordering of the class of ordinal numbers is the relation
${\alpha\le \beta} \prisv \allowbreak
\alpha=\beta \logic\lor \alpha\in \beta$.
The class $\On$ is transitive and linearly ordered by the relation~$\le$.

There are some simple assertions about ordinal numbers:
\begin{enumerate}
\item
if $\alpha$ is an ordinal number, $a$ is a~set, and $a\in \alpha$,
then $a$ is an ordinal number;
\item
$\alpha+1\equiv \alpha\cup \{ \alpha\}$
is the smallest ordinal number that is greater
than~$\alpha$;
\item
every nonempty set of ordinal numbers has the smallest element.
\end{enumerate}

Therefore the ordered class $\On$ is well-ordered.
Thus $\On$ is an ordinal.

An ordinal number $\alpha$ is called a~\emph{successor} if
$\alpha=\beta+1$ for some ordinal number~$\beta$.
In the opposite case $\alpha$ is called a~\emph{limit ordinal number}.

The smallest (in the class $\On$)
nonzero limit ordinal is denoted by~$\omega$.
Ordinals which are smaller than~$\omega$
are called \emph{natural numbers}.

Classes~$F$ which are functions with domains equal to~$\omega$
are called \emph{infinite sequences}.
Functions with domains
equal to $n\in \omega$ are called \emph{finite sequences}.

Sets $a$~and~$b$ are called
\emph{equivalent} ($a \sim b$) if there exists a~bijective function
$u\colon a\to b$.

An ordinal number~$\alpha$ is called a~\emph{cardinal} if for every
ordinal number~$\beta$ the conditions
$\beta \leq \alpha$ and $\beta\sim\alpha$
imply $\beta =\alpha$. The class of all cardinal numbers
is denoted by~$\Cn$. The class~$\Cn$
with the order induced from the class~$\On$ is well-ordered.

The axiom of choice implies that
for every set~$a$ there exists a~unique cardinal number~$\alpha$
such that $a\sim \alpha$.
This number~$\alpha$ is called the \emph{power of the set}~$a$
(denoted by $|a|$, or $\card a$).
A~set of power~$\omega$ is called \emph{countable}.
A~set of power $n\in \omega$ is called \emph{finite}.
A~set is called \emph{infinite} if it is not finite.
A~set is called \emph{uncountable} if it is neither countable, nor finite.
The cardinal number $c\prisv |\mathcal P(\omega)|$ is called the power
of \emph{continuum}.

To denote cardinals we shall use
small Greek letters (as in the case of ordinals):
$\xi$th infinite cardinal will be denoted
by~$\omega_\xi$
(i.e., the cardinal number~$\omega$ will also be denoted by~$\omega_0$).

A~set~$X$ is said to be \emph{cofinal} in~$\alpha$
if $X\subset \alpha$ and $\alpha=\ucup X$.
The cofinality of~$\alpha$, written $\cf \alpha$,
is the least cardinal~$\beta$
such that a~set of power~$\beta$ is cofinal in~$\alpha$.

A~cardinal~$\varkappa$ is said to be \emph{regular} if
$\cf \varkappa=\varkappa$, i.e., for every ordinal number~$\beta$
for which there exists a~function $f\colon \beta\to \varkappa$ such that
$\ucup \rng  f=\varkappa$ the inequality $\varkappa \le \beta$ holds,
where $\ucup \rng f=\varkappa$ means that for every $y\in \varkappa$
there exists $x\in \beta$ such that $y< f(x)$.
A~cardinal~$\varkappa$ is said
to be \emph{singular} if it is not regular.

The continuum hypothesis states that $|\mathcal P(\omega)|=\omega_1$,
i.e., the power of continuum is the smallest uncountable cardinal
number.

We shall assume the continuum hypothesis if we need it.

\subsection{Models, Satisfaction, and Elementary Equivalence}\label{ss1.3}
We now suppose that all our constructions are made in the theory NBG.

A~\emph{model of a~first order language~$\mathcal L$}
is a~pair~$\mathcal U=\langle A,I\rangle$
consisting of a~universe~$A$ (i.e., some class or set of the theory NBG)
and some correspondence~$I$ that
assigns to every relation symbol~$Q^n$ some
$n$-place relation $R\subset A^n$ on~$A$,
to every function symbol~$F^n$ some
$m$-place function $G\colon A^m\to A$, and to every constant symbol~$c$
some element of~$A$.

A simple example of a~model of the group language is the set
$1\prisv \{\varnothing\}=\{ 0\}$, where $I(Q^2)=\{ \langle 0,0,0\rangle\}$.
Another simple example of a~model of the group language is
the set $2\prisv\{ \varnothing,\{ \varnothing\}\}=\{ 0,1\}$,
where
$I(Q^2)=\{ \langle 0,0,0\rangle, \langle 0,1,1\rangle,
\langle 1,0,1\rangle, \langle 1,1,0\rangle\}$.

The \emph{power} of a~model $\mathcal U=\langle A,I\rangle$
is the cardinal number~$|A|$ (if the universe~$A$ is a~set).
For all models which will be considered in this paper,
the universe~$A$ is a~set.
A~model~$\mathcal U$ is called finite, countable, or uncountable
if $|A|$ is finite, countable, or uncountable, respectively.

Models $\mathcal U$~and~$\mathcal U'$ of a~language~$\mathcal L$
are called \emph{isomorphic},
if there exists a~bijective mapping~$f$ of the set (universe)~$A$
onto the set~$A'$
satisfying the following conditions:
\begin{enumerate}
\item
for each $n$-place relation symbol~$Q^n$ and any $a_1,\dots,a_n$
from~$A$
$$
\langle a_1,\dots,a_n\rangle\in I(Q^n)\text{ if and only if }
\langle f(a_1),\dots,f(a_n)\rangle \in I'(Q^n);
$$
\item
for each $m$-place function symbol~$F^m$
of the language~$\mathcal L$ and any $a_1,\dots,a_m \in A$
$$
f(I(F^m)(\langle a_1,\dots,a_m\rangle))=
I'(F^m)(\langle f(a_1),\dots,f(a_m)\rangle);
$$
\item
for each constant symbol~$c$ of the language~$\mathcal L$
$$
f(I(c))=I'(c).
$$
\end{enumerate}

Every mapping~$f$ satisfying these conditions is called
an \emph{isomorphism of the model~$\mathcal U$ onto the model~$\mathcal U'$}
or an \emph{isomorphism between the models $\mathcal U$~and~$\mathcal U'$}.
The fact that $f$ is an
isomorphism of the model~$\mathcal U$ onto the model~$\mathcal U'$
will be denoted
by $f\colon \mathcal U\cong \mathcal U'$, and the formula
$\mathcal U\cong \mathcal U'$ means that
the models $\mathcal U$~and~$\mathcal U'$ are isomorphic.
For convenience we use~$\cong$ to denote the isomorphism relation
between models.

Indeed, unless we wish to consider the particular structure of each element
of $A$~or~$A'$, for all practical purposes $\mathcal U$~and~$\mathcal U'$
are the same if they are isomorphic.

Now we shall give a~formal definition of satisfiability.
Let $\varphi$ be an arbitrary
formula of a~language~$\mathcal L$,
let all its variables, free and bound,
be contained in the set $x_1,\dots,x_q$, and let
$a_1,\dots,a_q$ be an arbitrary sequence of elements
of the set~$A$. We define the predicate
$$
\text{\emph{$\varphi$ is true on the sequence $a_1,\dots,a_q$
in the model~$\mathcal U$}, or \emph{$a_1,\dots,a_q$ satisfies the
formula~$\varphi$ in~$\mathcal U$}.}
$$

The definition proceeds in three stages.
Let $\mathcal U$ be a~fixed model for~$\mathcal L$.

1. The value of a~term $t(x_1,\dots,x_q)$ at
$a_1,\dots,a_q$ is defined as follows
(we let $t[a_1,\dots,a_q]$ denote this value):
\begin{enumerate}
\item
if $t$ is a~variable $x_i$, then $t[a_1,\dots,a_q]=a_i$;
\item
if $t$ is a~constant symbol~$c$, then $t[a_1,\dots,a_q]=I(c)$;
\item
if $t$ is $F^m(t_1,\dots,t_m)$, where
$t_1(x_1,\dots,x_q),\dots, \allowbreak t_m(x_1,\dots,x_q)$ are terms, then
$$
t[a_1,\dots,a_q]=
I(F^m)(\langle t_1[a_1,\dots,a_q],\dots,t_m[a_1,\dots,a_q]\rangle).
$$
\end{enumerate}

2.
\begin{enumerate}
\item
Suppose that $\varphi(x_1,\dots\ttm,x_q)$ is an
elementary formula $t_1\!=\!t_2$, where $t_1(x_1,\dots\ttm,x_q)$ and
$t_2(x_1,\dots,x_q)$ are terms. Then $a_1,\dots,a_q$ satisfies $\varphi$
if and only if
$$
t_1[a_1,\dots,a_q]=t_2[a_1,\dots,a_q].
$$
\item
Suppose that $\varphi(x_1,\dots,x_q)$ is an elementary formula
$Q^n(t_1,\dots,t_n)$, where
$Q^n$ is an $n$-place relation symbol and
$t_1(x_1,\dots,x_q),\dots, t_n(x_1,\dots,x_q)$ are terms.
Then $a_1,\dots,a_q$ satisfies~$\varphi$ if and only if
$$
\langle t_1[a_1,\dots, a_q],\dots , t_n[a_1,\dots,a_q]\rangle \in I(Q^n).
$$
\end{enumerate}
For brevity, we write
$$
\mathcal U\vDash \varphi [a_1,\dots, a_q]
$$
for: $a_1,\dots,a_q$ satisfies~$\varphi$ in~$\mathcal U$.

3. Now suppose that $\varphi$ is any formula of~$\mathcal L$
and all its
free and bound variables are among $x_1,\dots,x_q$.
\begin{enumerate}
\item
If $\varphi$ is $\theta_1\land \theta_2$, then
$$
\mathcal U\vDash \varphi[a_1,\dots,a_q]\text{ if and only if }
\mathcal U\vDash \Theta_1 [a_1,\dots,a_q]\text{ and }\mathcal U\vDash
\Theta_2 [a_1,\dots, a_q].
$$
\item
If $\varphi$ is $\neg \Theta$, then
$$
\mathcal U\vDash \varphi[a_1,\dots,a_q]
\text{ if and only if it is not true that }
\mathcal U\vDash \Theta [a_1,\dots,a_q].
$$
\item
If $\varphi$ is $\forall x_i\losp \psi$, where $i\le q$, then
$$
\mathcal U\vDash \varphi[a_1,\dots,a_q]\text{ if and only if }
\mathcal U\vDash \psi [a_1,\dots,a_{i-1},a, a_{i+1},\dots,a_q]
\text{ for any } a\in A.
$$
\end{enumerate}

It is easy to check that the abbreviations
$\lor$, $\Rightarrow$, $\Leftrightarrow$, and
$\exists$ have their usual meanings. In particular,
if $\varphi$ is $\exists x_i\losp \psi$, where $i\le q$, then
$\mathcal U\vDash \varphi[a_1,\dots,a_q]$
if and only if there exists $a\in A$ such that
$$
\mathcal U\vDash \psi[a_1,\dots,a_{i-1},a,a_{i+1},\dots,a_q].
$$

The following proposition shows that
the relation
$$
\mathcal U\models \varphi(x_1,\dots,x_p)[a_1,\dots,a_q]
$$
depends only on $a_1,\dots,a_p$, where $p< q$.

\begin{proposition}\label{prop2}
\begin{enumerate}
\item
Let $t(x_1,\dots,x_p)$ be a~term, and let $a_1,\dots,a_q$ and
$b_1,\dots,b_r$ be two sequences of elements such that
$p\le q$, $p\le r$, and $a_i=b_i$ whenever $x_i$ is a~free
variable of the term~$t$. Then
$$
t[a_1,\dots,a_q]=t[b_0,\dots,b_r].
$$
\item
Let $\varphi$ be a~formula, let all its variables,
free and bound, belong to the set $x_1,\dots,x_p$, and let
$a_1,\dots,a_q$ and $b_1,\dots,b_r$
be two sequences of elements such that
$p\le q$, $p\le r$, and $a_i=b_i$ whenever $x_i$ is a~free variable
in the formula~$\varphi$. Then
$$
\mathcal U\models \varphi [a_1,\dots,a_q]\quad
\text{if and only if}\quad
\mathcal U\models \varphi [b_1,\dots,b_r].
$$
\end{enumerate}
\end{proposition}

This proposition allows us to give the following definition.
Let $\varphi(x_1,\dots,x_p)$ be a~formula, and let all its variables,
free and bound, be contained in the set $x_1,\dots,x_q$, where $p\le q$.
Let $a_1,\dots, a_p$ be a~sequence of elements of the set~$A$.
We shall say that \emph{$\varphi$~is true in~$\mathcal U$ on
$a_1,\dots,a_p$},
$$
\mathcal U\models \varphi[a_1,\dots,a_p]
$$
if $\varphi$ is true in~$\mathcal U$ on $a_1,\dots,a_p,\dots,a_q$
with some (or, equivalently, any) sequence $a_{p+1},\dots,a_q$.

Let $\varphi$ be a~sentence, and let all its bound variables
be contained in the set $x_1,\dots,x_q$. We shall say that
\emph{$\varphi$~is true in the model~$\mathcal U$}
(notation: $\mathcal U\models \varphi$) if $\varphi$~is true
in~$\mathcal U$ on some (equivalently, any) sequence
$a_1,\dots,a_q$.

This last phrase is equivalent to each of the following phrases:
\begin{align*}
& \text{$\varphi$ holds in $\mathcal U$;}\\
& \text{$\mathcal U$ satisfies $\varphi$;}\\
& \text{$\mathcal U$ is a~model of $\varphi$.}
\end{align*}

In the case where $\sigma$ is not true in~$\mathcal U$,
we say that $\sigma$ is \emph{false} in~$\mathcal U$,
or that $\sigma$ \emph{does not hold} in~$\mathcal U$,
or that $\mathcal U$~\emph{is a~model of the sentence} $\neg \sigma$.
If we have a~set~$\Sigma$
of sentences, we say that $\mathcal U$ is a~model of this set
if $\mathcal U$ is a~model of every sentence $\sigma\in \Sigma$.
It is useful to denote this concept by $\mathcal U\models \Sigma$.

As we have said above, a~\emph{theory}~$T$ of the language~$\mathcal L$ is
a~collection of sentences of the language~$\mathcal L$.
A~\emph{theory} of a~model~$\mathcal U$
(of the language~$\mathcal L$)
is the set of all sentences which hold in~$\mathcal U$.

Two models $\mathcal U$~and~$\mathcal V$ for~$\mathcal L$ are called
\emph{elementarily equivalent} if every sentence that is true
in~$\mathcal U$
is true in~$\mathcal V$, and vice versa.
We express this relationship between models
by~$\equiv$. It is easy to see that $\equiv$ is
indeed an equivalence relation.
We can note that two models are elementarily equivalent if and only if
their theories coincide.

Any two isomorphic models of the same language are elementarily equivalent.
If two models of the same language are elementarily equivalent
and one of them is finite, then these models are isomorphic.
If models are infinite
and elementarily equivalent, they are not necessarily isomorphic.
For example, the field~$\mathbb C$ of all complex numbers and the field
$\bar{\mathbb Q}$ of all algebraic numbers are elementarily equivalent,
but not isomorphic.

Together with first order languages we need to consider
second order languages, in which we can also quantify
relation symbols, i.e., use relation symbols as variables.
Such languages will be described in Sec.~\ref{ss1.4}.

\subsection{Second Order Languages and Models}\label{ss1.4}
Now we shall introduce all notions, similar to
the notions of Secs.\ \ref{ss1.1}~and~\ref{ss1.3},
for second order languages and models.

A~\emph{second order language}~$\mathcal L_2$ is a~collection of symbols,
consisting of
(1)~parentheses ${(}$,~${)}$;
(2)~connectives $\land$~(``and'') and $\neg$~(``not'');
(3)~the quantifier $\forall$~(for all);
(4)~the binary relation symbol ${=}$~(identity);
(5)~a~countable set of object variables~$x_i$;
(6)~a~countable set of predicate variables~$P_i^l$;
(7)~a~finite or countable set of relation symbols~$Q_i^n$ ($n\ge 1$);
(8)~a~finite or countable set of function symbols~$F_i^n$ ($n\ge 1$);
(9)~a~finite or countable set of constant symbols~$c_i$.

\emph{Terms} of the language $\mathcal L_2$ are defined as follows:
\begin{enumerate}
\item
a~variable is a~term;
\item
a~constant symbol is a~term;
\item
if $F^n$ is an $n$-place function symbol and
$t_1,\dots,t_{n}$ are terms, then $F^n(t_1,\dots,t_{n})$ is a~term;
\item
a~symbol-string is a~term only if it can be shown to be a~term
by a~finite number of applications of (1)--(3).
\end{enumerate}

Therefore terms of the language~$\mathcal L_2$ coincide with terms
of the language~$\mathcal L$.

\emph{Elementary formulas} of the language~$\mathcal L_2$ are
symbol-strings of the form given below:
\begin{enumerate}
\item
if $t_1$~and~$t_2$ are terms of the language~$\mathcal L_2$,
then $t_1=t_2$ is an elementary formula;
\item
if $P^l$ is a~predicate variable and $t_1,\dots,t_l$ are terms,
then the symbol-string $P^l(t_1,\dots,t_l)$ is an
elementary formula;
\item
if $Q^n$ is an $n$-place relation symbol,
and $t_1,\dots,t_{n}$ are terms, then the symbol-string
$Q^n(t_1,\dots,t_{n})$ is an elementary formula.
\end{enumerate}

Therefore elementary formulas of the second order group language have
the form $x_i=x_j$, $x_i=x_j\cdot x_k$, and
$P^l(x_{i_1},\dots,x_{i_l})$, where $l\ge 1$.

\emph{Formulas} of the language~$\mathcal L_2$ are defined as follows:
\begin{enumerate}
\item
an elementary formula is a~formula;
\item
if $\varphi$ and $\psi$ are formulas and $x$~is an object variable,
then
$(\neg \varphi)$, $(\varphi\land \psi)$, and
$(\forall x \losp \varphi)$ are formulas;
\item
if $P^l$ is a~predicate variable and $\varphi$ is a~formula, then
the symbol-string
$(\forall P^l(v_1,\dots,v_l)\varphi)$ is a~formula;
\item
a~symbol-string is a~formula only
if it can be shown to be a~formula by a~finite number
of applications of (1)--(3).
\end{enumerate}

Let us introduce the following abbreviations:
\begin{itemize}
\item[]
$\exists P^l(v_1,\dots,v_l)\losp \varphi$ is an abbreviation for
$\neg (\forall P^l(v_1,\dots,v_l)\losp (\neg \varphi))$;
\item[]
$\forall P^{l_1}_1(v_1,\dots,v_{l_1})\dots
\forall P_n^{l_n}(v_1,\dots,v_{l_n})\losp \varphi$ is an abbreviation for
$$
(\forall P^{l_1}_1(v_1,\dots,v_{l_1}))\dots
(\forall P_n^{l_n}(v_1,\dots,v_{l_n}))\losp \varphi;
$$
\item[]
$\exists P^{l_1}_1(v_1,\dots,v_{l_1})\dots
\exists P_n^{l_n}(v_1,\dots,v_{l_n})\losp \varphi$ is an abbreviation for
$$
(\exists P^{l_1}_1(v_1,\dots,v_{l_1}))\dots
(\exists P_n^{l_n}(v_1,\dots,v_{l_n}))\losp \varphi.
$$
\end{itemize}

Introduce the notions of \emph{free} and \emph{bound} occurrence
of a~predicate variable in a~formula of the language~$\mathcal L_2$.

\begin{enumerate}
\item
All occurrences of all predicate variables in elementary formulas
are free occurrences.
\item
Every free (bound) occurrence of a~variable~$P^l$ in a~formula~$\varphi$
is a~free (bound) occurrence of a~variable~$P^l$ in the formulas
$(\neg \varphi)$, $(\varphi\land \psi)$, and $(\psi \land \varphi)$.
\item
For any occurrence of a~variable~$P^l$ in a~formula~$\varphi$,
the occurrence of the variable~$P^l$ in the formula
$\forall P^l(v_1,\dots,v_l)\losp \varphi$ is bound.
If an occurrence of a~variable~$P_1^l$ in a~formula~$\varphi$
is free (bound), then the occurrences of~$P_1^l$
in the formulas $\forall x\losp \varphi$
and $\forall P_2^m (v_1,\dots,v_m)\losp \varphi$ are free (bound).
\end{enumerate}

As in Sec.~\ref{ss1.1}, any formula such that all its \emph{free} object
and predicate variables are among the set
$\{ x_1,\dots,x_n, P_1^{l_1},\dots,P_k^{l_k}\}$
will be denoted by $\varphi(x_1,\dots,x_n,P_1^{l_1},\dots,P_k^{l_k})$.

To our five purely logical axioms from Sec.~\ref{ss1.1} we shall
add the sixth purely logical axiom:
\begin{enumerate}
\setcounter{enumi}{\value{lastla}}
\renewcommand{\labelenumi}{\theenumi.}
\item
$(\forall P^n(v_1,\dots,v_n)\losp (\psi\Rightarrow \varphi)\Rightarrow
(\psi\Rightarrow (\forall P^n(v_1,\dots,v_n)\losp \varphi))$ if
$\psi$ does not contain any free occurrences of the variable~$P^n$.
\end{enumerate}

To the identity axioms we add the fourth identity axiom:
\begin{enumerate}
\setcounter{enumi}{\value{lastar}}
\renewcommand{\labelenumi}{\theenumi.}
\item
$\forall P^n(v_1,\dots,v_n)\losp (y=z\Rightarrow
(P^n(x_1,\dots,x_{i-1},y,x_{i+1},\dots,x_n)\Leftrightarrow
P^n (x_1, \dots, x_{i-1},z,x_{i+1},\dots,x_n))$.
\end{enumerate}

The rule of generalization can be changed
to ``from~$\varphi$ infer $\forall x\losp \varphi$ and
$\forall P^n(v_1,\dots,v_n)\losp \varphi$.''

A \emph{model of a~second order language~$\mathcal L_2$} (see
Sec.~\ref{ss1.3}) is a~pair~$\mathcal U=\langle A,I\rangle$
consisting of an object~$A$ (i.e., some class or set of the theory NBG)
and some correspondence~$I$
that assigns to every relation symbol~$Q^n$
some $n$-place relation in~$A$, to every function
symbol~$F^n$ some $n$-place function in~$A$,
and to every constant symbol~$c$ some element of~$A$.

Now we shall give a~definition of satisfaction.
Let $\varphi$ be any formula
of the language~$\mathcal L_2$ such that all its free and bound variables
are among $x_1,\dots,x_q,P_1^{l_1},\dots,P_s^{l_s}$,
and let $a_1,\dots,a_q, b_1^{l_1},\dots,b_s^{l_s}$ be any
sequence, where $a_1,\dots,a_q$ are elements of the set~$A$,
$b_i^{l_i} \subset A^{l_i}$.
We define the predicate
$$
\text{\emph{$\varphi$ is satisfied by the sequence
$a_1,\dots,a_q,b_1^{l_1},\dots,b_s^{l_s}$ in the model~$\mathcal U$}}.
$$

1.
The value of a~term $t(x_1,\dots,x_q,P_1^{l_1},\dots,P_s^{l_s})$
at $a_1,\dots,a_q,\allowbreak b_1^{l_1},\dots,b_s^{l_s}$
is defined as follows (we let
$t[a_1,\dots,a_q,b_1^{l_1},\dots,b_s^{l_s}]$):
\begin{enumerate}
\item
if $t$ is a~variable $x_i$, then
$t[a_1,\dots,a_q,b_1^{l_1},\dots,b_s^{l_s}]=a_i$;
\item
if $t$ is a~constant symbol~$c$,
then $t[a_1,\dots,a_q,b_1^{l_1},\dots,b_s^{l_s}]=I(c)$;
\item
if $t$ is $F^m(t_1,\dots,t_m)$, where
$t_1(x_1,\dots,x_q,P_1^{l_1},\dots,P_s^{l_s}),\dots,\allowbreak
t_m(x_1,\dots,x_q, P_1^{l_1},\dots,P_s^{l_s})$ are terms, then
$t[a_1,\dots,a_q,b_1^{l_1},\dots,b_s^{l_s}]=
I(F^m)(\langle t_1[a_1,\dots,a_q,b_1^{l_1},\dots,b_s^{l_s}],\dots,
t_m[a_1,\dots,a_q,b_1^{l_1},\dots,b_s^{l_s}]\rangle)$.
\end{enumerate}

2.
\begin{enumerate}
\item
Suppose that $\varphi(x_1,\dots,x_q,P_1^{l_1},\dots,P_s^{l_s})$
is an elementary formula $t_1=t_2$,
where $t_1(x_1,\dots,x_q,\allowbreak P_1^{l_1},\dots,P_s^{l_s})$ and
$t_2(x_1,\dots,x_q,P_1^{l_1},\dots,P_s^{l_s})$ are terms. Then
$a_1,\dots,a_q,b_1^{l_1},\dots,b_s^{l_s}$ satisfies $\varphi$ if and only if
$$
t_1[a_1,\dots,a_q,b_1^{l_1},\dots,b_s^{l_s}]=
t_2[a_1,\dots,a_q,b_1^{l_1},\dots,b_s^{l_s}].
$$
\item
Suppose that $\varphi(x_1,\dots,x_q,P_1^{l_1},\dots,P_s^{l_s})$ is
an elementary formula $Q^n(t_1,\dots,t_n)$, where
$Q^n$ is an $n$-place relation symbol and
$t_1(x_1,\dots,x_q,P_1^{l_1},\dots,P_s^{l_s}),\dots,\allowbreak
t_n(x_1,\dots,x_q,P_1^{l_1},\dots,P_s^{l_s})$ are terms.
Then $a_1,\dots,a_q,b_1^{l_1},\dots,b_s^{l_s}$ satisfies $\varphi$
if and only if
$$
\langle t_1[a_1,\dots, a_q,b_1^{l_1},\dots,b_s^{l_s}],\dots,
t_n[a_1,\dots,a_q,b_1^{l_1},\dots,b_s^{l_s}]\rangle \in I(Q^n).
$$
\item\sloppy
Suppose that
$\varphi(x_1,\dots,x_q,\allowbreak P_1^{l_1},\dots,P_s^{l_s})$ is an
elementary formula $P_i^{l_i}(t_1,\dots,t_{l_i})$, where
$t_1(x_1,\dots,x_q,\allowbreak P_1^{l_1},\dots,P_s^{l_s}),\dots,\allowbreak
t_n(x_1,\dots,x_q,P_1^{l_1},\dots,P_s^{l_s})$ are terms.
Then $a_1,\dots,a_q,b_1^{l_1},\dots,b_s^{l_s}$ satisfies $\varphi$
if and only if
$$
\langle t_1[a_1,\dots, a_q,b_1^{l_1},\dots,b_s^{l_s}],\dots,
t_n[a_1,\dots,a_q,b_1^{l_1},\dots,b_s^{l_s}]\rangle \in b_i^{l_i}.
$$
\end{enumerate}

3. Now suppose that $\varphi$ is any formula
such that all its free and bound
variables are among $x_1,\dots,x_q,\allowbreak P_1^{l_1},\dots,P_s^{l_s}$.
\begin{enumerate}
\item
If $\varphi$ is $\theta_1\land \theta_2$, then
\begin{multline*}
\mathcal U\vDash \varphi[a_1,\dots,a_q,b_1^{l_1},\dots,b_s^{l_s}]
\text{ if and only if}\\[-\jot]
\mathcal U\vDash \Theta_1 [a_1,\dots,a_q,b_1^{l_1},\dots,b_s^{l_s}]
\text{ and }\mathcal U\vDash
\Theta_2 [a_1,\dots, a_q,b_1^{l_1},\dots,b_s^{l_s}].
\end{multline*}
\item
If $\varphi$ is $\neg \Theta$, then
$$
\mathcal U\vDash \varphi[a_1,\dots,a_q,b_1^{l_1},\dots,b_s^{l_s}]
\text{ if and only if it is false that }
\mathcal U\vDash \Theta [a_1,\dots,a_q,b_1^{l_1},\dots,b_s^{l_s}].
$$
\item
If $\varphi$ is $\forall x_i \psi$, where $i\le q$, then
\begin{multline*}
\mathcal U\vDash \varphi[a_1,\dots,a_q,b_1^{l_1},\dots,b_s^{l_s}]
\text{ if and only if}\\[-\jot]
\mathcal U\vDash
\psi [a_1,\dots,a_{i-1},a, a_{i+1},\dots,a_q,b_1^{l_1},\dots,b_s^{l_s}]
\text{ for any } a\in A.
\end{multline*}
\item
If $\varphi$ is $\forall P_i^{l_i}(v_1,\dots v_{l_i}) \psi$,
where $i\le s$, then
\begin{multline*}
\mathcal U\vDash \varphi[a_1,\dots,a_q,b_1^{l_1},\dots,b_s^{l_s}]
\text{ if and only if}\\[-\jot]
\mathcal U\vDash
\psi [a_1,\dots,a_q,b_1^{l_1},\dots,b_{i-1}^{l_{i-1}},
b,b_{i+1}^{l_{i+1}},\dots,b_s^{l_s}]
\text{ for any } b\subset A^{l_i}.
\end{multline*}
\end{enumerate}

The proposition that
$$
\mathcal U\models
\varphi(x_1,\dots,x_p,P_1^{l_1},\dots,P_t^{l_t})
[a_1,\dots,a_q,b_1^{l_1},\dots,b_s^{l_s}]
$$
depends only on the values $a_1,\dots,a_p,b_1^{l_1},\dots,b_t^{l_t}$,
where $p< q$, $s<t$, is the same as Proposition~\ref{prop2}.

All other definitions are also similar to definitions from Sec.~\ref{ss1.3}.

We shall say that two models of the second order language~$\mathcal L_2$
are equivalent in~$\mathcal L_2$ if
for every sentence of this language
the sentence is true in one model
if and only if it is true in the other model.

\section{Basic Concepts about Abelian Groups}
\subsection{Preliminaries}\label{ss2.1}
The word ``group'' will mean, throughout, an additively written
Abelian (i.e., commutative) group. That is, by group
we shall understand a~set~$A$
such that with every pair of elements $a,b\in A$
there is associated an element $a+b$ of~$A$, which is called the \emph{sum}
of elements $a$~and~$b$; there is an element $0\in A$, the \emph{zero},
such that $a+0=a$ for every $a\in A$; for each $a\in A$ there is an $x\in A$
with the property $a+x=0$,
this $x=-a$ is called the \emph{inverse} (\emph{opposite})
to~$a$; finally, we have both commutative and associative laws:
$a+b=b+a$, $(a+b)+c=a+(b+c)$ for every $a,b,c\in A$.

A~sum $a+\dots+a$ ($n$~times) is abbreviated as $na$,
and $-a-\dots-a$ ($n$~times) as~$-na$. By the \emph{order} of a~group~$A$
we mean the power~$|A|$ of the set of its elements.
If the power~$|A|$ is a~finite (countable) cardinal, then the group~$A$
is called \emph{finite} (\emph{countable}).

A~subset~$B$ of~$A$ is a~\emph{subgroup} if
$\forall b_1,b_2\in B\losp (b_1+b_2\in B)$.
If $B$ is a~subgroup consisting of the zero alone
or of all elements of~$A$, then $B$ is a~\emph{trivial}
subgroup of~$A$; but a~subgroup of~$ A $ that is different from~$A$
is called a~\emph{proper}
subgroup of~$A$. We shall write $B\lhd A$ to indicate that $B$ is
a~subgroup of~$A$. Let $B\lhd A$ and $a\in A$.
The set $a+B=\{ a+b\mid b\in B\}$ is said to be a~\emph{coset
of~$A$ modulo~$B$}. An element of the coset is called a~\emph{representative}
of this coset. A~set consisting of representatives,
one for each coset of~$A$ modulo~$B$,
is called the \emph{complete system of representatives of cosets modulo~$B$}.
Its power is called the index of~$B$ in~$A$, and denoted by $|A:B|$.

The cosets of~$A$ modulo~$B$ form a~group $A/B$
known as the \emph{quotient group} (of~$A$ with respect to~$B$). If $S$
is any subset in~$A$, then by $\langle S\rangle$
we denote the subgroup of~$A$,
\emph{generated} by~$S$, i.e., the intersection of all subgroups of~$A$
containing~$S$. In particular,
if $S$ consists of the elements $a_i$ ($i\in I$),
we also write
$\langle S\rangle=\langle \dots, a_i,\dots \rangle_{i\in I}$ or simply
$\langle S\rangle =\langle a_i\rangle_{i\in I}$.
The subgroup $\langle S\rangle$ consists exactly
of all finite linear combinations of the elements of~$S$, i.e., of all sums
$n_1a_1+\dots+n_ka_k$ with $a_i\in S$, $n_i$ integers, and $k$ an arbitrary
nonnegative integer.
If $S$ is empty, then we put $\langle S\rangle=0$.
In the case $\langle S\rangle=A$
we say $S$ is a~\emph{generating system} of~$A$ and the elements of~$S$
are \emph{generators} of~$A$.
If there is a~finite generating system, then $A$~is said to be
a~\emph{finitely generated} group.

If $B$ and $C$ are two subgroups of~$A$,
then the subgroup $\langle B,C\rangle$
generated by them consists of all elements of~$A$ of the form $b+c$,
where $b\in B$ and $c\in C$. Therefore we shall write
$\langle B,C\rangle=B+C$. Similarly, for some set of subgroups~$B_i$ of~$A$
we shall write $B=\langle B_i\rangle_{i\in I}=\sum\limits_{i\in I} B_i$.

The group $\langle a\rangle$ is the \emph{cyclic} group
generated by~$a$. The order of $\langle a\rangle$ is also called the
\emph{order} of~$a$ (notation: $o(a)$). If every element of~$A$ is of
finite order, then $A$
is called a~\emph{torsion} or \emph{periodic} group.
If all elements of~$A$, except~$0$,
are of infinite order.
then $A$ is \emph{torsion free}.
\emph{Mixed groups} contain both nonzero elements of
finite order and elements of infinite order.

A~\emph{primary group} or \emph{$p$-group} is defined to be
a~group, the orders of whose elements are powers of a~fixed prime~$p$.

Given $a\in A$, the greatest nonnegative integer~$r$ for which
$p^rx=a$ has a~solution $x\in A$
is called the \emph{$p$-height} $h_p(a)$ of~$a$.
If $p^r x=a$ is solvable whatever $r$ is, then $a$ is of
\emph{infinite $p$-height}, $h_{p}(a)=\infty$. If it is completely clear
from the context which prime~$p$ is meant, then we call $h_p(a)$ simply
the \emph{height} of~$a$ and write $h(a)$. For a~group~$A$ and an
integer $n> 0$, let $nA=\{ na\mid a\in A\}$ and
$A[n]=\{ a\mid a\in A,\ na=0\}$.

A~map $\alpha\colon A\to B$ is called a~\emph{homomorphism}
(of~$A$ into~$B$) if
$$
\forall a_1,a_2\in A \quad \alpha(a_1+a_2)=\alpha a_1+\alpha a_2.
$$
The \emph{kernel} of~$\alpha$ ($\Ker \alpha \lhd A$)
is the set
of all elements $a\in A$ such that $\alpha a=0$. The \emph{image} of~$\alpha$
($\Image \alpha \lhd B$) consists of all $b\in B$
that for some $a\in A$ satisfy $\alpha a=b$. If $\Image\alpha =B$, then
$\alpha$ is called a~\emph{surjective} homomorphism, or an \emph{epimorphism}.
If $\Ker\alpha=0$, then $\alpha$ is called an \emph{injective} homomorphism,
or a~\emph{monomorphism}. A~homomorphism that is injective and
surjective simultaneously is called an \emph{isomorphism}.

Now we consider the most important types of Abelian groups.

Cyclic groups were defined above as groups $\langle a \rangle$
for some~$a$. Note that all subgroups of cyclic groups are likewise
cyclic.

For a~fixed prime~$p$, consider the $p^n$th complex roots
of unity, $n\in \mathbb N \cup \{ 0\}$.
They form an infinite multiplicative group;
in accordance with our convention, we switch to the additive notation.
This group, called
a~\emph{quasicyclic} group or a~\emph{group of type~$p^\infty$} (notation:
$\mathbb Z(p^\infty)$) can be defined as follows. It is
generated by elements
$c_1,c_2,\dots,c_n,\ldots$, such that
$pc_1=0$, $pc_2=c_1,\dots,\allowbreak pc_{n+1}=c_n,\dots$.
Here $o(c_n)=p^n$, and every element of $\mathbb Z(p^\infty)$ is
a~multiple of some~$c_n$.
It is clear that all the quasicyclic groups
corresponding to the same prime~$p$
are isomorphic.

Let $p$ be a~prime, and $Q_p$ be the ring of rational numbers
whose denominators are prime to~$p$. The nonzero ideals of~$Q_p$
are principal ideals generated by~$p^k$ with $k=0,1,\ldots$.
If one considers the ideals $(p^k)$ as
a~fundamental system of neighborhoods of~$0$,
then $Q_p$ becomes a~topological ring, and we may form the completion
$Q_p^*$ of~$Q_p$ in this topology. $Q_p^*$~is again
a~ring whose ideals are $(p^k)$ with $k=0,1,\ldots$,
and which is complete in the topology defined by its ideals.

The elements of $Q_p^*$ may be represented as follows:
let $\{ t_0,t_1,\dots,t_{p-1}\}$ be a~complete set of representatives
of cosets of $p^kQ_p$ modulo $p^{k+1}Q_p$. Let $\pi \in Q_p^*$,
and let $a_n\in Q_p$ be a~sequence tending to~$\pi$.
According to the definition of fundamental sequence, almost all its elements
(i.e., all with a~finite number of exceptions)
belong to the same coset modulo $pQ_p$, say, to that represented by~$s_0$.
Almost all differences $a_n-s_0$ belonging to $pQ_p$
belong to the same coset of the ring $pQ_p$ modulo $p^2Q_p$, say, to
that represented by~$ps_1$. So proceeding, $\pi$~uniquely defines a~sequence
$s_0,ps_1,\ldots$, and we associate with~$\pi$
the formal infinite series $s_0+s_1p+\ldots$. Its partial sums
$b_n=s_0+s_1p+\dots+s_np^n$ form in $Q_p$ a~fundamental sequence
which converges in $Q_p^*$ to~$\pi$, in view of $\pi-b_n\in p^kQ_p^*$.
From the uniqueness of limits it follows
that, in this way, different elements of $Q_p^*$
are associated with different series. Therefore we can identify
the elements~$\pi$ of the ring~$Q_p^*$ with the formal series
$s_0+s_1p+s_2p^2+\dots$ with coefficients
from $\{ 0,1,\dots, p-1\}$ and write
$$
\pi=s_0+s_1p+s_2p^2+\dots \quad (s_n=0,1,\dots,p-1).
$$
The arising ring is a~\emph{commutative integral domain}
(i.e., the commutative ring without zero divisors) and is called
the \emph{ring of $p$-adic integers}.

\subsection{Direct Sums}\label{ss2.3}
Let $B$ and~$C$ be two subgroups of~$A$ such that:
\begin{enumerate}
\item
$B+C=A$;
\item
$B\cap C=0$.
\end{enumerate}

Then we call $A$ the \emph{direct sum} of its subgroups $B$~and~$C$
($A=B\oplus C$).

If the condition~(2) is satisfied, then we say that the groups
$B$~and~$C$ \emph{are disjoint}.

If $B_i$ ($i\in I$) is a~family of subgroups of~$A$ such that
\begin{enumerate}
\item
$\sum\limits_{i\in I} B_i=A$;
\item
$\forall i\in I\ B_i\cap \sum\limits_{j\ne i} B_j=0$,
\end{enumerate}
then we say that the group~$A$ is a~\emph{direct sum}
of its subgroups~$B_i$,
and write $A=\bigoplus\limits_{i\in I} B_i$,
or $A=B_1\oplus \dots \oplus B_n$,
if $I=\{ 1,\dots,n\}$. A~subgroup~$B$ of the group~$A$ is called
a~\emph{direct summand} of~$A$ if there exists a~subgroup $C\lhd A$
such that $A=B\oplus C$.
In this case, $C$ is called
the \emph{complementary direct summand} or simply
the \emph{complementation}.

Two direct decompositions of~$A$, $A=\bigoplus\limits_i B_i$
and $A=\bigoplus\limits_j C_j$, are called \emph{isomorphic} if
the components $B_i$~and~$C_j$
may be brought into a~one-to-one correspondence
such that the corresponding components are isomorphic.

If we have two groups $B$~and~$C$, then the set of all pairs
$(b,c)$, where $b\in B$ and $c\in C$,
forms a~group~$A$ if we set $(b_1,c_1)=(b_2,c_2)$ if and only if
$b_1=b_2$ and $c_1=c_2$, and $(b_1,c_1)+(b_2,c_2)=(b_1+b_2,c_1+c_2)$.

The correspondences $b\mapsto (b,0)$ and $c\mapsto (0,c)$ are isomorphisms
between the groups $B$,~$C$ and the subgroups $B'$,~$C'$ of~$A$.
We can write $A=B\oplus C$
and call $A$ the (\emph{external}) \emph{direct sum} of $B$~and~$C$.

Let $B_i$ ($i\in I$) be a~set of groups. A~\emph{vector} $(\dots,b,\dots)$
over the groups~$B_i$ is a~vector with $i$th coordinate for every $i\in I$
equal to some $b_i\in B_i$. Equality and addition of vectors is defined
coordinatewise. In this way, the set of all vectors becomes a~group~$C$,
called the \emph{direct product} of the groups~$B_i$,
$$
C=\prod_{i\in I} B_i.
$$

The correspondence
$\rho_i\colon b_i\mapsto (\dots,0,b_i,0,\dots)$, where $b_i$
stands on the $i$th place, and $0$ everywhere else, is an
isomorphism of the group~$B_i$ with a~subgroup~$B_i'$ of~$C$.
The groups~$B_i'$ ($i\in I$)
generate in~$C$ the group~$A$ of all vectors $(\dots,b_i,\dots)$
with $b_i=0$ for almost all $i\in I$,
and $A=\bigoplus\limits_{i\in I} B_i'$.
The group~$A$ is also called the (\emph{external})
\emph{direct sum} of~$B_i$.
If $A$ is a~group and $\varkappa$ is a~cardinal number, then by
$\bigoplus\limits_{\varkappa} A$ we shall denote
the direct sum of $\varkappa$~groups isomorphic to~$A$, and by
$\smash[b]{\prod\limits_{\varkappa}}A=A^\varkappa$
we shall denote the direct product
of $\varkappa$~such groups.

\begin{proposition}[\cite{Fuks}]\label{p2.1}
Every torsion group~$A$ is a~direct sum
of $p$-groups~$A_p$ belonging to different primes~$p$.
The groups~$A_p$ are uniquely determined by~$A$.
\end{proposition}

The subgroups~$A_p$ are called the $p$-\emph{components} of~$A$.
By virtue of Proposition~\ref{p2.1},
the theory of torsion groups is essentially reduced
to that of primary groups.

\begin{proposition}[\cite{Fuks}]\label{p2.2}
If there exists a~projection~$\pi$ of the group~$A$ onto its
subgroup~$B$, then $B$ is a~direct summand in~$A$.
\end{proposition}

\subsection{Direct Sums of Cyclic Groups}\label{ss2.4}
\begin{proposition}[\cite{Fuks}]\label{p2.3}
For a~group~$A$
the following conditions are equivalent\textup{:}
\begin{enumerate}
\item
$A$ is a~finitely generated group\textup{;}
\item
$A$ is a~direct sum of a~finite number of cyclic groups.
\end{enumerate}
\end{proposition}

A~system $\{ a_1,\dots,a_k\}$ of nonzero elements of a~group~$A$ is called
\emph{linearly independent} if
$$
n_1a_1+\dots+n_ka_k=0\quad (n_i\in \mathbb Z)
$$
implies
$$
n_1a_1=\dots=n_ka_k=0.
$$
A~system of elements is \emph{dependent} if it is not
independent.

An infinite
system $L=\{ a_i\}_{i\in I}$ of elements of the group~$A$
is called \emph{independent} if every finite subsystem of~$L$ is
independent. An independent system~$M$
of~$A$ is \emph{maximal} if there is no independent system in~$A$
containing~$M$ properly. By the \emph{rank} $r(A)$ of a~group~$A$
we mean the cardinality of a~maximal independent system
containing only elements of infinite and prime power orders.

\begin{proposition}\label{p2.4}
The rank $r(A)$ of a~group~$A$ is an invariant of this group.
\end{proposition}

\begin{theorem}[\cite{Prufer2,Baer1}]\label{t2.1}
A~bounded group is a~direct sum of cyclic groups.
\end{theorem}

\begin{theorem}[\cite{Fuks}]\label{t2.2}
Any two decompositions of a~group into direct sums of cyclic groups
of infinite and prime power orders are isomorphic.
\end{theorem}

\begin{theorem}[\cite{Kulikov2}]\label{t2.3}
Subgroups of direct sums of cyclic groups
are again direct sums of cyclic groups.
\end{theorem}

\subsection{Divisible Groups}\label{ss2.5}
We shall say that an element~$a$ of~$A$ is \emph{divisible by}
a~positive integer~$n$ ($n\mid a$)
if there is an element $x\in A$ such that $nx=a$.
A~group~$D$ is called \emph{divisible} if $n\mid a$
for all $a\in D$ and all natural~$n$.
The groups $\mathbb Q$ and $\mathbb Z(p^\infty)$ are examples
for divisible groups.

\begin{theorem}\label{t2.4}
Any divisible group~$D$ is a~direct sum of quasicyclic groups
and full rational groups~$\mathbb Q$. The powers of the sets of components
$\mathbb Z(p^\infty)$ \textup{(}for every~$p$\textup{)}
and~$\mathbb Q$ form a~complete and independent system of invariants
for the group~$D$.
\end{theorem}

\begin{corol}
Any divisible $p$-group $D$
is a~direct sum of the groups $\mathbb Z(p^\infty)$.
The power of the set of components $\mathbb Z(p^\infty)$ is
the only invariant of~$D$.
\end{corol}

\begin{theorem}[\cite{Fuks}]\label{t2.5}
For a~group~$D$ the following conditions are equivalent\textup{:}
\begin{enumerate}
\item
$D$ is a~divisible group\textup{;}
\item
$D$ is a~direct summand in every group containing~$D$.
\end{enumerate}
\end{theorem}

\subsection{Pure Subgroups}\label{ss2.6}
A~subgroup $G$ of~$A$ is called \emph{pure} if the equation $nx=g\in G$ is
solvable in~$G$ whenever it is solvable in the whole group~$A$.
In other words, $G$~is pure if and only if
$$
\forall n\in \mathbb Z\quad nG=G\cap nA.
$$

\begin{proposition}[\cite{Sele6}]\label{p2.5}
Assume that a~subgroup~$B$ of~$A$ is a~direct sum of cyclic
groups of the same order~$p^k$.
Then the following statements are equivalent\textup{:}
\begin{enumerate}
\item
$B$ is a~pure subgroup of~$A$\textup{;}
\item
$B$ is a~direct summand of~$A$.
\end{enumerate}
\end{proposition}

\begin{corol}
Every element of order~$p$ and of finite height can be embedded in a~finite
cyclic direct summand of the group.
\end{corol}

\begin{theorem}[\cite{Kulikov1}]\label{t2.6}
A~bounded pure subgroup is a~direct summand.
\end{theorem}

\begin{corollary}[\cite{Erdeli1}]
A~$p$-subgroup~$B$ of a~group~$A$
can be embedded in a~bounded direct summand
of~$A$ if and only if the heights of the nonzero elements of~$B$
\textup{(}relative to~$A$\textup{)} are bounded.
\end{corollary}

\begin{corollary}
An element~$a$ of prime power order belongs to a~finite
direct summand of~$A$ if and only if
$\langle a\rangle$ contains no elements of infinite height.
\end{corollary}

\subsection{Basic Subgroups}\label{ss2.7}
A~subgroup~$B$ of a~group~$A$ is called
a~$p$-\emph{basic subgroup} if it satisfies the following conditions:
\begin{enumerate}
\item
$B$ is a~direct sum of cyclic $p$-groups and infinite
cyclic groups;
\item
$B$ is pure in~$A$;
\item
$A/B$ is $p$-divisible.
\end{enumerate}

According to this definition, $B$ possesses a~basis which is said to be
a~\emph{$p$-basis} of~$A$.

Every group, for every prime~$p$, contains $p$-basic
subgroups~\cite{Fuks9}.

We now focus our attention on $p$-groups,
where $p$-basic subgroups are particularly
important. If $A$ is a~$p$-group and $q$ is a~prime different from~$p$,
then evidently $A$ has only one $q$-basic subgroup, namely~$0$.
Therefore, in $p$-groups we may refer to the $p$-basic subgroups
simply as \emph{basic} subgroups, without danger of confusion.

\begin{theorem}[\cite{Sele7}]\label{t2.7}
Assume that $B$ is a~subgroup of a~$p$-group~$A$,
$B=\bigoplus\limits_{n=1}^\infty B_n$, and $B_n$ is a~direct
sum of groups $\mathbb Z(p^n)$.
Then $B$~is a~basic subgroup of~$A$ if and only if
for every integer $n> 0$ the subgroup $B_1\oplus \dots\oplus B_n$
is a~maximal $p^n$-bounded direct summand of~$A$.
\end{theorem}

\begin{theorem}[Bear, Boyer \cite{Bojer1}]\label{Fuks32.4}
Suppose that $B$ is a~subgroup of a~$p$-group~$A$,
$$
B=B_1\oplus B_2\oplus\dots\oplus B_n\oplus \dots,
$$
where
$$
B_n\cong \bigoplus_{\mu_n} \mathbb Z(p^n).
$$
The subgroup $B$ is a~basic subgroup of~$A$ if and only if
$$
A=B_1\oplus B_2\oplus \dots\oplus B_n\oplus (B_n^*+p^n A),
$$
where $n\in \mathbb N$,
$$
B_n^*=B_{n+1}\oplus B_{n+2}\oplus \dots.
$$
\end{theorem}

Since the group $B$ has a~basis
and the factor group $A/B$ is a~direct sum
of groups isomorphic to $\mathbb Z(p^\infty)$ (i.e., $A/B$ also has
a~generating system which can be easily described),
it is natural to combine
these generating systems and to obtain one for~$A$.

We write
$$
B=\bigoplus_{i\in I} \langle a_i\rangle\ \
\text{and}\ \ A/B=\bigoplus_{j\in J} C_j^*,\ \
\text{where}\ \ C_j^*=\mathbb Z(p^\infty).
$$
If the direct summand $C_j^*$ is generated by cosets
$c_{j1}^*,\dots, c_{jn}^*,\dots$ modulo~$B$
with $pc_{j1}^*=0$ and $pc_{j,n+1}^*=c_{jn}^*$ ($n=1,2,\dots$),
then, by the purity of~$B$ in~$A$,
in the group~$A$ we can pick out $c_{jn}\in c_{jn}^*$
of the same order as~$c_{jn}^*$.
Then we get the following set of relations:
$$
pc_{j1}=0,\ \ pc_{j,n+1}=c_{jn}=b_{jn}\quad
(n\ge 1,\ b_{jn}\in B),
$$
where $b_{jn}$ must be of order at most~$p^n$, since
$o(c_{jn})=p^n$.

The system $\{ a_i,c_{jn}\}_{i\in I,\ j\in J,\ n\in \omega}$
will be called a~\emph{quasibasis} of~$A$.

\begin{proposition}[\cite{Fuks2}]\label{p2.6}
If $\{ a_i,c_{jn}\}$ is a~quasibasis of the $p$-group~$A$, then every
$a\in A$ can be written in the following form:
\begin{equation}\label{e2.1}
a=s_1a_{i_1}+\dots+s_ma_{i_m}+t_1 a_{j_1n_1}+\dots+t_r a_{j_rn_r},
\end{equation}
where $s_i$ and $t_j$ are integers, no $t_j$ is divisible by~$p$,
and the indices $i_1,\dots,i_m$ as well as $j_1,\dots,j_r$
are distinct. The expression~\eqref{e2.1} is unique in the sense
that $a$ uniquely defines the terms $s a_i$ and $tc_{jn}$.
\end{proposition}

\begin{theorem}[Kulikov \cite{Kulikov3}]\label{Fuks34.4s}
If $B$ is a~basic subgroup of a~reduced $p$-group~$A$, then
$$
|A|\le |B|^\omega.
$$
\end{theorem}

The \emph{final rank} of a~basic subgroup~$B$ of a~$p$-group~$A$ is the
infimum of the cardinals $r(p^nB)$. Note that if the rank of~$B$
is equal to~$\mu_1$ and the final rank of~$B$ is equal to~$\mu_2$, then
$A=A_1\oplus A_2$, where the group~$A_1$
is bounded and has the rank~$\mu_1$,
and a~basic subgroup of the group~$A_2$
has the rank~$\mu_2$ which coincides with its final rank.

\begin{theorem}\label{4.endom}
If two endomorphisms of a~reduced Abelian group coincide on some of its
basic subgroups, then they are equal.
\end{theorem}

\subsection{Endomorphism Rings of Abelian Groups}\label{ss2.8}
If we have an Abelian group~$A$, then its endomorphisms
form a~ring with respect to the operations of
addition and composition of homomorphisms. This ring will be denoted
by $\Endom(A)$.

We need some facts about the ring $\Endom(A)$.
\begin{enumerate}
\item
There exists a~one-to-one correspondence between
finite direct decompositions
$$
A=A_1\oplus \dots\oplus A_n
$$
of the group~$A$ and decompositions of the ring $\Endom(A)$
in finite direct sums of left ideals
$$
\Endom(A)=L_1\oplus \dots \oplus L_n;
$$
namely, if $A_i=e_iA$, where $e_1,\dots,e_n$ are mutually orthogonal
idempotents, then $L_i =\Endom(A)e_i$.
\item
An idempotent $e\ne 0$ is called primitive if it can not be represented
as a~sum of two nonzero orthogonal idempotents.
If $e\ne 0$ is an idempotent of the ring $\Endom(A)$, then $eA$
is an indecomposable direct summand
of~$A$ if and only if $e$ is a~primitive idempotent.
\item
Let $A=B\oplus C$ and $A=B'\oplus C'$ be direct decompositions
of the group~$A$, and let $e\colon A\to B$ and $e'\colon A\to B'$
be the corresponding projections.
Then $B\cong B'$ if and only if there exist elements
$\alpha,\beta\in \Endom(A)$ such that
$$
\alpha\beta=e\text{ and } \beta\alpha=e'.
$$
\end{enumerate}

\begin{theorem}[Baer \cite{Baer9}, Kaplansky \cite{Kaplans2}]\label{t2.8}
If $A$~and~$C$ are torsion groups, and $\Endom(A)\cong \Endom(C)$, then
$A\cong C$.
\end{theorem}

\begin{theorem}[Charles \cite{Sharl1}, Kaplansky \cite{Kaplans2}]\label{t2.9}
\hspace{-4pt}%
The center of the endomorphism ring $\Endom(A)\ttm\ttm$ of a~$p$-group~$A$
is the ring of $p$-adic integers or
the residue class ring of the integers modulo~$p^k$,
depending on whether $A$ is unbounded or bounded with $p^k$
as the least upper bound for the orders of its elements.
\end{theorem}

\section{Beautiful Linear Combinations}\label{sec3}
In this section, we shall completely follow
the paper~\cite{Shelah} of S.~Shelah.

Suppose that we have some fixed Abelian
$p$-group $A\cong \smash[b]{\bigoplus\limits_{\mu} \mathbb Z(p^l)}$,
where $\mu$ is an infinite cardinal number
and $\Endom(A)$ is its endomorphism ring.

For any $f\in \Endom(A)$ let $\Rng f$ be its range in~$A$ and
$\Cl B$ (or $\langle B\rangle$) be the closure of~$B\subset A$ in~$A$,
i.e., the minimal subgroup in~$A$ containing~$B$.

As usual, $\vec x$ denotes a~finite sequence of variables
$\vec x=\langle x_1,\dots, x_n\rangle$.
A~linear combination $k_1 x_1+\dots +k_n x_n$,
where $k_i\in \mathbb Z$, will also be denoted by
$\tau (x_1,\dots,x_n)$, or $\tau(\vec x)$. Such a~combination
will be called \emph{reduced} if all~$k_i$ are distinct and
not equal to zero.

Let $\{ a_i\mid i\in I\}$ be some independent subset
of elements of order~$p^l$ in the group~$A$.
It is clear that every function $h\colon \{ a_i\mid i\in I\}\to A$
has a~unique extension
$$
\tilde h\in \Hom (\Cl\{ a_i\mid i\in I\},A).
$$

Let $B$ be some set and $h$ be a~function from~$B$ into~$B$.
For every $x\in B$, we define its \emph{depth}
($\Dp(x)=\Dp(x,h)$) as the least ordinal number (or infinity)
satisfying the
following conditions:
\begin{enumerate}
\item
$\Dp(x)\ge 0$ if and only if $x\in B$;
\item
$\Dp(x)\ge \delta$ if and only if $\Dp(x)\ge \alpha$ for every
$\alpha\in \delta$ such that $\delta$ is a~limit ordinal number;
\item
$\Dp(x)\ge \alpha+1$ if and only if for some $y\in B$
we have $h(y)=x$ and $\Dp(y)\ge \alpha$.
\end{enumerate}

\begin{lemma}\label{shel_l1_1}
Let $\{ a_i\mid i\in I\}\subset A$
be an independent set consisting of elements of order~$p^l$,
and $h$ be a~function from~$I$ into~$I$. Define~$\tilde h$ by the formula
$$
\tilde h (k_1 a_{i_1}+\dots+k_n a_{t_n})=
k_1 a_{h(t_1)}+\dots +k_n a_{h(t_n)}.
$$

Then
\begin{enumerate}
\renewcommand{\theenumi}{\alph{enumi}}
\item
$\tilde h\in \Endom(\Cl\{ a_i\mid i\in I\})$\textup{;}
\item
$\Dp(k_1 a_{t_1}+\dots +k_n a_{t_n},\tilde h)\ge
\min\limits_{i\in \{ 1,\dots,n\}}\Dp(t_i,h)$\textup{;}
\item
if in~\textup{(b)} the linear combination
$k_1 a_{t_1}+\dots+k_n a_{t_n}$ is reduced and $t_i$
are distinct, then the equality holds.
\end{enumerate}
\end{lemma}

\begin{proof}
(a) Immediate.

(b) We prove by induction on ordinal numbers~$\alpha$ that
$$
\min_{i\in \{ 1,\dots,n\}} \Dp(t_i,h)\ge \alpha \Rightarrow
\Dp(k_1 a_{t_1}+\dots+k_n a_{t_n},\tilde h)\ge \alpha.
$$
This suffices for the proof of~(b).
If $\alpha=0$ or $\alpha$ is a~limit ordinal,
then this is trivial. For $\alpha=\beta+1$,
by the assumption and definition of depth,
there are $s_i\in I$ such that
$ h(s_i)=t_i$ and $\Dp(s_i,h)\ge \beta$. Then
$\min\limits_i \Dp(s_i,h)\ge \beta$,
whence by the induction hypothesis
$$
\Dp(k_1 a_{s_1}+\dots+k_n a_{s_n},\tilde h)\ge \beta,
$$
but
$$
\tilde h(k_1 a_{s_1}+\dots+k_n a_{s_n})=
k_1 a_{t_1}+\dots+k_na_{t_n},
$$
whence
$$
\Dp (k_1 a_{t_1}+\dots+k_n a_{t_n},\tilde h)\ge \beta+1=\alpha.
$$

(c) It suffices to prove by induction on~$\alpha$ that
$$
\Dp(k_1 a_{t_1}+\dots+k_n a_{t_n},\tilde h)\ge \alpha
\Rightarrow \Dp (t_i,h)\ge \alpha
$$
for all~$i$. If $\alpha=0$ or $\alpha$ is a~limit ordinal,
then this is trivial.
For $\alpha =\beta+1$, by the definition of depth, there are
a~reduced linear combination $l_1 a_{s_1}+\dots+l_m a_{s_m}$ and
distinct~$s_i$ such that
\begin{enumerate}
\item
$\tilde h(l_1 a_{s_1} +\dots +l_m a_{s_m})=
k_1 a_{t_1}+\dots+ k_n a_{t_n}$ and
\item
$\Dp(l_1 a_{s_1}+\dots +l_m a_{s_m},\tilde h)\ge\beta$.
\end{enumerate}

By (2) and the induction hypothesis, for $i=1,\dots,m$ we have
$\Dp (s_i,h)\ge \beta$.

By (1) and the definition of~$\tilde h$,
$$
l_1 a_{h(s_1)}+\dots +l_m a_{h(s_m)}=
k_1 a_{t_1}+\dots+k_n a_{t_n}.
$$

As the linear combination $k_1 a_{t_1}+\dots+k_n a_{t_n}$ is reduced
and the~$t_i$ are distinct,
$$
\{ t_1,\dots,t_n\}\subseteq \{ h(s_1),\dots,h(s_m)\}.
$$

Thus, for each $i=1,\dots,n$ there is $k_i$, $1\le k_i\le m$, such that
$t_i=h(s_{k_i})$. Hence
$$
\Dp(t_i,h)\ge \Dp(s_{k_i},h)+1\ge \beta+1=\alpha.\qed
$$
\renewcommand{\qed}{}
\end{proof}

\begin{lemma}\label{shel_l1_2}
Let $h_1$ and~$h_2$ be commuting functions from~$B$ into~$B$.
Then for any $x\in B$ we have
$$
\Dp(x,h_1)\le \Dp(h_2(x),h_1).
$$
\end{lemma}

\begin{proof}
We prove by induction on~$\alpha$ that
$$
\Dp (x,h_1)\ge \alpha\Rightarrow \Dp(h_2(x),h_1)\ge \alpha.
$$
If $\alpha=0$ or $\alpha$ is a~limit ordinal, this is immediate.

Now we consider $\alpha=\beta+1$.

If $\Dp(x,h_1)\ge \beta+1$, then for some $y\in B$ we have
$h_1(y)=x$ and $\Dp(y,h_1)\ge \beta$. Thus
$h_1(h_2(y))=h_1\circ h_2(y)=h_2\circ h_1(y)=h_2(x)$,
and by the induction hypothesis
$\Dp(h_2(y),h_1)\ge \beta$ (since $\Dp(y,h_1)\ge \beta)$), whence
$\Dp(h_2(x),h_1)\ge \Dp(h_2(y),h_1)+1\ge \beta+1=\alpha$.
\end{proof}

\begin{lemma}\label{shel_l1_3}
Let $\{ a_i\mid i\in I\}\subset A$ be an independent set with elements
of the same order~$p^l$ and let $A'=\Cl\{ a_i\mid i\in I\}$.
Let $J\subseteq I$, $|I\setminus J|=|I|$,
$J=\smash[b]{\bigcup\limits_{\alpha\in \alpha(0)} J_\alpha}$, and
$B=\Cl\{ a_i\mid i\in J\}$. Then we can find
$f\in \Endom(A')$ such that
\begin{enumerate}
\renewcommand{\theenumi}{\alph{enumi}}
\item
$i\in J_\alpha\Rightarrow \Dp (a_i,f)=\alpha$\textup{;}
\item
if $g\in \Endom(A')$, $g\circ f=f\circ g$, and $g$ maps~$B$
into~$B$, then for every $\alpha\in \alpha(0)$
$$
i\in J_\alpha\Rightarrow
g(a_i)\in
\Cl\biggl\{ a_j \biggm|
j\in \bigcup\limits_{\alpha\le\beta < \alpha(0)} J_\beta\biggr\};
$$
\item
every function $g\colon \{ a_i\mid i\in J\}\to B$ satisfying
the condition~\textup{(b)}
can be extended to an endomorphism of $\Endom (A')$
that commutes with~$f$ and maps~$B$ into~$B$.
\end{enumerate}
\end{lemma}

\begin{proof}
By renaming we can assume that
$I\setminus J=I_0\cup \{ \langle 0,t,\eta\rangle \mid
t\in J_\alpha,\ \alpha< \alpha(0),\ l(\eta)> 0$,
$\eta$~is a~decreasing sequence of ordinals,
$\eta_0< \alpha\} \cup \{ \langle 1,t,n\rangle \mid 0< n< \omega\}$
and we identify $\langle 0,t,\langle\ \rangle \rangle$ and
$\langle 1,t,0\rangle$ with~$t$.
Let us define a~function $h$ on~$I$:
\begin{enumerate}
\item
for $s\in I_0$
$$
h(s)=s;
$$
\item
for $t\in J$, $\langle 0,t, \langle \eta, \beta\rangle \rangle\in I$,
where $\beta$ is an ordinal number,
$$
h(\langle 0,t,\langle \eta, \beta\rangle \rangle)=\langle 0,t,\eta\rangle;
$$
\item
for $t\in J$
$$
h(\langle 1,t,n\rangle)=\langle 1,t,n+1\rangle .
$$
\end{enumerate}

Let $f=\tilde h$. It is easy to prove that if $t\in J$ and $l(\eta)=n$, then
$\Dp (\langle 0,t,\eta\rangle, h)=\eta_n$. Hence (a)~follows immediately
by~(c) of Lemma~\ref{shel_l1_1},
and $f\in \Endom(A')$ by~(a) of Lemma~\ref{shel_l1_1}.

Let now $g\in \Endom(A')$ commute with~$f$, $g$ map~$B$ into~$B$,
$i\in J_\alpha$, $\alpha\in \alpha_0$. Suppose that
the linear combination $g(a_i)=k_1 a_{s_1}+\dots+k_n a_{s_n}$
is reduced. Since $g$ maps~$B$ into~$B$ and $g(a_i)\in B$, we see that
$s_k\in J$ for $k=1,\dots, n$. Therefore we can assume that
$s_k\in J_{\alpha_k}$.
By~(c) from Lemma~\ref{shel_l1_1} and the remarks above,
$$
\Dp(g(a_i),f)=\Dp(k_1 a_{s_1}+\dots+k_n a_{s_n}, f)=
\min_k \Dp (s_k, h)=\min_k \alpha_k.
$$

On the other hand, by Lemma~\ref{shel_l1_2}, as $g$ commutes with~$f$,
$$
\alpha\le \Dp(a_i,f)\le \Dp(g(a_i),f).
$$

Combining both, we get $\alpha \le \alpha_k$ for
$k=1,\dots,n$. Hence
$g(a_i)\in \Cl \{ a_i\mid i\in J_\beta,\ \alpha \le \beta< \alpha_0\}$,
so we have proved~(b).

(c) Extend $g$ to a~function $g_1\colon \{ a_i\mid i\in I\}\to A'$
in the following way:
if $g(a_t)=k_1 a_{t_1}+\dots+k_na_{t_n}$, then let
\begin{align*}
g_1(a_{\langle 0,t,\eta\rangle}) &=
k_1 a_{\langle 0,t_1, \eta\rangle}+ \dots+
k_n a_{\langle 0,t_n,\eta\rangle},\\
g_1(a_{\langle a,t,m\rangle}) &=
k_1 a_{\langle 1,t_1,m\rangle}+\dots+
k_n a_{\langle 1,t_2,m\rangle},\\
g_1(a_s) &=a_s \text{ for } s\in I_0.
\end{align*}

It is easy to check that $g_1$ is well defined
(because $t\in J_\alpha$ and $t_k\in J_\beta$
imply $\alpha\le \beta$) and it has a~unique extension to
$g_2\in \Endom(A')$. In order to check that $g_2$~and~$f$
commute, it suffices to prove that for every $s\in I$
$$
f\circ g_2(a_s)=g_2\circ f(a_s),
$$
and this is quite easy.
\end{proof}

A~linear combination $\tau(x_1,\dots,x_n)=k_1 x_1+\dots+k_nx_n$
is called \emph{beautiful} (see~\cite{Shelah}, where similar terms
were called \emph{beautiful terms}) if
\begin{enumerate}
\item
we have
$\tau(\tau(x_1^1,\dots,x_n^1),\tau(x_1^2,\dots, x_n^2),\dots
\tau(x_1^n,\dots, x_n^n))=
\tau(x_1^1,\dots,x_n^n)$;
\item
we have
$\tau(x,\dots,x)=x$.
\end{enumerate}

The condition~(2) implies
$k_1+\dots+k_n=1$.
It is clear that the condition~(1) implies
$k_i k_j=\delta_{ij} k_i$
for all $i,j=1,\dots,n$.

We can see that all $k_j$, except one (let it be~$k_i$),
are equal to zero, and $k_i=1$, i.e., all beautiful linear combinations
have the form~$x_i$ for some $i\in I$.
Therefore the following lemma is trivial.

\begin{lemma}\label{shel_l3_1}
The set of beautiful linear combinations is closed under substitution.
\end{lemma}

\begin{theorem}\label{shel_t4_1}
There is a~formula $\tilde \varphi({\dots})$ such that the
following holds.
Let $\{ f_i\}_{i\in \mu}$ be a~set of elements of $\Endom(A')$.
Then we can find a~vector~$\bar g$ such that the formula
$\tilde \varphi( f,\bar g)$
holds in $\Endom(A')$ if and only if $ f=f_i$ for some $i\in \mu$.
\end{theorem}

\begin{proof}
Suppose that a~set $\{ a_i\mid i\in I^*\}\subset A$
is independent, consists
of elements of the same order, and $\Cl\{ a_i\mid i\in I^*\}=A'$.
Let $J\subseteq I^*$ and
$\mu =|J|=|I^*\setminus J|$. For notational simplicity let
$J=\{ \langle \alpha,\beta\rangle\mid {\alpha,\beta\in \mu}\}$,
$a_{\langle \alpha,\beta\rangle}=a_\alpha^\beta$.

\begin{lemma}\label{shel_l4_2}
There is a~formula $\varphi(f)$ with one free variable~$f$
such that $\varphi(f)$ holds in $\Endom(A')$ if and only if
there is an ordinal number $\alpha\in \mu$ such that for all
$\beta \in \mu$
$$
f(a_0^\beta)=a_\alpha^\beta.
$$
\end{lemma}

Now we shall show the proof of Theorem~\ref{shel_t4_1} with the help
of Lemma~\ref{shel_l4_2}.

Let a~function $f_0^*\colon A'\to A'$
map the set $\{ a_0^\alpha\mid \alpha\in \mu\}$
onto the set $\{ a_t\mid t\in I^*\}$,
and let $f_0^*(a_\alpha^\beta)=f_0^*(a_0^\beta)$.
Suppose that we have a~set $\{ f_i\}_{i\in \mu}$ and let the function
$f^*\colon A'\to A'$ be such that
$$
f^*(a_\alpha^\beta)=f_\alpha\circ f_0^*(a_\alpha^\beta).
$$
Let $f_1^*\colon A\to A$ map the set $\{ a_t\mid t\in I^*\}$ onto the set
$\{ a_0^\beta\mid \beta\in \mu\}$.
Let the formula $\tilde \varphi(f',{\dots})$ say that there exists
$f\in \Endom(A')$ such that
\begin{enumerate}
\item
$\varphi(f)$;
\item
$f'\circ f_0^*\circ f_1^*=f^*\circ f\circ f_1^*$.
\end{enumerate}

Then $\vDash \varphi(f)$ if and only if
there exists $\alpha\in \mu$ such that
$$
\forall \beta\in \mu\
f(a_0^\beta)=a_{\alpha}^{\beta}.
$$
Therefore
\begin{multline*}
f'\circ f_0^*\circ f_1^*(a_t)=f^*\circ f\circ f_1^*(a_t)
\Leftrightarrow
f'\circ f_0^*(a_0^\beta)=f^*\circ f(a_0^\beta)\\
\Leftrightarrow f'\circ f_0^*(a_0^\beta)=f^*(a_{\alpha}^\beta)
\Leftrightarrow f'\circ f_0^*(a_0^\beta)=f_{\alpha}\circ f_0^*(a_0^\beta).
\end{multline*}

Let $f_0^*(a_0^\beta)=a_{t_\beta}$. Then
$$
f'(a_{t_\beta})= f_{\alpha}(a_{t_\beta}),
$$
what we needed.
\end{proof}

\begin{proof}[Proof of Lemma~\ref{shel_l4_2}]
We partition our proof to three cases.

Case I. $\mu=\omega$.

Introduce the following mappings
\begin{enumerate}
\item
$f_2^*\in \Endom(A')$ such that $f_2^*(a_t)=a_0^0$ for all $t\in I^*$;
\item
$f_3^*\in \Endom(A')$ such that
$t\in I^*\setminus J\Rightarrow f_3^*(a_t)=a_0^0$,
$t\in J\Rightarrow f_3^*(a_n^m)=a_{n+1}^m$;
\item
$f_4^*\in \Endom(A')$ such that
$t\in I^*\setminus J\Rightarrow f_4^*(a_t)=a_0^0$,
$t\in J\Rightarrow f_4^*(a_n^m)=a_n^{m+1}$;
\item
$f_5^*\in \Endom(A')$ such that
$t\in I^*\setminus J\Rightarrow f_5^*(a_t)=a_0^0$,
$t\in J\Rightarrow f_5^*(a_n^m)=a_m^n$;
\item
$f_6^*\in \Endom(A')$ such that
$t\in I^*\setminus J\Rightarrow f_6^*(a_t)=a_0^0$,
$t\in J\Rightarrow f_6^* (a_n^m)=a_n^n$;
\item
$f_7^*\in \Endom(A')$ such that
$t\in I^*\setminus J\Rightarrow f_7^*(a_t)=a_0^0$,
$t\in J\Rightarrow f_7^*(a_n^m)=a_n^0$.
\end{enumerate}

Suppose that
\begin{align*}
B_0&=\Cl\{ a_0^0\},\\
B_1&= \Cl\{ a_n^0\mid n\in \omega\},\\
B_2&= \Cl\{ a_n^m\mid n,m\in \omega\},\\
B_3&= \Cl\{ a_0^m\mid m\in \omega\}.
\end{align*}

Introduce now the following first order formulas.
\begin{enumerate}
\renewcommand{\labelenumi}{\theenumi.}
\item
$\varphi_1(f;\dots)$ says that the restriction of the function~$f$
on the set~$B_1$ is a~function with the image in~$B_1$, commuting
with $f_3^*|_{B_1}$.
As a~formula we have the following:
$$
\varphi_1(f;\dots)=
(\rho_{B_1}\circ f\circ \rho_{B_1}=f\circ \rho_{B_1}) \logic\land
(f\circ \rho_{B_1}\circ f_3^*\circ \rho_{B_1}=
f_3^*\circ \rho_{B_1}\circ f\circ \rho_{B_1}),
$$
where $\rho_{B_1}$ is the projection onto~$B_1$.
\item
Similarly, $\varphi_2(f;\dots)$ says that the formula
$\varphi_1(f;\dots)$ holds and the image of $f|_{B_2}$
is included in~$B_2$ and $f_{B_2}$ commutes with $f_4^*|_{B_2}$.
\item
The formula $\varphi_3(f;\dots)$ says that the formula $\varphi_2(f;\dots)$
holds and $(f_5^*\circ f\circ f_5^*\circ f)|_{B_0}=(f_6^*\circ f)|_{B_0}$.
\item
The formula $\varphi_4(f;\dots)$ says that the formula
$\varphi_3(f;\dots)$ holds
and $(f_2^*\circ f)|_{B_0}=f_2^*|_{B_0}$.
\item
The formula $\varphi_5^*(f;\dots)$ says that there exists
$f'\in \Endom(A')$ such that $f'|_{B_3}=f|_{B_3}$ and the formula
$\varphi_4(f';\dots)$ holds.
\end{enumerate}

Now we note the following.

1. The formula $\varphi_1(f;\dots)$ holds if and only if
$$
f(a_n^0)=q_1a_{n+l_1}^0+\dots+q_k a_{n+l_k}^0
$$
for some $q_1,\dots,q_k\in \mathbb Z$, $l_1,\dots,l_k\in \omega$ for any
$n\in \omega$.

Actually, $f\colon B_1\to B_1$ means that
$f(a_n^0)=s_1 a_{n_1}^0+\dots+s_k a_{n_k}^0$,
and $f|_{B_1}\circ f_3^*|_{B_1}=f_3^*|_{B_1}\circ f|_{B_1}$ is equivalent
to the condition that for any $n\in \omega$
$$
f\circ f_3^*(a_n^0)=f_3^*\circ f(a_n^0).
$$
Let $f(a_0^0)=r_1 a_{l_1}+\dots+r_k a_{l_k}$. Then
$f(a_1^0)=f_3^*\circ f(a_0^0)=r_1 a_{1+l_1}+\dots+r_k a_{1+l_k}$ and so on
by induction.

2. The formula $\varphi_2(f;\dots)$ holds if and only if
$$
f(a_n^m)=r_1 a_{n+l_1}^m+\dots+r_k a_{n+l_k}^m.
$$
Actually, from $\varphi_1(f)$ it follows that for any $n\in \omega$
$$
f(a_n^0)=r_1 a_{n+l_1}^0+\dots+ r_k a_{n+l_k}^0,
$$
and
$$
f(a_n^1)=f\circ f_4^*(a_n^0)=f_4^* \circ f(a_n^0)=
f_4^* (r_1 a_{n+l_1}^0+\dots+ r_k a_{n+l_k}^0)=
r_1 a_{n+l_1}^1+\dots+r_k a_{n+l_k}^1,
$$
and so on by induction.

3. The formula $\varphi_3(f;\dots)$ holds if and only if
the formula $\varphi_2(f;\dots)$ holds, where $l_1,\dots,l_k$ are distinct
and the corresponding linear combination~$\tau$
satisfies the condition~(1) of the definition of
beautiful linear combination.

Actually,
let the formula $\varphi_2(f)$ be true and $\tau$ be the corresponding
linear combination. Then
\begin{align*}
& f_5^*\circ f\circ f_5^*\circ f|_{B_0}=
f_6^*\circ f|_{B_0}\Leftrightarrow
f_5^*\circ f\circ f_5^* \circ f(a_0^0)=f_6^*\circ f(a_0^0)\\
& \quad {}\Leftrightarrow
f_5^*\circ f\circ f_5^*(\tau(a_{l_1}^0,\dots,a_{l_k}^0))=
f_6^*(\tau(a_{l_1}^0,\dots,a_{l_k}^0))\Leftrightarrow
f_5^*\circ f(\tau(a_0^{l_1},\dots,a_0^{l_k}))=
\tau(a_{l_1}^{l_1},\dots,a_{l_k}^{l_k})\\
& \quad {}\Leftrightarrow
f_5^*(\tau(\tau(a_{l_1}^{l_1},\dots,a_{l_1}^{l_k}),\dots,
\tau(a_{l_k}^{l_1},\dots, a_{l_k}^{l_k}))=
\tau(a_{l_1}^{l_1},\dots,a_{l_k}^{l_k})\\
& \quad {}\Leftrightarrow
\tau (\tau(a_{l_1}^{l_1},\dots, a_{l_k}^{l_1}),\dots,
\tau(a_{l_1}^{l_k},\dots,a_{l_k}^{l_k}))=
\tau(a_{l_1}^{l_1},\dots,a_{l_k}^{l_k}).
\end{align*}
It is equivalent to the condition~(1)
from the definition of linear combination.
The converse condition is proved in the same manner.

4. The formula $\varphi_4(f;\dots)$ holds if and only if
$$
f(a_n^m)=\tau(a_{n+l_1}^m,\dots,a_{n+l_k}^m),
$$
where $\tau(x_1,\dots,x_k)$ is a~beautiful linear combination, i.e.,
$f(a_n^m)=a_{n+l_s}^m$.

Since $\varphi_4(f)\Rightarrow \varphi_3(f)$, we only need to show
that $a_0^0=\tau(a_0^0,\dots,a_0^0)$.

Actually,
$$
f_2^*\circ f|_{B_0}=f_2^*|_{B_0}\Leftrightarrow f_2^*\circ f(a_0^0)=
f_2^*(a_0^0)
\Leftrightarrow f_2^*(\tau(a_{l_1}^0,\dots,a_{l_k}^0))=a_0^0\Leftrightarrow
\tau(a_0^0,\dots,a_0^0)=a_0^0.
$$

5. The formula $\varphi_5(f)$ holds if and only if
$$
f(a_0^n)=\tau(a_{l_1}^n,\dots,a_{l_k}^n)
$$
for some beautiful linear combination~$\tau$ and $l_1,\dots,l_k\in \omega$
for all $n\in \omega$.

This follows immediately from
$f_5^*\circ f_7^*\circ f_5^*(a_0^m)=f_5^*\circ f_7^*(a_m^0)=a_0^m$.

\smallskip

Case II. The cardinal number $\mu=|J|$ is regular and $\mu> \omega$.

Let $I^*\setminus J=I_0\cup
\{ \langle \alpha,\delta,n\rangle\mid \alpha\in \mu,\
\delta\in \mu,\ \cf\delta=\omega,\ n\in \omega\}$,
$|I_0|=\mu$, and let us denote
$a_\alpha^{\beta,n}=a_{\langle \alpha,\beta,n\rangle}$.

For every limit ordinal $\delta\in \mu$ such that $\cf \delta=\omega$,
choose an increasing sequence $(\delta_n)_{n\in \omega}$ of ordinals
less than~$\delta$
such that their limit is~$\delta$ and for each
$\beta\in \mu$, $n\in \omega$ the set
$\{ \delta\in \mu\mid \beta=\delta_n\}$
is a~stationary subset in~$\mu$ (see~\cite{Solovay}).

Let us define some $f_i^*$ by defining $f_i^*(a_t)$
for some $t\in I^*$ and understanding that when $f_i^*(a_t)$ is
not explicitly defined, it is $a_0^0$.

So let for $\alpha,\beta\in \mu$
$$
f_2^*(a_\alpha^\beta)=a_0^0;\quad
f_3^*(a_\alpha^\beta)=a_\alpha^0;\quad
f_4^*(a_\alpha^\beta)=a_0^\beta;\quad
f_5^*(a_\alpha^\beta)=a_\beta^\alpha;\quad
f_6^*(a_\alpha^\beta)=a_\alpha^\alpha;
$$
let for $\delta\in \mu$, $\cf \delta=\omega$
$$
f_7^*(a_\alpha^\delta) =a_\alpha^{\delta,0};
$$
and let for $\delta\in\mu$, $\cf \delta\ne\omega$
$$
f_7^*(a_\alpha^\delta)=a_\alpha^\delta;
$$
and, further,
\begin{gather*}
f_8^*(a_\alpha^\beta)= a_\alpha^\beta,\quad
f_8^*(a_\alpha^{\delta,n})=a_\alpha^{\delta,n+1};\\
f_9^*(a_\alpha^\beta)=a_\alpha^\beta,\quad
f_9^*(a_\alpha^{\delta,n})=a_\alpha^{\delta_n}.
\end{gather*}

Let
\begin{align*}
B_0&= \Cl\{ a_0^0\},\\
B_1&= \Cl\{ a_\alpha^0\mid  \alpha\in \mu\},\\
B_2&= \Cl\{ a_0^\beta\mid \beta\in \mu\},\\
B_3&= \Cl\{ a_\alpha^\beta\mid \alpha,\beta\in \mu\},\\
B_4&= \Cl\{ a_\alpha^\beta,a_\alpha^{\beta,n}\mid
\alpha,\beta\in \mu,\ n\in \omega\},\\
B_5&= \Cl\{ a_0^\beta,a_0^{\beta,n}\mid \beta\in \mu,\ n\in \omega\},\\
B_6&= \Cl\{ a_0^\beta\mid \beta\in \mu,\ \cf \beta=\omega\},\\
B_7&= \Cl\{ a_\alpha^\beta\mid \alpha,\beta\in \mu,\ \cf \beta=\omega\}.
\end{align*}
Clearly $f_3^*$, $f_4^*$, and $f_9^*$ are projections onto $B_1$, $B_2$,
and $B_3$, respectively.
Let $f_{10}^*$, $f_{11}^*$, and $f_{12}^*$ be projections onto
$B_4$, $B_5$, and $B_6$, respectively.

Now we apply Lemma~\ref{shel_l1_3} with
$J=\{ \langle \alpha,\beta\rangle\mid \alpha,\beta \in \mu\}$,
$J_\beta=\{ \langle \alpha,\beta\rangle\mid \alpha\in \mu\}$, $I=I^*$,
and $f=f_{13}^*\in \Endom(A')$.

Let the first order formula $\varphi^1(f,g;\dots)$ say that
\begin{enumerate}
\item
$f$ and $g$ are conjugate to~$f_2^*$;
\item
$\Rng f, \Rng g\subset B_3$;
\item
$\exists h\in \Endom(A')\
(h\circ f_{13}^*=f_{13}^*\circ h \logic\land
\Rng h|_{B_3}\subseteq B_3 \logic\land h\circ f =g)$.
\end{enumerate}

We shall write $\varphi^1(f,g;\dots)$ also in the form $f\le g$.

If $f$ and $g$ are conjugate to $f_2^*$,
$f(a_0^0)=a_\alpha^\beta$, and
$g(a_0^0)=\tau(a_{\alpha_1}^{\beta_1},\dots, a_{\alpha_k}^{\beta_k})$,
then $f\le g$
if and only if
$\beta\le \beta_1,\dots,\allowbreak \beta\le \beta_k$,
which is easy to see from Lemma~\ref{shel_l1_3}.

Let the first order formula $\varphi_2(f)$ say that
\begin{enumerate}
\item
$\Rng f|_{B_2}\subseteq B_3$; $\Rng f|_{B_5}\subseteq B_4$;
$\Rng f|_{B_0}\subseteq B_1$; $\Rng f|_{B_6}\subseteq B_7$;
\item
for any $g\in \Endom (A')$ conjugate to~$f_2^*$, if
$\Rng g\subseteq B_2$, then $g\le f\circ g$;
\item
$f|_{B_5}$ commutes with $f_7^*$, $f_8^*$, and $f_9^*$;
\item
$f|_{B_3}$ commutes with $f_3^*$.
\end{enumerate}

\begin{statement}\label{shel_s4_3}
The formula $\varphi_2(f)$ holds in $\Endom(A')$ if and only if
for any $\beta\in \mu$
$$
f(a_0^0)=\tau(a_{\alpha_1}^\beta,\dots,a_{\alpha_k}^\beta)\ \
\text{and}\ \
f(a_0^{\beta,n})=\tau(a_{\alpha_1}^{\beta,n},\dots,a_{\alpha_k}^{\beta,n})
$$
for some linear combination~$\tau$ and ordinal numbers
$\alpha_1,\dots,\alpha_k \in \mu$
\textup{(}which do not depend on~$\beta$\textup{)}.
\end{statement}

\begin{proof}
Assume that $\varphi_2(f)$ holds for a~given function~$f$,
and let
$$
f(a_0^\beta)=
\tau_\beta(a_{\alpha_{\beta,1}}^{\gamma_{\beta,1}},\dots,
a_{\alpha_{\beta,k_\beta}}^{\gamma_{\beta,k_\beta}}).
$$
Without loss of generality we can assume that
the ordinal numbers $\alpha_{\beta,i}$ grow
as $i$~grows.

Choose $g$ such that $g(a_t)=a_0^\beta$ for all $t\in I^*$.
Then $g$ is conjugate to
$f_2^*$, $\Rng g\subseteq B_2$, and therefore,
by~(2), $g\le f\circ g$.
Since $g(a_0^0)=a_0^\beta$ and
$f(a_0^0)=\tau_\beta(a_{\alpha_{\beta,1}}^{\gamma_{\beta,1}},\dots,
a_{\alpha_{\beta,k_\beta}}^{\gamma_{\beta,k_\beta}})$,
we have that $\beta\le \gamma_{\beta,i}$ for all $i=1,\dots,k_\beta$.

Since the cardinal number~$\mu$ is regular and $\mu> \omega$, we have
that for any $\beta_0< \mu$
$$
\sup \{ \gamma_{\beta,i}\mid i=1,\dots,k_\beta;\ \beta\in \beta_0\}< \mu .
$$
Hence the set
$$
S=\{ \beta_0\mid \beta_0\in \mu \logic\land
\forall \beta\losp (\beta\in \beta_0 \logic\land
1\le i\le k_\beta\Rightarrow \gamma_{\beta,i}\in \beta_0)\}
$$
is an unbounded subset in~$\mu$; and by its definition it is closed
(indeed,
$\beta_1=
\sup\{ \gamma_{\beta,i}\mid i=1,\dots,k_\beta,\ \beta\in \beta_0\}$,
$\beta_2=
\sup \{ \gamma_{\beta,i}\mid i=1,\dots,k_\beta,\ \beta\in \beta_1\}$,~\dots,
$\bar \beta=\bigcup\limits_{l\in \omega} \beta_l$,
whence $\bar \beta< \mu$, $\bar \beta \in S$, and so on).

Now we shall prove that for $\delta\in S$, $\cf\delta=\omega$ implies
$\gamma_{\delta,i}=\delta$. Suppose that
$\gamma=\gamma_{\delta,i_0}\ne \delta$
for some~$i_0$,
which yields $\delta< \gamma_{\delta,i_0}$ as was said above.
As $\delta_n$ is increasing (as a~function of~$n$)
and its limit is~$\delta$, for some $n\in \omega$ big enough,
\begin{equation}\label{shel1}
\delta< \gamma_{\delta,i_0}(n)
\end{equation}
and $\gamma_{\delta,i_1}\ne \gamma_{\delta,i_2}\Rightarrow
\gamma_{\delta,i_1}(n)\ne \gamma_{\delta,i_2}(n)$.

Since $\Rng f|_{B_6}\subseteq B_7$,
we see that necessarily $\gamma_{\delta,i}$ has
cofinality~$\omega$, and as $f|_{B_5}$ commutes with~$f_7^*$,
we see that
$$
f(a_0^{\delta,0})=f\circ f_7^*(a_0^\delta)=
f_7^*\circ f(a_0^\delta)=
\tau_\delta(a_{\alpha_{\delta,1}}^{\gamma_{\delta,1},0},\dots,
a_{\alpha_{\delta,s}}^{\gamma_{\delta,s},0}).
$$
As $f|_{B_5}$ commutes with $f_8^*$,
$$
f(a_0^{\delta,n})=
f\circ f_8^*(a_0^{\delta,n-1})=
f_8^*\circ f(a_0^{\delta,n-1}),
$$
whence
$$
f(a_0^{\delta,n})=
\tau_\delta (a_{\alpha_{\delta,1}}^{\gamma_{\delta,1},n},\dots).
$$
Since $f|_{B_5}$ commutes with $f_9^*$, we have
$$
f(a_0^{\delta_n})=
\tau_\delta(a_{\alpha_{\delta,1}}^{\gamma_{\delta,1}(n)},\dots)
$$
and the linear combination
$\tau_\delta(a_{\alpha_{\delta,1}}^{\gamma_{\delta,1}(n)},\dots)$
is reduced.

But $f(a_0^{\delta_n})=
\tau_{\delta_n}(a_{\alpha_{\delta_n,1}}^{\gamma_{\delta_n,1}},\dots)$,
where the linear combination on the right-hand side of
the equality is also reduced. Hence
$$
\gamma_{\delta,i_0}(n)\in
\{ \gamma_{\delta_n,i}\mid i\le i\le k_{\delta_n}\},
$$
but on the one hand, $\delta< \gamma_{\delta,i_0}(n)$ by~\eqref{shel1},
and on the other hand,
$\delta< \delta_n\Rightarrow \gamma_{\delta_n,i}< \delta$, as
$\delta\in S$, contradiction. Consequently, $\gamma_{\delta,i}=\delta$
for all $\delta\in S$.

We know that for all $\beta\in \mu$ and $n\in \omega$ the set
$\{ \delta\in \mu\mid  \cf \delta=\omega \logic\land \delta_n=\beta\}$
is stationary, so there is $\delta\in S$, where
$\cf \delta=\omega$, such that $\delta_n=\beta$.

As before, we can show that
$$
f(a_0^{\delta_n})=
\tau_\delta(a_{\alpha_{\delta,1}}^{\delta_n},\dots,
a_{\alpha_{\delta,k_\delta}}^{\delta_n})=
\tau_{\delta_n}(a_{\alpha_{\delta_n,1}}^{\gamma_{\delta_n,1}},\dots,
a_{\alpha_{\delta_n,k_{\delta_n}}}^{\gamma_{\delta_n,k_{\delta_n}}}),
$$
where the linear combination $\tau_{\delta_n}$ is reduced.
Therefore $\gamma_{\delta_n,i}\in \{ \delta_n\}$ for all~$i$,
so $\gamma_{\delta_n,i}=\delta_n$,
i.e., $\gamma_{\beta,i}=\delta_n$ for all~$i$. As the linear combination
$\tau_\beta(a_{\alpha_{\beta,1}}^{\gamma_{\beta,1}},\dots)$
is reduced, we have that the ordinals
$\alpha_{\beta,i}$, where $1\le i\le k_\beta$, are distinct.

As $f|_{B_3}$ commutes with~$f_3^*$, for every~$\beta$
$$
\tau_\beta(a_{\alpha_{\beta,1}}^0,\dots, a_{\alpha_{\beta,k_\beta}}^0)=
\tau_0(a_{\alpha_{0,1}}^0,\dots, a_{\alpha_{\alpha_0,k_0}}^0).
$$
As the $\alpha_{\beta,i}$ are distinct, necessarily
$$
\{ \alpha_{\beta,i}\mid 1\le i\le k_\beta\}=
\{ \alpha_{0,i}\mid 1\le i\le k_0\},
$$
but as $\alpha_{\beta,i}$ is increasing with~$i$
(for each~$\beta$, by the choice of $\alpha_{\beta,i}$),
necessarily $\alpha_{\beta,i}=\alpha_{0,i}$ and
$k_\beta=k_0$, i.e., $\tau=\tau_0$.

Therefore $f(a_0^{\beta,n})=\tau_0(a_{\alpha_{0,1}}^{\beta,n},\dots)$.

The other direction of this assertion is immediate.
\end{proof}

Let $\varphi_3(f)$ say that
\begin{enumerate}
\item
$\Rng f|_{B_2}\subseteq B_3$;
\item
$\exists f_1\in \Endom(A')\ (f_1|_{B_2}=f|_{B_2}\logic\land \varphi_2(f_1))$.
\end{enumerate}

The formula $\varphi_3(f)$ holds if and only if
$$
f(a_0^\beta)=\tau(a_{\gamma_1}^\beta,\dots,a_{\gamma_k}^\beta)
$$
for any $\beta\in \mu$ and some $\tau, \gamma_1,\dots,\gamma_k$.
This follows immediately from Statement~\ref{shel_s4_3}.

Let the formula $\varphi_4(f)$ say that
\begin{enumerate}
\item
$\Rng f|_{B_2}\subseteq B_3$;
\item
$\varphi_3(f)$;
\item
$\forall g\in \Endom(A')\
(\varphi_3(g)\Rightarrow g\circ f_5\circ f|_{\tilde a_0^0}=
f_5\circ f\circ f_5\circ g|_{\tilde a_0^0}
\logic\land
f_5^*\circ f\circ f_5^*\circ f|_{\tilde a_0^0}=
f_6^*\circ f|_{\tilde a_0^0}
\logic\land
f_2^*\circ f|_{\tilde a_0^0}=f_2^*|_{\tilde a_0^0}$.
\end{enumerate}

As in case~I, we can check that the formula $\varphi_4(f)$ holds in
$\Endom(A')$ if and only if
$$
f(a_0^\beta)=\tau(a_{\gamma_1}^\beta,\dots, a_{\gamma_k}^\beta),
$$
where $\tau$ is a~beautiful linear combination.

\smallskip

Case III. $\mu$ is a~singular cardinal number.

Let $\mu_1< \mu$, where $\mu_1$ is a~regular cardinal number
and $\mu_1> \omega$.
Let
$I^*\setminus J=I_0\cup
\{ \langle \alpha,\delta,n\rangle\mid \alpha\in \mu,\allowbreak\
{\delta\in \mu_1},\allowbreak\ \cf \delta=\omega,\ n\in \omega\}$,
$|I_0|=\mu$, $a_\alpha^{\beta,n}=a_{\langle \alpha,\beta,n\rangle}$.

For every limit ordinal $\delta\in \mu_1$ such that $\cf \delta=\omega$,
similarly to the previous case, choose an increasing sequence
$(\delta_n)_{n\in\omega}$ of ordinal numbers less than $\delta$,
with limit~$\delta$, such that
for any $\beta\in \mu_1$ and $n\in \omega$ the set
$\{ \delta\in \mu_1\mid \beta=\delta_n\}$
is a~stationary subset in~$\mu_1$.

Let
\begin{align*}
B_1&=\Cl\{ a_{\alpha}^0\mid \alpha\in \mu\},\\*
B_2&= \Cl\{ a_0^\beta\mid \beta\in \mu_1\},\\*
B_3&=\Cl\{ a_\alpha^\beta\mid \alpha\in \mu,\ \beta\in \mu_1\}.
\end{align*}

As in case~II, we can define the functions~$f_i^*$ in such a~way that
for some $\varphi^0(\dots)$ the formula $\varphi^0(f;\dots)$ holds in
$\Endom(A')$ if and only if
there exist a~linear combination~$\tau$ and distinct ordinal numbers
${\alpha_1,\dots,\alpha_n\in \mu}$ such that for every $\beta\in \mu_1$
$$
f(a_0^\beta)=\tau(a_{\alpha_1}^\beta,\dots,a_{\alpha_n}^\beta).
$$

Let the formula $\varphi^1(f)$ say that
\begin{enumerate}
\item
$\Rng f|_{\tilde a_0^0}\subseteq B_2$;
\item
for every $g\in \Endom(A')$
$\varphi^0(g)\Rightarrow
(f\circ g)|_{\tilde a_0^0}=(g\circ f)|_{\tilde a_0^0}$.
\end{enumerate}

It is easy to check that the formula $\varphi^1(f)$ holds if and only if
there exist a~linear combination~$\sigma$ and distinct ordinal numbers
$\beta_1,\dots,\beta_m\in \mu_1$ such that for any $\alpha\in \mu$
$$
f(a_\alpha^0)=\sigma(a_\alpha^{\beta_1},\dots,a_\alpha^{\beta_m}).
$$

As the cardinal number $\mu_1$ is regular, we can use case~II.
Thus, there is a~formula $\varphi^2(f;\dots)$ such that $\varphi^2(f)$
holds in $\Endom(A')$ if and only if there exist
a~beautiful linear combination~$\sigma$ and distinct $\beta_i\in \mu_1$ such
that for any $\alpha\in \mu$
$$
f(a_\alpha^0)=\sigma(a_\alpha^{\beta_1},\dots,a_\alpha^{\beta_m}).
$$

Let $\mu=\bigcup\limits_{i\in \cf \mu} \mu_i$, where $\mu_i\in \mu$ and
the sequence $(\mu_i)$
increases. We just prove that for every $\gamma\in \cf \mu$
there is a~function $\bar f_\gamma^*$ such that
\begin{enumerate}
\item
the formula $\varphi^2[f,\bar f_\gamma^*]$ holds in $\Endom(A')$
if and only if
there exist a~beautiful linear combination~$\sigma$
and distinct $\beta_i\in \mu_\gamma^+$
such that for all $\alpha\in \mu$
$$
f(a_\alpha^0)=\sigma(a_\alpha^{\beta_1},\dots,a_\alpha^{\beta_m});
$$
\item
${f_\gamma^*}_0$ is a~projection onto
$$
\Cl\{ a_\alpha^\beta\mid \alpha\in \mu,\ \beta\in \mu_\gamma^+\}.
$$
\end{enumerate}

Choose $\mu_1=\mu_\gamma^+$ and consider the same $\varphi^2(f)$ as
in the case~II, with ${f_\gamma^*}_0$ being a~projection onto
$\Cl\{ a_\alpha^\beta\mid {\alpha\in \mu},\allowbreak\ \beta\in \mu_1\}$.

Let $\tau$ and $\sigma$ be beautiful linear combinations,
$\beta_1,\dots,\beta_m \in \mu_{\gamma_k}^+$, $k=1,\dots,n$,
and for every $\alpha\in \mu$
$$
f(a_\alpha^0)=\sigma(a_\alpha^{\beta_1},\dots,a_\alpha^{\beta_n}).
$$
We show that in this case we have
$$
\varphi^2(f, \tau(\bar f_{\gamma_1}^*,\dots,\bar f_{\gamma_n}^*)).
$$
The formula
$\varphi^2(f,\tau(\bar f_{\gamma_1}^*,\dots, \bar f_{\gamma_n}^*))$
holds if and only if
$f(a_\alpha^0)=\sigma(a_\alpha^{\beta_1},\dots,a_\alpha^{\beta_m})$
for all $\alpha\in \mu$ and distinct
$\beta_1,\dots,\beta_n\in \mu_{\gamma_1}$, which is true.
Further, we have that
$\tau({f_{\gamma_1}^*}_0,\dots,{f_{\gamma_n}^*}_0)$
is a~projection onto the set
$\Cl \{ a_\alpha^\beta\mid {\alpha\in \mu},\allowbreak\
\beta\in \mu_{\gamma_1}^+\}$
because
$$
\tau(x_1,\dots,x_n)=x_s.
$$

Recall how we proved Theorem~\ref{shel_t4_1} from Lemma~\ref{shel_l4_2}.
This proof easily implies that there exist
a~formula $\varphi^3$ and a~vector of functions $g^*$ such that the formula
$\varphi^3(\bar f,\bar g^*)$ holds
if and only if
$\bar f =\bar f_\gamma^*$ for some $\gamma\in \mu$.

Let now the formula $\varphi^4(f,\bar g^*)$ say that there exists
$\bar f_1$ such that $\varphi^3(\bar f_1,\bar g^*)$ and for
every~$\bar f_2$ satisfying the formulas
$\varphi^3(\bar f_2,\bar g^*)$ and
$\Rng (\bar f_1)_0 \subseteq \Rng (\bar f_2)_0$,
also $\varphi^2(f,\bar f_2)$ holds.
If the formula $\varphi^4(f,g^*)$ holds, then there exists
$\bar f_1=\bar f_{\gamma}^*$ for some $\gamma\in \mu$,
and for every $\bar f_2=\bar f_\lambda^*$ (where $\lambda\ge \gamma$)
we have the formula
$$
f(a_\alpha^0)=\sigma(a_\alpha^{\beta_1},\dots,a_\alpha^{\beta_m}),
$$
where $\beta_1,\dots,\beta_m< \mu_{\lambda}^+$.

Let $f$ be such that
$$
f(a_\alpha^0)=\sigma(a_\alpha^{\beta_1},\dots,a_\alpha^{\beta_m}),\quad
\beta_1,\dots,\beta_m\in \mu.
$$
Then $\beta_1,\dots,\beta_m\in \mu_\gamma^+$ for $\gamma\in \cf \mu$ and
therefore the formula $\varphi^4(f,g^*)$ holds for some~$g^*$.

Now we only need to consider the formula
$\varphi^4(f_5^*\circ f\circ f_5^*)$, which is the required formula.
\end{proof}

\section{Formulation of the Main Theorem,
Converse Theorems, Different Cases}
\subsection{Second Order Language of Abelian Groups}
As we mentioned above, we shall consider second order models of Abelian
groups, i.e., consider the second order group language, where the
3-place symbol will denote not multiplication, but addition
(i.e., we shall write $x_1=x_2+x_3$ instead of $P^3(x_1,x_2,x_3)$).

As we see, formulas $\varphi({\dots})$ of the language $\mathcal L_2$
consist of the following subformulas:
\begin{enumerate}
\item
$\forall x$ ($\exists x$);
\item
$x_1=x_2$ and $x_1=x_2+x_3$,
where every variable $x_1$, $x_2$, and~$x_3$ either is a~free variable
of the formula~$\varphi$ or is defined in the formula~$\varphi$
with the help of the subformulas $\forall x_i$ or $\exists x_i$, $i=1,2,3$;
\item
$\forall P(v_1,\dots,v_n)$ ($\exists P(v_1,\dots,v_n)$), $n> 0$;
\item
$P(x_1,\dots,x_n)$, where every variable $x_1,\dots,x_n$,
and also every ``predicate'' variable $P(v_1,\dots,v_n)$ either is a~free
variable of the formula~$\varphi$,
or is defined in this formula with the help of the subformula
$\forall x_i$, $\exists x_i$,
$\forall P(v_1,\dots,v_n)$, $\exists P(v_1,\dots,v_n)$.
\end{enumerate}

Equivalence of two Abelian groups $A_1$~and~$A_2$
in the language~$\mathcal L_2$ will be denoted by
$$
A_1\equiv_{\mathcal L_2} A_2,\ \ \text{or}\ \ A_1\equiv_2 A_2.
$$

As we remember,
the \emph{theory} of a~model~$\mathcal U$ of a~language~$\mathcal L$
is the set of all sentences of the language~$\mathcal L$ which are true
in this model.
In some cases we shall consider together with theories
$\Th_2(A)=\Th_{\mathcal L_2}(A)$
also theories $\Th_2^\varkappa (A)$, which contain
those sentences~$\varphi$ of the language~$\mathcal L_2$
that are true for arbitrary sequence
$$
\langle a_1,\dots, a_q,b_1^{l_1},\dots,b_s^{l_s}\rangle,
$$
where $a_1,\dots,a_q\in A$, $b_i^{l_i}\subset A^{l_i}$,
and $|b_i^{l_i}|\le \varkappa$.
If
$\varkappa\ge |A|$, then $\Th_2(A)$ and $\Th_2^\varkappa(A)$ coincide.

\subsection{Formulation of the Main Theorem}\label{sec4_2}
If $A=D\oplus G$, where the group~$D$ is divisible and the group~$G$
is reduced, then the \emph{expressible rank of the group~$A$}
is the cardinal number
$$
r_{\mathrm{exp}}=\mu=\max(\mu_D,\mu_G),
$$
where $\mu_D$ is the rank of~$D$ and $\mu_G$ is the rank of the
basic subgroup of~$G$.

We want to prove the following theorem.

\begin{theorem}\label{t4.1}
For any infinite $p$-groups $A_1$~and~$A_2$
elementary equivalence of endomorphism rings $\Endom(A_1)$
and $\Endom(A_2)$ implies coincidence of the second order theories
$\Th_2^{r_{\mathrm{exp}}(A_1)}(A_1)$ and
$\Th_2^{r_{\mathrm{exp}}(A_2)}(A_2)$ of the groups $A_1$~and~$A_2$,
bounded by the cardinal numbers
$r_{\mathrm{exp}}(A_1)$ and $r_{\mathrm{exp}}(A_2)$, respectively.
\end{theorem}

In Secs.~\ref{sec5}--\ref{sec7}, we shall separately prove this
theorem for Abelian groups $A_1$~and~$A_2$ with various properties
and in Sec.~\ref{sec8} gather them and prove the main theorem.

Note that if the group~$A$ is finite,
then the ring $\Endom(A)$ is also finite.
Since in the case of finite models elementary equivalence
(and also equivalence in the language~$\mathcal L_2$)
is equivalent to isomorphism of them, then in the case where
one of the groups $A_1$~and~$A_2$
is finite Theorem~\ref{t4.1} follows from Theorem~\ref{t2.8}.
Therefore we now suppose that
the groups $A_1$~and~$A_2$ are infinite.

\subsection{Proofs of ``Converse'' Theorems}
Let us prove two theorems which are, in some sense, converse
to our main theorem.

\begin{theorem}\label{t4.2}
For any Abelian groups $A_1$~and~$A_2$, if the groups $A_1$~and~$A_2$
are equivalent in the second order logic~$\mathcal L_2$,
then the rings $\Endom(A_1)$ and $\Endom(A_2)$
are elementarily equivalent.
\end{theorem}

\begin{proof}
Every 2-place predicate variable $P(v_1,v_2)$ will be called
a~\emph{correspondence on the group~$A$}.
A~correspondence $P(v_1,v_2)$ on the group~$A$ will be
called a~\emph{function on the group~$A$} (notation:
$\Func(P(v_1,v_2))$, or simply $\Func(P)$)
if it satisfies the condition
$$
(\forall x\losp \exists y\losp P(x,y)) \logic\land
(\forall x\losp \forall y_1\losp \forall y_2\losp P(x,y_1) \logic\land
P(x,y_2)\Rightarrow y_1=y_2).
$$
A~function $P(v_1,v_2)$ will be called an
\emph{endomorphism on the group~$A$}
(notation: $\rEndom(P(v_1,v_2))$, or simply $\rEndom(P)$)
if it satisfies the additional condition
$$
\forall x_1\losp \forall x_2\losp \forall y_1\losp \forall y_2\losp
P(x_1,y_1) \land P(x_2,y_2)\Rightarrow P(x_1+x_2,y_1+y_2).
$$

Now consider an arbitrary sentence~$\varphi$
of the first order ring language.
This sentence can contain the subformulas
\begin{enumerate}
\item
$\forall x$;
\item
$\exists x$;
\item
$x_1=x_2$;
\item
$x_1=x_2+x_3$;
\item
$x_1=x_2\cdot x_3$.
\end{enumerate}

Let us translate this sentence to a~sentence~$\tilde \varphi$
of the second order group language by the following algorithm:
\begin{enumerate}
\item
the subformula $\forall x\losp ({\dots})$ is translated to the subformula
$$
\forall P^x(v_1,v_2)\losp (\rEndom(P^x)\Rightarrow {\dots});
$$
\item
the subformula $\exists x\losp ({\dots})$ is translated to the subformula
$$
\exists P^x(v_1,v_2)\losp (\rEndom(P^x)\land {\dots});
$$
\item
the subformula $x_1=x_2$ is translated to the subformula
$$
\forall y_1\losp \forall y_2\losp
(P^{x_1}(y_1,y_2) \Leftrightarrow P^{x_2}(y_1,y_2));
$$
\item
the subformula $x_1=x_2+x_3$ is translated to the subformula
$$
\forall y\losp \forall z_1\losp \forall z_2\losp
\forall z_3\losp
(P^{x_2}(y,z_2)\land P^{x_3}(y,z_3) \Rightarrow
(P^{x_1}(y,z_1)\Leftrightarrow z_1=z_2+z_3));
$$
\item
the subformula $x_1=x_2\cdot x_3$ is translated to the subformula
$$
\forall y\losp \forall z\losp (P^{x_1}(y,z) \Rightarrow
\exists t\losp (P^{x_2}(y,t)\land P^{x_3}(t,z)).
$$
\end{enumerate}

We need to show that a~sentence~$\varphi$ holds in the model $\Endom(A)$
if and only if
the sentence~$\tilde \varphi$ holds in the model~$A$.

If $A$ is a~model of an Abelian group, then the model $\Endom(A)$
consists of sets of pairs of elements of the model~$A$,
$x=\{ \langle u_1,u_2\rangle \mid u_1,u_2\in A\}$, with the conditions
\begin{enumerate}
\item
$\forall u_1\losp \exists u_2\losp \langle u_1,u_2\rangle\in x$;
\item
$\forall u_1\losp \forall u_2\losp \forall u_3\losp
(\langle u_1,u_2\rangle \in x \logic\land
\langle u_1,u_3\rangle \in x\Rightarrow u_2=u_3)$;
\item
$\forall u_1\losp \forall u_2\losp \forall u_3\losp \forall u_4\losp
(\langle u_1,u_3\rangle\in x \logic\land
\langle u_2,u_4\rangle \in x \Rightarrow
\langle u_1+u_2,u_3+u_4\rangle \in x)$.
\end{enumerate}

Therefore a~sequence $a_1,\dots,a_q$ for which the formula~$\varphi$
is satisfied in the model $\Endom(A)$ is a~sequence
consisting of sets of pairs of elements from~$A$
satisfying the conditions (1)--(3).

Let us establish the identity bijection between elements of $\Endom(A)$
and the corresponding sets of pairs
of the model~$A$. Let an element~$a_i$ of the model $\Endom(A)$
correspond to a~set $A_i\subset A\times A$.

1. If the formula~$\varphi$ has the form $x_i=x_j$, then $\varphi$
holds for a~sequence $a_1,\dots,a_q$ if and only if $a_i=a_j$,
i.e., $a_i$~and~$a_j$ are equal endomorphisms of the model $\Endom(A)$,
and the sets $A_i$~and~$A_j$ consist of the same elements, i.e.,
in the model~$A$ for the sequence $A_1,\dots,A_q$ the formula
$$
\forall y_1\losp \forall y_2\losp
(P^{x_i}(y_1,y_2)\Leftrightarrow P^{x_j}(y_1,y_2))
$$
is true.

2. If the formula $\varphi$ has the form $x_i=x_j+x_k$,
then $\varphi$ holds for a~sequence
$a_1,\dots,a_q$ if and only if $a_i=a_j+a_k$,
i.e., an endomorphism~$a_i$ is the sum of endomorphisms $a_j$~and~$a_k$,
and this means that
in the model~$A$ for every element $b\in A$ and for every $b_1,b_2,b_3\in A$
such that $\langle b,b_1\rangle\in A_i$, $\langle b,b_2\rangle\in A_j$,
$\langle b,b_3\rangle\in A_k$, we have $b_1=b_2+b_3$
(i.e., formally speaking, $\langle b_1,b_2,b_3\rangle\in I(Q_1^3)$).
It is equivalent to $A\vDash \tilde \varphi$.

3. If the formula~$\varphi$ has the form $x_i=x_j\cdot x_k$,
then the formula~$\varphi$ is true for a~sequence $a_1,\dots,a_q$
if and only if
$a_i=a_j\cdot a_k$,
i.e., the endomorphism~$a_i$ is a~composition of endomorphisms
$a_j$~and~$a_k$, and this means that
in the model~$A$ for every $b_1\in A$ and for every $b_2\in A$
such that $\langle b_1,b_2\rangle\in A_i$,
there exists $b_3\in A$ such that
$\langle b_1,b_3\rangle\in A_j$ and $\langle b_3,b_2\rangle \in A_k$.
This is equivalent to $A\vDash \tilde \varphi$.

4. If $\varphi$ has the form $\theta_1\land \theta_2$,
$\theta_1$~and~$\theta_2$
are true in the model $\Endom(A)$ for a~sequence $a_1,\dots,a_q$
if and only if
$\tilde \theta_1$~and~$\tilde \theta_2$
are true in the model~$A$ for the sequence $A_1,\dots,A_q$,
then it is clear that
the formula~$\varphi$ is true in the model $\Endom(A)$ for a~sequence
$a_1,\dots,a_q$ if and only if the formula~$\tilde \varphi$
is true in the model~$A$ for the sequence $A_1,\dots,A_q$,
because
$$
\widetilde{\theta_1\land \theta_2}=\tilde \theta_1\land \tilde \theta_2.
$$

5. A~similar case is the formula~$\varphi$ having the form $\neg \theta$,
because
$$
\widetilde{\neg \theta}=\neg \tilde \theta.
$$

6. Finally, suppose that the formula~$\varphi$
has the form $\forall x_i\losp \psi$.
The formula~$\varphi$ is true in the model $\Endom(A)$
for a~sequence $a_1,\dots,a_q$ if and only if the formula~$\psi$
is true in the model $\Endom(A)$ for the sequence
$a_1,\dots,a_{i-1},a,a_{i+1},\dots,a_q$ for any $a\in \Endom(A)$,
i.e., the formula~$\tilde \psi$ is true in the model~$A$ for
the sequence $A_1,\dots,A_{i-1}$, $\bar A,A_{i+1},\dots,A_q$
for every set $\bar A\subset A\times A$ which is an
endomorphism of the ring~$A$, i.e., which satisfies the formula $\rEndom$.
Therefore the formula~$\varphi$ is true in the model $\Endom(A)$
for a~sequence $a_1,\dots,a_q$ if and only if the formula
$$
\widetilde{\forall x_i\losp \psi} \prisv
\forall P^{x_i}(v_1,v_2)\losp (\rEndom(P^{x_i})\Rightarrow \tilde \psi)
$$
is true for the sequence $A_1,\dots,A_q$ in the model~$A$.

Suppose now that Abelian groups $A_1$~and~$A_2$ are equivalent in the
language~$\mathcal L_2$. Consider an arbitrary sentence~$\varphi$
of the first order ring language, which is true in $\Endom(A_1)$. Then
the sentence~$\tilde \varphi$ is true in the group~$A_1$
and hence in the group~$A_2$.
Consequently, the sentence~$\varphi$ is true in the ring $\Endom(A_2)$.
Therefore the rings $\Endom(A_1)$ and $\Endom(A_2)$
are elementarily equivalent.
\end{proof}

For the next theorem we need some formulas.

1. The formula
$$
\Gr(P(v))\prisv \forall a\losp \forall b\losp
(P(a)\land P(b) \Rightarrow \exists c\losp (c=a+b\logic\land P(c))
\logic\land P(0) \logic\land
\forall a\losp (P(a)\Rightarrow \exists b\losp (b=-a \logic\land P(b)))
$$
holds for those sets $\{ a\in A\mid P(a)\}$ that are subgroups in~$A$,
and only for them.

2. The formula
$$
\mathrm{Cycl}(P(v))\prisv \Gr(P(v)) \logic\land
\exists a\losp
(P(a) \logic\land
\forall P_a(v)\losp
(\Gr(P_a(v)) \land P_a(a)\Rightarrow
\forall b\losp (P(b)\Rightarrow P_a(b)))
$$
characterizes cyclic subgroups in~$A$.

3. The formula
\begin{multline*}
\mathrm{DCycl}(P(v))\prisv
\Gr(P(v)) \logic\land
\forall a\losp
(P(a)\Rightarrow \exists P_1(v)\losp \exists P_2(v)\losp
(P_1(a) \logic\land \mathrm{Cycl}(P_1(v))\\
{}\logic\land \forall b\losp \neg (P_1(b) \land P_2(b)) \logic\land
\forall b\losp (P(b)\Rightarrow
\exists b_1\losp \exists b_2\losp
(P_1(b_1) \logic\land P_2(b_2) \logic\land b=b_1+b_2)))
\end{multline*}
characterizes those subgroups in~$A$
that are direct sums of cyclic subgroups.

4. For every $a,a_1,a_2\in A$ the formula
$$
\Gr_a(P_a(v)) \prisv P_a(a) \logic\land \Gr(P_a(v)) \logic\land
\forall P(v)\losp
(P(a)\land \Gr(P(v))\Rightarrow \forall b\losp (P_a(b)\Rightarrow P(b)))
$$
defines in~$A$ the subgroup $\{ b\in A\mid P_a(b)\}$
of all powers (exponents) of the element~$a$;
the formula
\begin{align*}
& (o(a_1)\le o(a_2))\\
& \quad \prisv \exists P_1(v)\losp \exists P_2(v)\losp
\exists P(v_1,v_2)\losp
(\Gr_{a_1}(P_1) \logic\land \Gr_{a_2}(P_2) \logic\land
\forall b_1\losp
(P_1(b_1)\Rightarrow \exists b_2\losp (P_2(b_2)\land P(b_1,b_2)))\\
& \quad \logic\land \forall b_1\losp \forall b_2\losp
\forall c_1\losp \forall c_2\losp
(P_1(b_1) \logic\land
P_1(c_1) \logic\land b_1\ne c_1 \logic\land
P_2(b_2) \logic\land P_2(c_2) \logic\land P(b_1,b_2) \logic\land
P(c_1,c_2)\Rightarrow b_2\ne c_2))
\end{align*}
holds if and only if the order of~$a_1$
is not greater than the order of~$a_2$;
the formula
$$
(o(a_1)=o(a_2)) \prisv (o(a_1)\le o(a_2)) \logic\land (o(a_2)\le o(a_1))
$$
shows that the orders of $a_1$~and~$a_2$ coincide;
the formula
$$
(o(a_1)< o(a_2))\prisv (o(a_1)\le o(a_2)) \logic\land
\neg (o(a_2)\le o(a_1))
$$
shows that the order of~$a_1$
is strictly smaller than the order of~$a_2$.

5. For every $a\in A$ the formula
$$
\mathrm{GOrd}_a(P(v))\prisv
\Gr(P) \logic\land \forall b\losp (P(b)\Rightarrow o(b)\le o(a))
$$
holds for those subgroups that are bounded by the order of~$a$,
and only for them.

6. The formula
\begin{align*}
& \mathrm{Mult}_a(x,b) \prisv
\exists P(v)\losp \exists P_{x,b}(v_1,v_2)\losp
(
\mathrm{Cycl}(P) \logic\land P(x) \logic\land P(b) \logic\land{}
\forall b_1\losp
(
P(b_1)\Rightarrow \exists b_2\losp (P(b_2) \land P_{x,b}(b_1,b_2))
\\
& \quad {}\logic\land
(
\forall b_1\losp \forall b_2\losp \forall b_3\losp
P(b_1) \land P_{x,b}(b_1,b_2) \land P_{x,b}(b_1,b_3) \Rightarrow b_2=b_3
)
\\
& \quad {}\logic\land
(
\forall b_1\losp \forall b_2\losp \forall b_3\losp
\forall c_1\losp \forall c_2\losp \forall c_3\losp
P(b_1) \logic\land P(b_2) \logic\land P(b_3)
\\
& \quad {}\logic\land b_3=b_1+b_2 \logic\land c_3=c_1+c_2 \logic\land
P_{x,b}(b_1,c_1) \logic\land P_{x,b}(b_2,c_2) \Rightarrow P_{x,b}(b_3,c_3)
)
\\
& \quad {}\logic\land P_{x,b}(x,0) \logic\land \forall y\losp
(
P(y) \land py=x \Rightarrow \neg P_{x,b}(y,0)
)
\logic\land
\exists c\losp
(
P(b,c) \land o(c)=o(a)
)
))
\end{align*}
holds for elements $x$ and~$b$ with the property
$x=o(a)\cdot b$, and only for them.

7. The formula
$$
\mathrm{Serv}(P(v)) \prisv \Gr(P) \logic\land
\forall a\losp \forall x\losp
(P(x)\Rightarrow\\
\Rightarrow \exists b\losp
(\mathrm{Mult}_a(x,b)\Rightarrow
\exists c\losp (P(c)\land \mathrm{Mult}_a(x,c))))
$$
holds for pure subgroups of the group~$A$, and only for them.

8. The formula
$$
\mathrm{FD}(P(v))\prisv
\Gr(P) \logic\land \forall a\losp\exists b\losp
\exists x_1\losp\exists x_2\losp
(P(x_1) \logic\land P(x_2) \logic\land a+x_1=p(b+x_2))
$$
holds for subgroups $G=\{ x\mid P(x)\}$ such that $A/G$
is a~divisible subgroup, and only for them.

9. The formula
$$
\Base(P(v)) \prisv \Gr(P) \land \mathrm{DCycl}(P) \land
\mathrm{Serv}(P) \land \mathrm{FD}(P)
$$
defines basic subgroups in~$A$.

It is clear that if we have some subgroup~$G'$ of the group~$G$,
then we similarly can write the formula $\Base_{G'}(P)$
which holds for basic subgroups of the group~$G'$, and only for them.

10. The formula
$$
\mathrm{D}(P(v)) \prisv \Gr(P) \logic\land
\forall a\losp ( P(a) \Rightarrow \exists b\losp ( P(b) \logic\land a=pb ))
$$
defines divisible subgroups in~$A$.

11. The sentence
\begin{align*}
& \mathrm{Exept} \prisv
\forall P\losp (\Gr(P) \Rightarrow \neg ( \mathrm{D}(P) ))
\\
& \quad \logic\land \forall P(v)\losp
(
\Base(P)\Rightarrow
\neg
(
\exists F(v_1,v_2)\losp
(
\forall a\losp
(
P(a)\Rightarrow \exists b\losp (F(a,b))
\\
& \quad \logic\land \forall b\losp \exists a\losp (P(a)\land F(a,b))
\logic\land \forall a\losp \forall b\losp (F(a,b) \Rightarrow P(a))
\\
& \quad \logic\land \forall a_1\losp \forall a_2\losp
\forall b_1\losp \forall b_2\losp
(
a_1\ne a_2 \logic\land F(a_1,b_1) \logic\land F(a_2,b_2)
\Rightarrow b_1\ne b_2
)
\\*
& \quad \logic\land \forall b_1\losp \forall b_2\losp
\forall a_1\losp \forall a_2\losp
(
b_1\ne b_2 \logic\land F(a_1,b_1) \logic\land F(a_2,b_2)
\Rightarrow a_1\ne a_2
)
))))
\end{align*}
is true for reduced $p$-groups such that their basic subgroups
have smaller power (and therefore are countable), and only for them.
Thus if $B_1$ is a~basic subgroup of the group~$A_1$, $B_2$~is a~basic
subgroup of~$A_2$, $\varkappa_1=|B_1|$, and $\varkappa_2=|B_2|$, then
$$
\Th_2^{\varkappa_1}(A_1)=\Th_2^{\varkappa_2}(A_2)
$$
implies that either the groups $A_1$~and~$A_2$ are reduced, their
basic subgroups are countable, and they themselves are uncountable,
or this is not true for both of the groups $A_1$~and~$A_2$.

In the first case $\varkappa_1=\varkappa_2=\omega$.

\begin{theorem}\label{t4.2d}
If Abelian groups $A_1$~and~$A_2$ are reduced
and their basic subgroups are countable,
then $\Th_2^\omega(A_1)=\Th_2^\omega(A_2)$
implies elementary equivalence of
the rings $\Endom(A_1)$ and $\Endom(A_2)$.
\end{theorem}

\begin{proof}
We know (see Theorem~\ref{4.endom}) that for a~reduced $p$-group~$A$
the action of any endomorphism $\varphi\in \Endom(A)$ is completely defined
by its action on a~basic subgroup~$B$. Furthermore, let
$A'\subset A$ and let $B$~be also a~basic subgroup of~$A'$. Then any
$\varphi\colon A'\to A$ is also completely defined by its action
on~$B$. Indeed, if $\varphi_1,\varphi_2\colon A'\to A$ and
$\varphi_1(b)=\varphi_2(b)$ for all $b\in B$, then for
$\varphi \prisv \varphi_1-\varphi_2\colon A'\to A$
we have $\varphi(b)=0$ for all $b\in B$. Hence $\varphi$ induces
a~homomorphism $\tilde \varphi\colon A'/B\to A$. But the group
$A'/B$ is divisible and the group~$A$ is reduced,
i.e., $\tilde \varphi= 0$. Consequently, $\varphi=0$.

Note that for every element $a\in A$ there exists a~countable
subgroup $A'\subset A$ containing~$a$ and the group~$B$ as a~basic subgroup.

Indeed, consider a~quasibasis of the group~$A$ having the form
$$
\{ a_i,c_{j,n}\}_{i\in \omega,\ j\in \varkappa,\ n\in \omega},
$$
where $\{ a_i\}$ is a~basis of~$B$, $p c_{j,1}=0$,
$pc_{j,n+1}=c_{j,n}-b_{j,n}$,
$b_{j,n}\in B$, $o(b_{j,n})\le p^n$, and $o(c_{j,n})=p^n$.

As we remember, every element $a\in A$ can be written in the form
$$
a=s_1 a_{i_1}+\dots+s_m a_{i_m}+t_1 c_{j_1,n_1}+\dots+t_r c_{j_r,n_r},
$$
where $s_i$~and~$t_j$ are integers, none of~$t_j$
is divisible by~$p$ and the indices
$i_1,\dots,i_m$, $j_1,\dots,j_r$ are all distinct.
Further, this form is unique in the sense that
the members $sa_i$ and $tc_{j,n}$ are uniquely defined.

Consider a~decomposition of our element~$a$ and the subgroup in~$A$
generated by the group~$B$ and all~$c_{k,n}$, where
$n\in \omega$ and $k\in \{ j_1,\dots,j_r\}$.
This group~$A'$ is countable, it contains~$a$, and $B\subset A'$ is
its basic subgroup.

Let now a~predicate $B(v)$ satisfy in~$A$ the formula $\Base(B)$, i.e.,
$B(v)$ defines in~$A$ a~basic subgroup $B=\{ x\mid B(x)\}$.

A~correspondence $P(v_1,v_2)$ is called a~\emph{homomorphism
of the group~$B$ into the group~$A$} (notation: $\Hom_B(P)$) if
\begin{multline*}
\forall x\losp
(B(x)\Leftrightarrow \exists y\losp (P(x,y))) \logic\land
\forall x\losp \forall y_1\losp \forall y_2\losp
(P(x,y_1)\land P(x,y_2) \Rightarrow y_1=y_2)
\\
\logic\land
\forall x_1\losp \forall x_2\losp \forall y_1\losp \forall y_2\losp
(P(x_1,y_1) \land P(x_2,y_2) \Rightarrow P(x_1+x_2,y_2+y_2)).
\end{multline*}
It is clear that such a~predicate $P(v_1,v_2)$ can be used in sentences
from $\Th_2^\omega(A)$, because the group~$B$ is countable.

Consider some $B(v)$ such that
the formula $\Base(B)$ is true, a~predicate $\Phi(v_1,v_2)$
such that $\Hom_B(\Phi)$, and $a\in A$.

We shall write $b=\Phi(a)$ if
\begin{enumerate}
\item
$B(a)\land \Phi(a,b)$ or
\item
$\neg B(a) \logic\land
\forall G(v)\losp
(
\Gr(G) \logic\land G(a) \logic\land \forall x\losp (G(x)\Rightarrow B(x))
\logic\land \Base_G(B) \Rightarrow \exists \varphi (v_1,v_2)\losp
(
\forall x\losp
(
G(x)\Leftrightarrow \exists y\losp (\varphi(x,y))
\logic\land
\forall x\losp \forall y_1\losp \forall y_2\losp
(\varphi(x,y_1)\land \varphi(x,y_2)\Rightarrow y_1=y_2)
\logic\land
\forall x_1\losp \forall x_2\losp \forall y_1\losp \forall y_2\losp
(\varphi(x_1,y_1)\land \varphi(x_2,y_2)\Rightarrow
\varphi(x_1+x_2,y_1+y_2))
\logic\land
\forall x\losp \forall y\losp (\Phi(x,y) \Rightarrow \varphi(x,y))
\logic\land \varphi(a,b)
)))
$.
\end{enumerate}

It is clear that for every $a\in A$ there exist no more than one
$b\in A$ such that
$b=\Phi(a)$, and if the homomorphism $\Phi\colon B\to A$ can be extended
to an endomorphism $A\to A$, then it necessarily exists.

Now we shall consider $\Phi(v_1,v_2)$ such that
\begin{multline*}
\rEndom_B(\Phi) \prisv
\Hom_B(\Phi) \logic\land
\forall a\losp \exists b\losp (b=\Phi(a)) \\
{}\logic\land \forall a_1\losp \forall a_2\losp
\forall b_1\losp \forall b_2\losp
(b_1=\Phi(a_1) \logic\land b_2=\Phi(a_2) \Rightarrow b_1+b_2=\Phi(a_1+a_2)).
\end{multline*}
In our case these $\Phi(v_1,v_2)$ define endomorphisms from $\Endom(A)$.

Let us show an algorithm of translation of formulas in this case.

A~sentence~$\varphi$ is translated to the sentence
$$
\tilde \varphi=\exists B(v)\losp (\Base(B)\land \varphi'(B)),
$$
where the formula~$\varphi'$ is obtained from the sentence~$\varphi$
in the following way:
\begin{enumerate}
\item
the subformula $\forall x\losp ({\dots})$ is translated to the
subformula
$$
\forall \Phi^x(v_1,v_2)\losp (\rEndom_B(\Phi^x)\Rightarrow \ldots);
$$
\item
the subformula $\exists x\losp ({\dots})$ is translated to the subformula
$$
\exists \Phi^x(v_1,v_2)\losp (\rEndom_B(\Phi^x)\land \ldots);
$$
\item
the subformula $x_1=x_2$ is translated to the subformula
$$
\forall y_1\losp \forall y_2\losp
(\Phi^{x_1}(y_1,y_2) \Leftrightarrow \Phi^{x_2}(y_1,y_2));
$$
\item
the subformula $x_1=x_2+x_3$ is translated to the subformula
$$
\forall y\losp \forall z_1\losp \forall z_2\losp \forall z_3\losp
(\Phi^{x_2}(y,z_2)\land \Phi^{x_3}(y,z_3) \Rightarrow
(\Phi^{x_1}(y,z_1)\Leftrightarrow z_1=z_2+z_3));
$$
\item
the subformula $x_1=x_2\cdot x_3$ is translated to the subformula
$$
\forall y\losp \forall z\losp
(\Phi^{x_1}(y,z) \Rightarrow \exists t\losp
(\Phi^{x_2}(y,t)\land z=\Phi^{x_3}(t)).
$$
\end{enumerate}

Now the proof is similar to the proof of the previous theorem.
\end{proof}

\subsection{Different Cases of the Problem}\label{sec4_4}
Following Theorem~\ref{t2.8}, we divide the class of all Abelian
$p$-groups into the following three subclasses:
\begin{enumerate}
\item
bounded $p$-groups;
\item
the groups $D\oplus G$, where $D$~is a~nonzero
divisible group and $G$~is a~bounded group;
\item
groups with unbounded basic subgroups.
\end{enumerate}

Now we shall show how to find sentences which distinguish
groups from different subclasses.

If a~group $A$ is bounded, then there exists a~natural number $n=p^k$
such that
$$
\forall a\losp (na=0)
$$
in the group~$A$. This means that in the ring $\Endom(A)$ we also have
$$
\forall x\losp (nx=0).
$$
But if a~group~$A$ is unbounded, then for every~$n$
this sentence is false.
Therefore, we distinguish the case~(1) from all other cases.
Moreover, we now can find the supremum of orders of elements from~$A$.
We shall denote this sentence $\forall x\losp (nx=0)$ by~$\varphi_n$.

Now consider the sentence
\begin{align*}
& \psi_n \prisv \exists \rho_1\losp \exists \rho_2\losp
(
\rho_1\rho_2=\rho_2\rho_1=0 \logic\land
\rho_1^2=\rho_1 \logic\land
\rho_2^2=\rho_2 \logic\land
\rho_1+\rho_2=1 \\
& \quad \logic\land
\forall x\losp (n\cdot \rho_2x\rho_2 =0) \logic\land
\forall \rho\losp \forall \rho'\losp
(
\rho^2=\rho \logic\land
{\rho'}^2=\rho' \logic\land
\rho\rho'=\rho'\rho=0 \logic\land
\rho+\rho'=\rho_1
\\
& \quad \logic\land
\forall \tau_1\losp \forall \tau_2\losp
(\tau_1^2=\tau_1 \logic\land \tau_2^2=\tau_2 \logic\land
\tau_1\tau_2=\tau_2\tau_1=0 \Rightarrow \tau_1+\tau_2\ne \rho)
\Rightarrow
\forall x\losp (\rho x\rho =0 \logic\lor p(\rho x \rho)\ne 0)
)).
\end{align*}

We shall explain by words what this sentence means.

1. There exist orthogonal projections $\rho_1$~and~$\rho_2$;
their sum is~$1$ in the ring $\Endom(A)$. This means that
$A=\rho_1A\oplus \rho_2 A$, $\rho_1 \Endom(A)\rho_1=\Endom(\rho_1 A)$, and
$\rho_2 \Endom(A)\rho_2=\Endom(\rho_2 A)$.

2. The condition $\forall x\losp (n\cdot \rho_2 x \rho_2=0)$ means that
in the ring $\Endom(\rho_2 A)$ all elements are bounded
by a~number $n=p^k$, i.e., the group $\rho_2 A$ is bounded.

3. The last part of the sentence~$\psi_n$ states that if in the ring
$\Endom(\rho_1A)$ we consider a~primitive idempotent~$\rho$ (i.e.,
a~projection onto an indecomposable direct summand
$\rho A=\rho \rho_1 A$), then this direct summand does not have any
endomorphisms of order~$p$. Therefore in the group $\rho_1 A$
there are no cyclic direct summands, and so it is divisible.

Consequently, a~ring $\Endom(A)$ satisfies the sentence~$\psi_n$
if and only if
$A=D\oplus G$, where the group~$D$ is divisible and the group~$G$
is bounded by the number~$n$.

Hence for any two groups $A_1$~and~$A_2$ from different classes there
exists a~sentence which distinguishes
the rings $\Endom(A_1)$ and $\Endom(A_2)$.

Thus we can now assume that if rings $\Endom(A_1)$
and $\Endom(A_2)$ are elementarily equivalent,
then the groups $A_1$~and~$A_2$
belong to one subclass, and if both of them belong
to the first or the second subclass, then their reduced parts
are bounded by the same number $n=p^k$, which is supposed to be fixed.

\section[Bounded $p$-Groups]{Bounded $\boldsymbol{p}$-Groups}\label{sec5}
\subsection{Separating Idempotents}\label{sec5_1}
As we have seen above (see Sec.~\ref{sec4_4}), the property of
$\rho \in \Endom(A)$ to be a~decomposable idempotent
which is a~direct summand
is a~first order property.
Let us denote the formula expressing this property
by $\Idem^*(\rho)$, while the formula expressing the property of
$\rho\in \Endom(A)$ to be simply an idempotent
(not necessarily indecomposable) will be denoted by $\Idem(\rho)$.

We consider the group $A=\sum\limits_{i=1}^k A_i$, where
$A_i\cong \bigoplus\limits_{\mu_i} \mathbb Z (p^i)$.
Since the group~$A$ is infinite, we have
that $\mu_l=\max\limits_{i=1,\dots,k} \mu_i$ is infinite and coincides
with~$|A|$.

Consider for every $i=1,\dots,k$ the formula
$$
\Idem_i^*(\rho)=\Idem^*(\rho)\land p^{i-1} \rho\ne 0\land p^i \cdot \rho=0.
$$
For every~$i$ this formula is true for projections on direct summands
of the group~$A$ which are isomorphic to $\mathbb Z(p^i)$,
and only for them.

Now consider the following formula:
\begin{align*}
& \mathrm{Comp}(\rho_1,\dots,\rho_k)=
(\rho_1+\dots+\rho_k=1) \logic\land
\biggl(\,\bigwedge_{i\ne j} \rho_i \rho_j=\rho_j\rho_i=0\biggr)
\\
& \quad \logic\land
\biggl(\,\bigwedge_{i=1}^k \rho_i^2=\rho_i\biggr)
\logic\land
\biggl(\,\bigwedge_{i=1}^k p^i \rho_i=0\biggr) \logic\land
\biggl(\,\bigwedge_{i=1}^k p^{i-1} \rho_i\ne 0\biggr)
\\
& \quad \logic\land
\biggl(\,\bigwedge_{i=1}^k \forall \rho\losp
(
\Idem^*(\rho) \logic\land
\exists \rho'\losp (\rho+\rho'=\rho_i \logic\land
\rho\rho'=\rho' \rho=0 \logic\land \Idem(\rho'))
\Rightarrow \Idem^*_i(\rho))
\biggr).
\end{align*}
We see that the group $A$ (which corresponds to this formula)
is decomposed into a~direct sum
${\rho_1A\oplus \rho_2 A\oplus\dots\oplus \rho_k A=A}$,
and in every subgroup
$\rho_iA$ all indecomposable direct summands have the order~$p_i$.
Therefore $\rho_1A\oplus \dots\oplus \rho_k A$ is a~decomposition
of~$A$, isomorphic to the decomposition $\sum\limits_{i=1}^k A_i$.
{\sloppy

}

Let us assume that the projections $\rho_1,\dots,\rho_k$ from the formula
$\mathrm{Comp}({\dots})$ are fixed.
To separate them from other idempotents we shall
denote them by $\bar \rho_1,\dots,\bar \rho_k$.
Having fixed idempotents
$\bar \rho_1,\dots,\bar \rho_k$ of the ring $\Endom(A)$,
we have also its fixed subrings
$$
\bar E_i=\bar \rho_i \Endom(A)\bar \rho_i,
$$
each of which is isomorphic to the ring
$$
\Endom(\rho_iA)\cong \Endom(A_i).
$$

Given an idempotent~$\rho$
satisfying the formula
$\Idem_i^*(\rho)$ ($\Idem_i(\rho)$),
the formula expressing for~$\rho$
the fact that it is a~direct
summand in the group $\bar \rho_i A$
(i.e., $\exists \rho'\losp
(\Idem(\rho') \logic\land \rho\rho'=\rho' \rho=0 \logic\land
\rho+\rho'=\bar \rho_i)$)
will be written as
$\overline{\Idem_i^*}(\rho)$
($\overline{\Idem_i}(\rho)$).
This formula means that the subgroup
$\rho A$ is a~direct summand in the group $\bar \rho_i A=A_i$.

The number~$l$ from the set $\{ 1,\dots,k\}$ which satisfies the sentence
$$
\mathrm{Card}_l=\bigwedge_{i=1,\,i\ne l}^k
\exists a\losp \forall \rho\losp
(\overline{\Idem_i^*}(\rho) \Rightarrow \rho a\bar \rho_l\ne 0)
$$
is the number of the group~$A_l$ with $|A_l|=|A|=\mu$, because this sentence
means that there exists an endomorphism $a\in \Endom(A)$
mapping~$A_l$ to~$A$
in such a~way that on every direct summand of~$A_i$ it is nonzero.
This means that $|A_l|\ge |A_i|$, i.e., the power $|A_l|$ is maximal.
Let us assume that also the number~$l$ is fixed.

The formula $\mathrm{Card}_l$ shows that we can write formulas
that determine for every two projections $\rho_1$~and~$\rho_2$
whether $|\rho_1A|< |\rho_2A|$,
or $|\rho_1A|> |\rho_2A|$, or $|\rho_1A|=|\rho_2A|$ is true.
Let us denote these formulas by $|\rho_1|< |\rho_2|$,
$|\rho_1|> |\rho_2|$, and $|\rho_1|=|\rho_2|$, respectively.

The formula
$$
\Fin(\rho) \prisv \forall \rho_1\losp \forall \rho_2\losp
(\Idem(\rho_1) \logic\land \Idem(\rho_2) \logic\land \Idem(\rho)
\logic\land \rho_1=\rho+\rho_2 \logic\land \rho\rho_2=\rho_2\rho=0
\Rightarrow |\rho|< |\rho_1|)
$$
means that the group $\rho A$ is finitely generated.
Respectively, the formula
$$
\mathrm{Inf}(\rho)\prisv \Idem(\rho) \logic\land \neg \Fin(\rho)
$$
holds for projections on infinitely generated groups.

The formula
$$
\mathrm{Count}(\rho) \prisv
\mathrm{Inf}(\rho) \logic\land
\forall \rho_1 \losp (\mathrm{Inf}(\rho_1)\Rightarrow |\rho|\le |\rho_1|)
$$
is true for projections on countably generated groups and only for them.

Finally, we need the formula
$$
\overline{\Idem_l^\omega}(\rho)=
\overline{\Idem_l}(\rho)\land \mathrm{Count}(\rho),
$$
which means that the group $\rho A$ is a~countably generated direct
summand of the group~$A_l$.

\subsection{Special Sets}\label{sec5_2}
At first we shall formulate what special sets we want to have.
We must obtain two sets. One of them must contain
$\mu_i$~independent indecomposable projections on direct summands
of~$A_i$, for every $i=1,\dots ,k$, the other set must contain
$\mu=\mu_l$ projections on independent countably generated direct summands
of the group~$A_l$.

By Theorem~\ref{shel_t4_1}, we see that there exists a~formula
$\varphi(\bar g;f)$ satisfying the following condition.
If $\{ f_i\}_{i\in \mu}$ is a~set of elements from $\Endom(A')$,
then there exists a~vector~$\bar g$ such that the formula
$\varphi(\bar g;f)$ is true in $\Endom(A')$ if and only if $f=f_i$
for some $i\in \mu$. We fix this formula~$\varphi$.

Suppose that we have some fixed $i\in \{ 1,\dots,k\}$.
We have already shown that from the ring $\Endom(A)$ we can transfer
to the ring $\Endom(A_i)$.
Suppose that we argue in the ring $\Endom(A_i)$
(which satisfies the conditions of Theorem~\ref{shel_t4_1}).
In this ring let us consider the following formula:
\begin{align*}
& \tilde \varphi_i(\bar g) \prisv
\forall f'\losp
(\varphi(\bar g,f')\Rightarrow \overline{\Idem_i^*}(f'))
\\
& \quad \logic\land \forall f'\losp
(\overline{\Idem_i}(f') \logic\land \forall f_1\losp
(\varphi(\bar g,f_1)
\Rightarrow
\exists f_2\losp
(\overline{\Idem_i}(f_2) \logic\land f_1f_2=f_2f_1=0 \logic\land
f_1+f_2=f')
)\Rightarrow |f'|=|\rho_i|)
\\
& \quad \logic\land \forall f'\losp
(
\varphi(\bar g,f')\Rightarrow
(
\exists f\losp
(
\overline{\Idem_i}(f)
\logic\land
\forall f_1\losp
(\varphi(\bar g,f_1) \logic\land f_1\ne f'\Rightarrow f_1 f=ff_1=f_1)
\logic\land ff'=f'f=0
))).
\end{align*}

The part $\forall f'\losp
(\varphi(\bar g,f')\Rightarrow \overline{\Idem_i^*}(f'))$
means that the vector~$\bar g$
is such that the formula
$\varphi(\bar g,f)$ is true only for projections~$f$
on indecomposable direct summands of the group~$A_i$.

The part $\forall f'\losp (\overline{\Idem_i}(f') \logic\land
\forall f_1\losp (\varphi(\bar g,f_1)\Rightarrow
\exists f_2\losp (\overline{\Idem_i}(f_2) \logic\land
f_1f_2=f_2f_1=0 \logic\land f_1+f_2=f'))\Rightarrow {|f'|=|\rho_i|})$
means that those subgroups of the group~$A_i$
that
contain all summands~$fA$ satisfying $\varphi (\bar g,f)$
have the same power as~$A_i$,
i.e., this part means that
the power of the set of these~$f$ is equal to~$\mu$.

The last part of the formula means that for every~$f'$ satisfying
$\varphi(\bar g,f')$,
the group generated by all other~$f$ satisfying $\varphi(\bar g,f)$
does not intersect with~$f'$,
i.e., the set of all $f$ satisfying $\varphi(\bar g,f)$ is independent.

This set will be denoted by~$\mathbf F_i$.
It consists of~$\mu_i$ independent projections
on indecomposable direct summands of the group~$A_i$. Naturally, this set
can be obtained for every vector~$\bar g_i$ satisfying the formula
$\tilde \varphi_i(\bar g_i)$,
therefore we have to write not~$\mathbf F_i$, but $\mathbf F_i(\bar g_i)$,
and we shall do so in what follows.
But in the cases where parameters are not so important we shall
omit them.

The union of all~$\mathbf F_i$ for $i=1,\dots,k$ will be denoted
by~$\mathbf F$. The set~$\mathbf F$ depends on the parameter
$\bar g=(\bar g_1,\dots,\bar g_k)$.

Now we need to obtain a~set~$\mathbf F'$ consisting of $\mu_l=\mu$
independent projections on countably generated
direct summands of~$A_l$.
It will be done similarly to the previous case, we only need to change
in the formula $\tilde \varphi_l(\bar g')$ the following parts:
$\overline{\Idem_l^*}(f)$ to $\overline{\Idem_l^\omega}(f)$;
besides, we shall consider vectors~$\bar g'$ such that
\begin{enumerate}
\item
$\forall f \in \mathbf F_l\losp
(\exists f'\losp (\varphi(\bar g',f')\logic\land ff'=f))$,
i.e., for every cyclic direct summand $fA$ (where $f\in \mathbf F_l$)
of~$A_l$ there exists a~countably generated
summand $f'A$ of~$A_l$ such that
$\varphi (\bar g',f')$ and $fA\subset f'A$;
\item
(we shall write it by words, because we do not want to write complicated
formulas) every direct summand in~$A_l$
which contains all $fA$ for all projections~$f$
such that $\varphi(\bar g_l,f)$ contains all~$f'A$ such that
$\varphi(\bar g',f')$.
\end{enumerate}

Denote the corresponding formula by
$\Tilde{\Tilde \varphi}_l(\bar g')$
and the obtained set of projections by
$\mathbf F'=\mathbf F'(\bar g')$.

\subsection[Interpretation of the Group $A$ for Every Element
$\protect\mathbf F'$]%
{Interpretation of the Group $\boldsymbol{A}$
for Every Element $\boldsymbol{\mathbf F'}$}\label{sec5_3}
By interpretation of the group~$A$ for every element from~$\mathbf F'$
we understand the following. We have $\mu$~independent
direct summands ${\mathcal F_i =f_iA}$ ($i\in \mu$) each of which
is a~direct sum of a~countable set of cyclic groups of order~$p^l$.
Every endomorphism of the group~$A$ acts independently on every
summand~$\mathcal F_i$, hence if for
every endomorphism $\varphi\in \Endom(A)$ we can map every element
of~$\mathcal F_i$ to some element of~$A$,
then we shall be able to map every endomorphism
$\varphi \in \Endom A$ to a~set of $\mu$ elements
of~$A$. This is what we need below
to obtain the second order theory of the group~$A$.
So in this section, we shall concentrate on a~bijective correspondence
between some homomorphisms from the group~$\mathcal F_i$
into the group~$A$, and elements of the group~$A$, and introduce
on the set of homomorphisms
an operation~$\oplus$
that
under this bijection
corresponds to the addition of the group~$A$.

Let us fix some projection $g\in \mathbf F'$.
Consider the set $\Endom_g$ of all those homomorphisms $h\colon gA\to A$
that satisfy the following conditions:
\begin{enumerate}
\item
$\forall f\in \mathbf F_l\losp
(fg=f\Rightarrow (hf=0 \logic\lor
\exists f'\in \mathbf F\losp (hf=f'hf\ne 0)))$;
this means that for every projection~$f$ from our special set~$\mathbf F_l$,
if the projection maps $A$ to the indecomposable direct summand $fA$
of the module~$gA$,
then either $h(fA)=0$ or $h(fA)\subset f'A$ for some
projection $f'\in \mathbf F$;
\item
$\exists f\losp (\Fin(f)\land \Idem(f)\land fh=h)$; this means that
the image of the subgroup~$gA$
under the endomorphism~$h$ is finitely generated;
\item
$\bigwedge\limits_{i=1}^k \forall f\in \mathbf F_i\losp
\neg \Bigl(\exists f_1\dots \exists f_{p^i}\losp
\Bigl(\,\bigwedge\limits_{q\ne s} f_q\ne f_s \logic\land
f_1,\dots,f_{p^i} \in F_l \logic\land f_1g=g_1 \logic\land
\dots \logic\land f_{p^i}g=f_{p^i} \logic\land hf_1=fhf_1\ne 0
\logic\land \dots \logic\land hf_{p^i}=fhf_{p^i}\ne 0\Bigr)\Bigr)$;
this means that for every $i=1,\dots,k$ the inverse image of each
$fA\subset A_i$, where $f\in \mathbf F_i$, can not contain more
than $p^i-1$ different elements $f_mA$, where $f_mA\subset gA$
and $f_m\in \mathbf F_l$.
\end{enumerate}

Two elements $h_1$~and~$h_2$ from the set $\Endom_g$
are said to be equivalent ($h_1\sim h_2$) if they satisfy
the following formula:
\begin{multline*}
\exists f_1\losp \exists f_2\losp
((gf_1g)\cdot (gf_2g)=(gf_2g)\cdot (gf_1g)=g \\
\logic\land \forall f\in \mathbf F_l\losp
(fg=f\Rightarrow \forall f'\in \mathbf F\losp
(h_1 f=f' h_1f\ne 0\Leftrightarrow (gf_1gh_2)f=f'(gf_1gh_2)f\ne 0)).
\end{multline*}

This means that there exists an automorphism $gf_1g$ of the group~$gA$
which maps $h_2$ to an endomorphism $(gf_1g\cdot h_2)$
such that for every $\rho\in F_l$, where $\Image \rho =gA$, both
endomorphisms $h_1$ and $gf_1gh_2$ map this subgroup either to zero or
to the same $f'A$ ($f'\in \mathbf F$).
The obtained set $\Endom_g/{\sim}$ will be denoted by $\tEnd_g$.
We can interpret elements of this set as
finite sets of projections from~$\mathbf F$ with the condition
that every projection from~$\mathbf F_i$ can belong
to this set at most $p^i-1$ times.
Respectively, every element of the set $\tEnd_g$
can be interpreted as a~set of pairs, where the first element in a~pair
is a~projection~$f$ from~$\mathbf F$ and the second element is an integer
from~$0$ to~$p^i-1$, where $i$ is such that $f\in \mathbf F_i$,
and almost all (all except for a~finite number)
second components of the pairs are equal to~$0$.
Now we can construct a~bijective mapping between the set
$\tEnd_g$ and the group~$A$,
where the image of the described set
$\{ \langle f_j,l_j\rangle\mid j\in J\}$ is the element
$\sum\limits_{j\in J} l_j\xi_j=a\in A$, where $\xi_j$ is some fixed
generator of the cyclic group~$f_jA$.

Now we only need to introduce addition on the set $\tEnd_g$
to make the obtained bijective mapping an isomorphism of Abelian groups.

We shall introduce addition by the formula ($h_1,h_2,h_3\in \tEnd_g$)
\begin{align*}
& (h_3=h_1\oplus h_2) \prisv
\bigwedge_{i=1}^k \forall f\in \mathbf F_i\losp
\biggl(\,\bigwedge_{j=0}^{p^i-1}
\exists g_1\dots \exists g_j\in \mathbf F_l\losp
\bigwedge_{q\ne s} (g_q\ne g_s \logic\land g_qg=g_q
\\
& \quad {}\logic\land h_3 g_q=fh_3g_q\ne 0) \logic\land
\neg \biggl(\exists g_1\dots \exists g_{j+1}\in \mathbf F_l\losp
\bigwedge_{q\ne s} (g_q\ne g_s \logic \land g_qg=g_q \logic\land
h_3g_q=fh_3g_q\ne 0 )\biggr)
\\
& \quad {}\Rightarrow
\biggl(\,\bigvee_{m=0}^j
\exists g_1\dots \exists g_m\in \mathbf F_l\losp
\bigwedge_{q\ne s}
(g_q\ne g_s \logic\land g_qg=g_q \logic\land h_1g_q=fh_1g_q\ne 0)
\\
& \quad {}\logic\land
\neg \biggl(\exists g_1\dots \exists g_{m+1}\in \mathbf F_l\losp
\bigwedge_{q\ne s}
(g_q\ne g_s \logic\land g_qg=g_q \logic\land h_1g_q= fh_1g_q\ne 0)\biggr)
\\
& \quad {}\logic\land
\exists g_1\dots \exists g_{\gamma(j,m)}\in \mathbf F_l\losp
\bigwedge_{q\ne s}
(g_q\ne g_s \logic\land g_qg=g_q \logic\land h_2 g_q=fh_2g_q\ne 0)
\\
& \quad {}\logic\land
\neg \biggl(\exists g_1\dots \exists g_{\gamma(j,m)+1}\in \mathbf F_l\losp
\bigwedge_{q\ne s}
( g_q\ne g_s \logic\land g_qg=g_q \logic\land
h_2 g_q=fh_2g_q\ne 0)\biggr)\!\biggr)\!\biggr),
\end{align*}
where $\gamma(j,m)=j-m$ if $j\ge m$, and $\gamma(j,m)=p^i+j-m$ if $j< m$.

Now we see that for every $g\in \mathbf F'$ we have a~definable set
$\tEnd_g$ with the addition operation~$\oplus$, which is isomorphic
to the group~$A$.

\subsection{Proof of the First Case in the Theorem}\label{sec5_4}
\begin{proposition}\label{p5.1}
For any two infinite Abelian
$p$-groups $A_1$~and~$A_2$ bounded by the number~$p^k$,
elementary equivalence of the rings $\Endom(A_1)$ and
$\Endom(A_2)$ implies equivalence of the groups
$A_1$~and~$A_2$ in the language~$\mathcal L_2$.
\end{proposition}

\begin{proof}
For every $\tilde g\in \mathbf F'$ by $\mathrm{Resp}_{\tilde g}(h)$
we shall denote the following formula:
$$
\mathrm{Resp}_{\tilde g}(h) \prisv
\forall g\in \mathbf F'\losp \exists h'\losp
((\tilde g hg)(gh'\tilde g)= \tilde g
\logic\land
(gh'\tilde g)(\tilde ghg)=g).
$$
This formula means that an endomorphism~$h$
isomorphically maps every summand $gA$ (where $g\in \mathbf F'$)
to the summand $\tilde gA$.

As above, let us consider an arbitrary sentence~$\varphi$
in the second order group language and show an algorithm
translating this sentence~$\psi$
to a~sentence~$\tilde \psi$ of the first order ring language
so that $\Endom(A)\vDash \tilde \psi$ if and only if $A\vDash \varphi$.

Let us translate the sentence $\psi$ to the sentence
$$
\exists \bar g_1\dots \exists \bar g_k\losp
(
\tilde \varphi_1(\bar g_1) \logic\land \dots \logic\land
\tilde \varphi_k(\bar g_k)
\\
{}\logic\land \exists \bar g'\losp
(
\Tilde {\Tilde \varphi}_l (\bar g',\bar g_l) \logic\land
\exists \tilde g\in \mathbf F'(\bar g')\losp
\exists h \losp
(\mathrm{Resp}_{\tilde g}(h) \logic\land
\psi'(\bar g_1,\dots,\bar g_k,\bar g',\tilde g,h))
)),
$$
where the formula $\psi'({\dots})$ is obtained from the sentence~$\psi$
with the help of the following translations of subformulas of~$\psi$:
\begin{enumerate}
\item
the subformula $\forall x$ is translated to the subformula
$\forall x\in \tEnd_{\tilde g}$;
\item
the subformula $\exists x$ is translated to the subformula
$\exists x\in \tEnd_{\tilde g}$;
\item
the subformula $\forall P_m(v_1,\dots,v_m)\losp ({\dots})$
is translated to the subformula
$$
\forall f_1^P\dots \forall f_m^P\losp
\biggl(\forall g\in \mathbf F'(\bar g')\losp
\biggl(\,\bigwedge_{i=1}^m (f_i^Pg\in \Endom_g)\biggr)\Rightarrow
\ldots\biggr);
$$
\item
the subformula $\exists P_m(v_1,\dots,v_m)\losp ({\dots})$
is translated to the subformula
$$
\exists f_1^P\dots \exists f_m^P\losp
\biggl(
\forall g\in \mathbf F'(\bar g')\losp
\biggl(\,\bigwedge_{i=1}^m (f_i^Pg\in \Endom_g)\biggr)\logic\land
\ldots\biggr);
$$
\item
the subformula $x_1=x_2$ is translated to the subformula $x_1\sim x_2$;
\item
the subformula $x_1=x_2+x_3$ is translated to the subformula
$x_1\sim x_2\oplus x_3$;
\item
the subformula $P_m(x_1,\dots,x_m)$
is translated to the subformula
$$
\exists g\in \mathbf F'(\bar g')\losp
\biggl(\, \bigwedge_{i=1}^m f_i^Pg=x_i hg\biggr).
$$
\end{enumerate}

We can explain by words what these translations mean. According to
existence of the set~$\mathbf F'$, we have
$\mu$~groups $\tEnd_g$ for $g\in \mathbf F'$, each of which
is isomorphic to the group~$A$. We fix one chosen element
$\tilde g\in \mathbf F'$, and therefore we
fix one group $\tEnd_{\tilde g}$, isomorphic to~$A$.
Naturally, all subformulas $\forall x$, $\exists x$,
$x_1=x_2$, $x_1=x_2+x_3$
(of first order logic) will be translated to the corresponding
subformulas for the group $\tEnd_{\tilde g}$. Now we need
to interpret an arbitrary relation $P_m(v_1,\dots,v_m)$
on~$A$ in the ring $\Endom(A)$.
Such a~relation is some subset in $A^m$, i.e.,
a~set of ordered $m$-tuples of elements from~$A$.
There are at most~$\mu$ such $m$-tuples,
therefore the set $P_m(v_1,\dots,v_m)$ can be
considered as a~set of~$\mu$
$m$-tuples of elements from~$A$ (some of them can coincide).
We consider $m$~endomorphisms $f_1^P,\dots,f_m^P\in \Endom(A)$ such that
the restriction of each of them on any $gA$ (where $g\in \mathbf F'$)
is an element of $\widetilde{\Endom}_g$. Thus
for every $g\in \mathbf F'$ the restriction of the endomorphisms
$f_1^P,\dots,f_m^P$ on~$gA$ is an $m$-tuple of elements
of the group $\tEnd_{\tilde g}$~($\cong A$),
where an isomorphism between $\tEnd_{\tilde g}$ and
$\tEnd_g$ is given by the fixed mapping~$h$
which isomorphically maps every module $gA$ to $\tilde gA$.

So we can see that the sentence~$\psi$ is true in~$A$ if and only if
the sentence~$\tilde \psi$ is true in the ring $\Endom(A)$.
Therefore, as in the previous section, we have the proof.
\end{proof}

\section[Direct Sums of Divisible and
Bounded $p$-Groups]{Direct Sums of Divisible and
Bounded $\boldsymbol{p}$-Groups}\label{sec6}
\subsection{Finitely Generated Groups}
Every infinite finitely generated Abelian $p$-group~$A$
has the form $D\oplus G$,
where $D$ is a~divisible finitely generated group and $G$~is
a~finite group. There is no need to prove the following proposition.

\begin{proposition}\label{p6.1}
If Abelian $p$-groups $A_1$~and~$A_2$ are finitely generated, then
elementary equivalence of their endomorphism rings
$\Endom(A_1)$ and $\Endom(A_2)$
implies that the groups $A_1$~and~$A_2$ are isomorphic.
\end{proposition}

\subsection{Infinitely Generated Divisible Groups}\label{sec6_2}
As in Sec.~\ref{sec5_1}, the formula $\Idem^*(\rho)$
will denote the property
of an endomorphism~$\rho$ to be an indecomposable idempotent, while
$\Idem(\rho)\prisv (\rho^2=\rho)$. If in a~divisible group~$D$ we have
$\Idem^*(\rho)$, then $\rho A\cong \mathbb Z(p^\infty)$.

Note that, despite the fact that in Sec.~\ref{sec3} we considered
direct sums of cyclic groups of the same order, the Shelah theorem
remains true also for divisible groups,
because a~divisible group is a~union of the groups
$$
\bigoplus_\mu \mathbb Z(p)\subset
\bigoplus_\mu \mathbb Z(p^2)\subset \dots \subset
\bigoplus_\mu \mathbb Z(p^n)\subset \ldots
$$
(proof of an even more general case is given later in Sec.~\ref{sec7}).
Therefore, similarly to Sec.~\ref{sec5_2}, we have a~definable set
$\mathbf F=\mathbf F(\bar g)$
consisting of $\mu$~indecomposable projections on linearly
independent direct summands of~$D$,
and also a~definable set $\mathbf F'=\mathbf F'(\bar g')$ consisting of
$\mu$~projections on linearly independent countably generated
direct summands of~$D$.

Let us fix some element $g\in \mathbf F'$ and construct
(as in Sec.~\ref{sec5}) an interpretation of the group~$D$ for this set.

Namely, let us consider the set $\Endom_g$
of all homomorphisms $h\colon gA\to A$
satisfying the following conditions:
\begin{enumerate}
\item
$\forall f\in \mathbf F\losp
(fg=f\Rightarrow
(hf=0 \logic\lor \exists \tilde f\in \mathbf F\losp
(\tilde hf=\tilde fhf\ne 0)))$; this means that for every
projection~$f$ from~$\mathbf F$ such that
$fA\subset gA$, we have either $h(fA)=0$ or $h(fA)\subset \tilde fA$
for some $\tilde f\in \mathbf F$;
\item
$\exists f\losp (\Fin(f) \logic\land \Idem(f) \logic\land fh=h)$;
this means that the image of~$gA$ on the endomorphism~$h$
is finitely generated;
\item
$\forall f\in \mathbf F\losp
(\exists f'\in \mathbf F\losp
(f'g=g' \logic\land hf'=fhf'\ne 0)\Rightarrow
\exists \tilde f\in \mathbf F\losp
(\tilde fg=\tilde f \logic\land h\tilde f=fh\tilde f\ne 0)
\logic\land
\forall {f'\in \mathbf F}\allowbreak\losp
(f'g=f' \logic\land hf'=fhf'\ne 0
\Rightarrow
\exists \alpha\losp (\alpha hf'=h\tilde f)))$;
this means that for every element $f\in \mathbf F$
either the inverse image of~$fA$ is empty or it
contains an element $\tilde f A\subset gA$
(where $\tilde f\in \mathbf F$) such that
the kernel $\tilde f A$ on the mapping~$h$ has the maximal order
among all kernels~$f'A$ such that $f'\in \mathbf F$ and $f'A\subset gA$.
\setcounter{slast}{\value{enumi}}
\end{enumerate}

Before the last condition we shall introduce some new notation. Let $h$
be some endomorphism and $f_1,f_2\in \mathbf F$.

We shall write $f_1\sim_h f_2$ if and only if the formula
$$
\exists \alpha\losp
(\alpha^2=1 \logic\land \alpha f_1=f_2\alpha \logic\land
\alpha f_2=f_1\alpha \logic\land hf_1=h\alpha f_1 \logic\land
hf_2=h\alpha f_2)
$$
is true.
This formula means that the images of the groups $f_1A$ and $f_2A$
on the endomorphism~$h$ coincide and the kernels of the groups $f_1A$ and
$f_2A$ on this endomorphism are isomorphic.

Now for an endomorphism~$h$ and projections $f_1,f_2\in \mathbf F$ we shall
introduce the formula
$$
\exists \alpha\losp
(\alpha^2=1 \logic\land \alpha f_1=f_2\alpha \logic\land
\alpha f_2=f_1\alpha \logic\land hf_1=ph\alpha f_1).
$$
This formula states that the images of the groups $f_1A$ and $f_2A$
on the endomorphism~$h$ coincide and the kernel of $f_1A$
is $p$~times greater than the kernel of $f_2A$.
This formula will be denoted
by $f_1\sim_h f_2+1$.

Now we shall introduce the last condition:
\begin{enumerate}
\setcounter{enumi}{\value{slast}}
\item
$ \neg \Bigl(\exists f_1,\dots,f_p\in \mathbf F\losp
\Bigl(\Bigl(\,\bigwedge\limits_{i\ne j} f_i\ne f_j\Bigr) \logic\land
f_1g=f_1 \logic\land \dots \logic\land f_pg=f_p \logic\land
hf_1=fhf_1\ne 0 \logic\land \dots \logic\land
hf_p= fhf_p\ne 0 \logic\land
\Bigl(\,\bigwedge\limits_{i\ne j} f_i\sim f_j\Bigr)\Bigr)\Bigr)$;
this means that there exist at most $p-1$ such projections
from~$\mathbf F$ onto subgroups of~$gA$
that their kernels are isomorphic.
\end{enumerate}

Two elements $h_1$~and~$h_2$ of the set $\Endom_g$ are said to be
equivalent ($h_1\sim h_2$)
if the following formula holds:
$$
\exists f_1\losp \exists f_2\losp
((gf_1g)\cdot (gf_2g)=(gf_2g)\cdot (gf_1g)=g
\logic\land \forall f\in \mathbf F\losp
(fg=f\Rightarrow h_1f=(gf_1g)h_2f \logic\land h_2f=(gf_2g)h_1f)).
$$

This means that there exist
mutually inverse automorphisms $gf_1g$ and $gf_2g$
of the group $gA$ which map $h_2$~and~$h_1$ to automorphisms
$(gf_1g) h_2$ and $(gf_2g)h_1$ such that
we have $h_2=(gf_2g)h_1$ and $h_1=(gf_1g)h_2$ on~$gA$
for every projection from~$\mathbf F$
projecting on the subgroup of~$gA$.

The obtained set $\Endom_g/{\sim}$ will be denoted by $\tEnd_g$.

Suppose that we have two quasicyclic groups $C$~and~$C'$, one of which
has generators $c_1,\dots,c_n,\dots$ ($pc_1=0$, $pc_{n+1}=c_n$)
and the other one has
$c_1',\dots,c_n',\dots$ ($pc_1'=0$, $pc_{n+1}'=c_n'$). Consider the set
of homomorphisms $\Hom(C,C')$.
Two homomorphisms $\alpha_1,\alpha_2\in \Hom(C,C')$
correspond to each other
under some automorphism of the group~$C$
if and only if
their kernels are isomorphic, i.e., have the same order.
Thus, all homomorphisms from $\Hom(C,C')$ are divided into
a~countable number of classes, and every class uniquely corresponds
to a~nonnegative integer~$i$ such that $|\Ker \alpha|=p^i$.

Consequently, every class $h\in \tEnd_g$ can be mapped
to a~finite set of finite sequences
$$
\langle f,m(f),l_1(f),\dots,l_{m(f)}(f)\rangle,
$$
where $f\in \mathbf F$, $m(f)\in \mathbb N$, and $l_i(f)=0,\dots,p-1$.
It is clear that endomorphisms from the same equivalence class are
mapped to the same sets, and endomorphisms from different classes
are mapped to different sets.
Moreover,
it is clear that every finite set of sequences is mapped to
some class of endomorphisms.
Now every such set of sequences
will be mapped to an element
$$
\sum_{f\in \mathbf F} l_1(f)c_1(f)+\dots+l_m(f)c_{m}(f),
$$
where $c_1(f),\dots,c_n(f),\dots$ are some fixed generators of~$fA$.

Therefore we have obtained a~bijection between the set
$\tEnd_g$ and the group $A\cong \bigoplus\limits_\mu \mathbb Z(p^\infty)$.

Now let us introduce on the set $\tEnd_g$ an addition
($h_3=h_1\oplus h_2$) in such a~way that this bijection becomes
an isomorphism between Abelian groups.

Let $h_1,h_2,h_3\in \tEnd_g$.
\begin{align*}
& (h_3=h_1\oplus h_2) \prisv
\forall f\in \mathbf F\losp
\biggl(
\forall f'\in \mathbf F\losp
\biggl(
f'g=f' \logic\land h_1f'=fh_1f'\ne 0
\\
& \quad {}\Rightarrow
\bigwedge_{i=0}^{p-1}
\biggl(
\exists f_1\dots \exists f_i\in \mathbf F\losp
\bigwedge_{q\ne s} (f_q\ne f_s \logic\land
f_q\sim_{h_1} f' \logic\land f_q g=f_q \logic\land
h_1f_q=fh_1f_q\ne 0)
\\
& \quad \logic\land
\neg \biggl(\exists f_1,\dots,\exists f_{i+1}\in \mathbf F\losp
\bigwedge_{q\ne s}
(f_q\ne f_s \logic\land f_q \sim_{h_1} f' \logic\land
f_qg=f_q \logic\land h_1f_q=fh_1f_q\ne 0)
\biggr)
\\
& \quad {}\Rightarrow
\bigvee_{j=0}^{p-1}
\biggl(\exists f_1\dots \exists f_j\in \mathbf F\losp
\bigwedge_{q\ne s}
(f_q\ne f_s \logic\land f_q\sim_{h_2} f' \logic\land
f_qg=f_q \logic\land h_2f_q=fh_2f_q\ne 0)
\\
& \quad \logic\land
\neg \biggl(\exists f_1\dots\exists f_{j+1}\in \mathbf F\losp
\bigwedge_{q\ne s}(f_q\ne f_s \logic\land f_q\sim_{h_2} f' \logic\land
f_qg=f_q \logic\land h_2f_q=fh_2f_q\ne 0)\biggr)
\\
& \quad \logic\land
\biggl(\!\biggl(
\exists f_1\dots \exists f_{(i+j) \bmod p}\in \mathbf F\losp
\bigwedge_{q\ne s}
(f_q \ne f_s \logic\land f_q\sim_{h_3} f' \logic\land
f_qg=f_q \logic\land h_3f_q=fh_3f_q\ne 0)
\\
& \quad \logic\land
\neg \biggl(\exists f_1\dots \exists f_{(i+j)\,mod\, p +1}\in
\mathbf F \bigwedge_{q\ne s}(f_q\ne f_s \logic\land f_q\sim_{h_3} f'
\logic\land f_qg=f_q \logic\land h_3f_q=fh_3f_q\ne 0)\biggr)
\\
& \quad \logic\land
\neg \biggl(\exists f_1\dots f_p\in \mathbf F\losp
\bigwedge_{q\ne s}
(f_q\ne f_s \logic\land f_q\sim_{h_3} f'+1
\\
& \quad \logic\land f_qg=f_q \logic\land
(h_1f_q=fh_1f_q \ne 0 \logic\lor h_2f_q=fh_2f_q\ne 0))\biggr)\!\biggr)
\\
& \quad \logic\lor
\biggl(\exists f_1\dots \exists f_{(i+j) \bmod p+1}\in \mathbf F\losp
\bigwedge_{q\ne s}
(f_q \ne f_s \logic\land f_q\sim_{h_3} f' \logic\land
f_qg=f_q \logic\land h_3f_q=fh_3f_q \ne 0)
\\
& \quad \logic\land
\neg
\biggl(\exists f_1\dots \exists f_{(i+j) \bmod p +2}\in \mathbf F\losp
\bigwedge_{q\ne s} (f_q\ne f_s \logic\land f_q\sim_{h_3} f' \logic\land
f_qg=f_q \logic\land h_3f_q=fh_3f_q\ne 0)\biggr)
\\
& \quad \logic\land \exists f_1\dots f_p\in \mathbf F\losp
\bigwedge_{q\ne s}
(f_q\ne f_s \logic\land f_q\sim_{h_3} f'+1 \logic\land f_qg=f_q
\\
& \quad \logic\land
(hf_q=fh_1f_q\ne 0 \logic\lor h_2f_q=fh_2f_q\ne 0))\biggr)\!
\biggr)\!\biggr)\!\biggr)\!\biggr)\!\biggr).
\end{align*}

Hence for every $g\in \mathbf F'$
there is a~definable set $\tEnd_g$ with the addition operation~$\oplus$,
which is isomorphic to~$A$.

\begin{proposition}\label{p6.2}
For two infinitely generated divisible
$p$-groups $A_1$~and~$A_2$,
elementary equivalence of the rings $\Endom(A_1)$ and
$\Endom(A_2)$ implies equivalence of the groups $A_1$~and~$A_2$
in the language~$\mathcal L_2$.
\end{proposition}

\begin{proof}
Since we have obtained an interpretation of the group~$A$
for every $g\in \mathbf F'$, the proof of this proposition
is completely similar to the proof of
Proposition~\ref{p5.1}.
\end{proof}

\subsection[Direct Sums of Divisible $p$-Groups and
Bounded $p$-Groups of Not Greater Power]{Direct Sums
of Divisible $\boldsymbol{p}$-Groups and
Bounded $\boldsymbol{p}$-Groups of Not Greater Power}\label{sec6_3}
In this section, we consider the groups of the form $D\oplus G$, where $D$
is an infinitely generated divisible group, $G$~is a~group bounded by
the number~$p^k$,
and $|G|\le |D|$.
This case is practically the union of the previous two cases.

Namely, let us have idempotents $\rho_D$~and~$\rho_G$
from the formula~$\psi_{p^k}$ of Sec.~\ref{sec4_4}, i.e.,
idempotents which are projections on divisible and bounded parts of~$A$,
respectively, and also idempotents $\rho_1,\dots,\rho_k$,
where ${\rho_1+\dots}+\rho_k=\rho_G$,
which are projections on direct summands of the form
$\smash[b]{\bigoplus\limits_{\mu_1} \mathbb Z(p),\dots,
\bigoplus\limits_{\mu_k}\mathbb Z(p^k)}$, respectively.
Let $|A|=|D|=\mu$. As before, we have the following definable sets:
\begin{enumerate}
\item
$\mathbf F=\mathbf F(\bar g)$ is a~set
of $\mu$~indecomposable projections on linearly independent
direct summands of the group~$D$;
\item
the set $\mathbf F'=\mathbf F'(\bar g')$ consists
of $\mu$~projections
on linearly independent countably generated direct summands of~$D$;
\item
for every $i=1,\dots,k$ the set $\mathbf F_i=\mathbf F_i(\bar g_i)$
consists of $\mu_i$~projections on independent
indecomposable direct summands of $\rho_i A$;
\item
an endomorphism $\varphi\in \Endom(A)$ satisfying the following formula:
\begin{align*}
& \forall f\in \mathbf F(\bar g)\losp (\varphi f\in \mathbf F)
\logic\land \bigwedge_{i=1}^k
(\rho_D\varphi \rho_i= \varphi \rho_i \logic\land
\forall f_i\in \mathbf F_i(\bar g_i)\losp
\exists f\in \mathbf F(\bar g)\losp (\varphi f_i=f\varphi f_i\ne 0))
\\
& \quad \logic\land
\forall f_1,f_2\in \mathbf F(\bar g)\losp
(f_1\ne f_2\Rightarrow \forall f_1',f_2'\in \mathbf F(\bar g)
\losp
(\varphi f_1=f_1'\varphi f_1\ne 0 \logic\land
\varphi f_1=f_2'\varphi f_2\ne 0\Rightarrow f_1'\ne f_2')
\\
& \quad \logic\land
\bigwedge_{i=1}^k \forall f\in \mathbf F(\bar g)\losp
\forall f_i\in \mathbf F_i(\bar g_i)\losp
\forall f'\in F(\bar g)\losp
(f'=\varphi f\Rightarrow f'\varphi f_i=0)
\\
& \quad \logic\land
\bigwedge_{i,j=1}^k \forall f_1,f_2\in \mathbf F(\bar g)\losp
\forall f_i\in \mathbf F_i(\bar g_i)\losp
\forall f_j\in \mathbf F_j(\bar g_j)
\\
& \quad
(f_i\ne f_j \logic\land \varphi f_i=f_1\varphi f_i\ne 0 \logic\land
\varphi f_j=f_2\varphi f_j\ne 0\Rightarrow f_1\ne f_2)).
\end{align*}
\end{enumerate}

We see that such an endomorphism~$\varphi$ embeds the set
$$
\mathbf F(\bar g)\cup \mathbf F_1(\bar g_1)\cup\dots \cup
\mathbf F_k(\bar g_k)
$$
into the set~$\mathbf F(\bar g)$. Therefore for a~given~$\varphi$
we can consider the sets
\begin{align*}
\mathbf F^D&=\mathbf F^D(\bar g, \varphi)=
\{ f\in \mathbf F(\bar g)\mid \exists f'\in \mathbf F(\bar g)\losp
(f\varphi f'=f\varphi=\varphi f'\ne 0)\},\\
\mathbf F^D_1&=\mathbf F^D_1(\bar g_1, \varphi)=
\{ f\in \mathbf F(\bar g)\mid \exists f'\in \mathbf F_1(\bar g_1)\losp
(f\varphi f'=f\varphi=\varphi f'\ne 0)\},\\
& \hspace{4.5cm}\ldots\\
\mathbf F^D_k&=\mathbf F_k^D(\bar g_k, \varphi)=
\{ f\in \mathbf F(\bar g)\mid \exists f'\in \mathbf F_k(\bar g_k)\losp
(f\varphi f'=f\varphi=\varphi f' \ne 0)\}.
\end{align*}
The sets $\mathbf F^D,\mathbf F_1^D,\dots,\mathbf F_k^D$
consist of $\mu,\mu_1,\dots,\mu_k$
projections on indecomposable linearly independent direct summands of
the group~$D$, respectively. We shall write
them in formulas, sometimes omitting parameters, but meaning that
they depend on the parameters $\bar g,\bar g_1,\dots,\bar g_k,\varphi$.

Let us fix some element $g\in \mathbf F'$ and construct an interpretation
of the group $A=D\oplus G$ for this set.

Namely, let us consider the set $\Endom_g$ of all homomorphisms
$h\colon gA\to A$ satisfying the following conditions:
\begin{enumerate}
\item
$\forall f\in \mathbf F\losp
\Bigl(fg=f\Rightarrow
\Bigl(hf=0 \logic\lor
\exists \tilde f\in \mathbf F\losp (hf=\tilde f hf\ne 0) \logic\lor
\bigwedge\limits_{i=1}^k \exists \tilde f\in \mathbf F_i^D\losp
(\tilde f h=\tilde f hf=hf\ne 0)\Bigr)\Bigr)$;
this means that for every projection~$f$ from~$\mathbf F$,
if $fA\subset gA$, then we have
either $h(fA)=0$, or $h(fA)\subset \tilde fA$
for some $\tilde f\in \mathbf F^D$, or $h(fA)=\tilde fA$
for some $\tilde f\in \mathbf F_i^D$;
\item
$\exists f\losp (\Fin(f) \logic\land \Idem(f) \logic\land fh=h)$;
this means that the image of~$gA$ under~$h$ is finitely generated;
\item
$\bigwedge\limits_{i=1}^k \forall f\in \mathbf F_i^D\losp
\neg
\Bigl(\exists f_1,\dots,\exists f_{p^i}\in \mathbf F\losp
\Bigl(\Bigl(\,\bigwedge\limits_{q\ne s} f_q\ne f_s\Bigr) \logic\land
f_1 g=f_1 \logic\land \dots \logic\land f_{p^i} g=f_{p^i} \logic\land
fhf_1=fh=hf_1\ne 0 \logic\land \dots \logic\land
fhf_{p^i}=fh=hf_{p^i}\ne 0\Bigr)\Bigr)$;
this means that for every $i=1,\dots,k$
the inverse image of every $fA\subset D$,
where $f\in \mathbf F_i^D$, contains at most $p^i-1$ distinct
elements $f_mA$ such that $f_mA\subset gA$ and $f_m\in \mathbf F$;
\item
$\forall f\in \mathbf F^D\losp \exists f'\in \mathbf F\losp
(f'g=g' \logic\land hf'=fhf'\ne 0\Rightarrow
\exists \tilde f\in \mathbf F \losp
(\tilde fg=\tilde f \logic\land h\tilde f= fh\tilde f\ne 0)
\logic\land
\forall f'\in \mathbf F\losp (f'g=f' \logic\land hf'=fhf'\ne 0\Rightarrow
\exists \alpha\losp (\alpha hf'=h\tilde f)))$;
this means that for every element $f\in \mathbf F$
either the inverse image $fA$ is empty or it contains
$\tilde f\in \mathbf F$ such that $\tilde f A\subset gA$
and the kernel of $\tilde f A$ under the mapping~$h$ has the maximal
order of all kernels of $f'A$ for $f'\in \mathbf F$
such that $f'A\subset gA$;
\item
$ \neg \Bigl(\exists f_1,\dots,f_p\in \mathbf F\losp
\Bigl(\Bigl(\,\bigwedge\limits_{q\ne s} f_q\ne f_s\Bigr)
\logic\land f_1g=f_1 \logic\land \dots \logic\land
f_p=gf_p \logic\land hf_1=fhf_1\ne 0 \logic\land \dots \logic\land
hf_p=fhf_p\ne 0 \logic\land
\Bigl(\bigwedge\limits_{q\ne s} f_q\sim_h f_s\Bigr)\Bigr)\Bigr)$,
this means that there exist at most $p-1$ projections
from~$\mathbf F$ onto subgroups from~$gA$
such that their images are included in~$\mathbf F^D$ and their kernels are
isomorphic.
\end{enumerate}

Two elements $h_1$~and~$h_2$ from the set $\Endom_g$ are said to be
equivalent ($h_1\sim h_2$) if and only if the formula
$$
\exists f_1\losp \exists f_2\losp
( (gf_1g)\cdot (gf_2g)=(gf_2g)\cdot (gf_1g)=g
\logic\land \forall f\in \mathbf F\losp
(fg=f\Rightarrow h_1f=(gf_1g)h_2f \logic\land h_2f=(gf_2g)h_1f))
$$
is true, i.e., there exist mutually inverse automorphisms $gf_1g$ and
$gf_2g$ of~$gA$ which
map $h_2$~and~$h_1$ to automorphisms $(gf_1g) h_2$ and
$(gf_2g)h_1$ such that
we have $h_2=(gf_2g)h_1$ and $h_1=(gf_1g)h_2$ on~$gA$
for every projection from~$\mathbf F$
that project onto a~subgroup of~$gA$.

Again the obtained set $\Endom_g/{\sim}$ will be denoted by $\tEnd_g$.

Every class $h\in \tEnd_g$ can be mapped
to a~set consisting of the following $k+1$ components:
its $i$th component (where $i=1,\dots,k$) is a~set of pairs
$$
M_i=\{ \langle f,m(f)\rangle\mid f\in \mathbf F_i^D,\ m=1,\dots,p^i-1\},
$$
where $m$ is the dimension of the inverse image $f\in \mathbf F_i^D$, if it
is not equal to zero;
its $(k+1)$th component is a~set of sequences
$$
M=\{ \langle f,m(f),l_0(f),\dots,l_{m(f)}(f)\rangle\mid f\in \mathbf F^D,\
m\in \mathbb N,\ l_1,\dots,l_m=0,\dots,p-1\},
$$
where $p^m$ is the maximal order of the kernel of the inverse image
of~$fA$ which is included in~$gA$
and has the form $f'A$ for $f'\in \mathbf F'$,
$l_i$~is the number of those inverse images of~$fA$ that belong to $gA$,
have the form $f'A$ for $f'\in \mathbf F'$, and their kernels
have the order $p^i$ under~$h$.

{\sloppy
It is clear that endomorphisms from one equivalence class
are mapped to
the same sets $M_1,\dots, M_k, M$,
and endomorphisms from different
classes are mapped to different sets.
Further, every sequence $M_1,\dots,M_k,M$ of finite sets of the described
form are mapped to some class from $\tEnd_g$.
Such a~sequence of sets $M_1,\dots,M_k,M$ is mapped to the element
$$
\sum_{\langle f,m(f)\rangle\in M_1} m(f)a(f)+\dots +
\sum_{\langle f,m(f)\rangle\in M_k} m(f)a(f)+
\sum_{\langle f,m(f),l_0,\dots,l_{m(f)}\rangle \in M}
l_1 c_1(f)+\dots+l_{m(f)} c_{m(f)}(f),
$$
where $a(f)$ is a~fixed generator of the cyclic group~$fA$
for $f\in \mathbf F_1\cup\dots\cup \mathbf F_k$, and
$c_1(f),\dots,c_n(f),\dots$ are fixed generators
of the quasicyclic group $fA$ for $f\in \mathbf F$.

}

Therefore we obtain a~bijection between the set $\tEnd_g$
and the group~$A$.

The addition ($h_3=h_1\oplus h_2$) on the set $\tEnd_g$
is introduced by a~formula which is similar to the union
of the formulas from Secs.\ \ref{sec5_3}~and~\ref{sec6_2},
so we shall not write it here.

\begin{proposition}\label{p6.3}
Let $A_1=D_1\oplus G_1$, $A_2=D_2\oplus G_2$,
the group $D_1$~and~$D_2$ be divisible and infinitely generated,
the groups $G_1$~and~$G_2$ be bounded by the number~$p^k$,
$|D_1|\ge |G_1|$, and $|D_2|\ge |G_2$.
Then $\Endom (A_1)\equiv \Endom(A_2)\Rightarrow
A_1\equiv_{\mathcal L_2} A_2$.
\end{proposition}

\begin{proof}
The proof is completely similar to the proof
of Proposition~\ref{p6.2}, therefore we shall not write it here.
\end{proof}

\subsection[Direct Sums of Divisible $p$-Groups
and Bounded $p$-Groups of Greater Power]{Direct Sums
of Divisible $\boldsymbol{p}$-Groups
and Bounded $\boldsymbol{p}$-Groups of Greater Power}
This case differs from the two previous cases, it is closer
to the case of bounded $p$-groups, but it is more complicated.

We shall consider groups of the form $D\oplus G$, where $D$
is an infinitely generated divisible group, $G$~is a~group bounded
by the number~$p^k$, $|G|> |D|$,
$\smash[b]{G=\sum\limits_{i=1}^k G_i}$,
$G_i\cong \smash[b]{\bigoplus\limits_{\mu_i} \mathbb Z(p^i)}$,
$\mu_l= \smash[b]{\max\limits_{i=1,\dots,k} \mu_i}$,
and
$D\cong \bigoplus\limits_{\mu} \mathbb Z(p^\infty)$, $\mu< \mu_l$.

Assume that (as in the previous section)
we have idempotents $\rho_D$~and~$\rho_G$
which are projections on the divisible and bounded parts of
the group~$A$, respectively; and also projections
$\rho_1,\dots,\rho_k$, where $\rho_1+\dots+\rho_k=\rho_G$,
on the direct summands $G_1,\dots,G_k$.
Further, assume that
the summand~$G_l$ with the maximal power~$\mu_l$,
equal to the power of the whole group~$A$, is known.

As above, we introduce various definable sets:
\begin{enumerate}
\item
the set $\mathbf F=\mathbf F(\bar g)$ consists of $\mu$~indecomposable
projections on linearly independent direct summands of~$D$;
\item
for every $i=1,\dots,k$ the set $\mathbf F_i=\mathbf F_i(\bar g_i)$
consists of $\mu_i$ projections on independent
indecomposable direct summands of the group $G_i=\rho_i A$;
\item
the set $\mathbf F'=\mathbf F'(\bar g')$ consists of
$\mu_l$~projections on linearly independent countably generated
direct summands of~$G_l$;
\item
an idempotent~$\gamma$ satisfies the following condition:
$$
\Gamma(\gamma) \prisv
(\gamma \rho_l=\gamma \logic\land \gamma^2=\gamma \logic\land
\forall f\in \mathbf F'\losp \exists \beta\losp
(f \gamma=\gamma f=\beta \logic\land \Idem^\omega(\beta))).
$$
This condition means that $\gamma$ is a~projection on such
a~direct summand in~$G_l$ that its intersection with every
subgroup $fA$, where $f\in \mathbf F'$,
is a~countably generated summand of~$G_l$;
\item
for every idempotent~$\gamma$ satisfying the formula $\Gamma(\gamma)$,
by $\Gamma_\gamma$ we shall denote the set
$\{ f\in \mathbf F_l\mid f\gamma=f\}$, and by $\Gamma_\gamma(g)$
for $g\in \mathbf F'$ we shall denote the set
$\{ f\in \mathbf F_l \mid f\gamma=f \logic\land fg=f\}$.
Let us fix two of these idempotents $\gamma_0$~and~$\gamma_1$
with the conditions:
(1)~$\Gamma_{\gamma_0}\cap \Gamma_{\gamma_1}=\varnothing$;
(2)~for every $g\in \mathbf F'$ the set
$ \mathbf F_l\setminus
(G_{\gamma_0}\cup G_{\gamma_1})\cap \{ f\in \mathbf F_l\mid fg=f\}$
is countable.
\setcounter{slast}{\value{enumi}}
\end{enumerate}

Denote $\Gamma_{\gamma_0}$ by~$\Gamma_0$
and $\Gamma_{\gamma_1}$ by~$\Gamma_1$.

\begin{enumerate}
\setcounter{enumi}{\value{slast}}
\item
Fix
an endomorphism $\varphi\in \Endom(A)$ satisfying the following formula:
\begin{align*}
& \Phi(\varphi) \prisv
\forall h\losp
(\Idem(h) \logic\land \forall g\in \mathbf F'\losp (hg=gh=0)
\Rightarrow \varphi h=h\varphi=0)
\\
& \quad \logic\land
\forall g\in \mathbf F'\losp
(\forall f\in \Gamma_0(g)\losp
(\varphi f=f \logic\land \forall f\in \Gamma_1(g)\losp
\exists f' \in \mathbf F_l\losp
(f'\notin \Gamma_1(g) \logic\land f'\notin \Gamma_0(g)
\\
& \quad \logic\land f'g=f' \logic\land f'\varphi f=\varphi f\ne 0)
\logic\land \forall f\in \mathbf F_l\losp
(f\notin \Gamma_0(g) \logic\land fg=f
\\
& \quad {}\Rightarrow \exists f'\in \mathbf F_l\losp
(f'\ne f \logic\land f'\notin \Gamma_0(g) \logic\land f'\notin \Gamma_1(g)
\logic\land f'g=f' \logic\land f'\varphi f=\varphi f\ne 0)
\\
& \quad \logic\land \forall f_1,f_2\in \mathbf F_l\losp
(f_1\ne f_2 \logic\land f_1 g=f_1 \logic\land f_2 g=f_2
\\
& \quad {}\Rightarrow
\neg
(\exists f'\in \mathbf F_l\losp
(f'g=f' \logic\land f'\varphi f_1=\varphi f_1 \logic\land
f'\varphi f_2=\varphi f_2))
\\
& \quad {}\logic\land \forall f'\in \mathbf F_l\losp
(f'g=f' \logic\land f'\in \Gamma_1(g)\Rightarrow
\exists f\in \mathbf F_l\losp
(fg=f \logic\land f'\varphi=f'\varphi f=\varphi f))
\\
& \quad \logic\land
\forall h\losp (\Idem(h) \logic\land hg=h \logic\land
h\gamma_0=\gamma_0 h=0
\logic\land \forall f\losp
(\Idem^*(f) \logic\land fg=f \logic\land fh=f
\\*
& \quad {}\Rightarrow \exists f'\losp
(\Idem^*(f') \logic\land f' g=f' \logic\land f'h=f' \logic\land
f'\varphi f=\varphi f)\Rightarrow \Idem^\omega (f)))).
\end{align*}
\setcounter{slast}{\value{enumi}}
\end{enumerate}

This condition introduces an endomorphism~$\varphi$ that maps the
complementary direct summand for $\sum\limits_{g\in \mathbf F'} gA$
to~$0$, i.e.,
acts only on $\sum\limits_{g\in \mathbf F'} gA$
in the following way: for every $g\in \mathbf F'$ the elements
$\alpha A$, where $\alpha\in \Gamma_0(g)$,
are mapped into themselves, and the elements
$\alpha A$, where $\alpha\in \Gamma_1(g)$,
are mapped somewhere to the elements of
$\mathbf F_lA\setminus (\Gamma_0(g)A\cup\Gamma_1 (g)A)$
which are included in~$gA$.
Further, $\varphi$ is a~monomorphism on~$gA$, and its image is
$$
\langle \Gamma_0(g)A\rangle\oplus
\langle \{ fA\mid f\in \mathbf F_l \logic\land fg=f \logic\land
f\notin \Gamma_1(g)A\}\rangle.
$$
Outside $\Gamma_0(A)$ there are no finite-dimensional proper
subspaces of this endomorphism.
Therefore, we can numerate all elements from $\mathbf F_l$ that project
on subgroups from~$gA$ (we shall denote this set by $F_l(g)$)
by the following:
$f_i^j=f_i^j(g)$, $i=0,1,\dots$, $j=1,\dots$,
and
\begin{enumerate}
\renewcommand{\theenumi}{\alph{enumi}}
\item
$f_0^j\in \Gamma_0(g)$;
\item
$f_1^j\in \Gamma_1(g)$;
\item
$\varphi(f_i^jA)=f_{i+1}^jA$ if $i> 0$;
\item
$\varphi(f_0^jA)=f_0^jA$.
\end{enumerate}

We shall denote the set $\{ f_i^j\}_{j=1,\dots}$ by $\Gamma_i(g)$
(note that for an arbitrary~$i$ this set is not definable).

\begin{enumerate}
\setcounter{enumi}{\value{slast}}
\item
The union $\bigcup\limits_{g\in \mathbf F'} \mathbf F_l(g)$
will be denoted by~$\mathbf F_l'$. This set is definable;
\item
note that on the group $B=\langle \mathbf F_l'A\rangle$
the endomorphism~$\varphi$ has a~left inverse endomorphism~$\psi$ such that
$\psi\circ \varphi=1_B$.
For every~$g$ we shall introduce $gA$ as follows:
$$
\begin{cases}
\psi(f_0^jA)=f_0^jA,\\
\psi(f_i^jA)=f_{i-1}^jA\text{ if }i>1,\\
\psi(f_1^jA)\text{ can be arbitrary}.
\end{cases}
$$
We shall consider $\psi$ with the condition $\psi(f_i^jA)=0$.
Then two elements $f_1,f_2\in \mathbf F_l(g)\setminus \Gamma_0(g)$
(or, more generally, $\mathbf F_l'\setminus \Gamma_0$)
will be called \emph{$\varphi$-equivalent}
($f_1\sim_{\varphi}f_2$) if
\begin{align*}
& \exists h_1\losp \exists h_2\losp \exists \alpha\losp
\biggl(h_1g=h_1 \logic\land h_2g=h_2 \logic\land \Idem(h_1) \logic\land
\Idem(h_2) \logic\land \alpha^2=1
\\[1mm]
& \quad \logic\land
\smash[t]{\bigwedge_{i=1}^2 \forall f}\losp
(\Idem^*(f) \logic\land fg=f \logic\land fh_i=f
\Rightarrow
\exists f'\losp
(\Idem^*(f') \logic\land f'g=f' \logic\land f'h_i=f' \logic\land
f'\psi f=\psi f))
\\
& \quad \logic\land f_1 h_1=f_1 \logic\land f_2 h_2=f_2 \logic\land
\smash[t]{\bigwedge_{i=1}^2 \forall h}\losp
(\Idem(h) \logic\land hg=h
\logic\land \forall f\losp
(\Idem^*(f) \logic\land fg=f \logic\land fh=f
\\
& \quad \Rightarrow \exists f'\losp
(\Idem^*(f') \logic\land f'g=f' \logic\land f'h=f' \logic\land
f'\psi f=\psi f) \logic\land f_i hf_i \Rightarrow h_ih=h_i))
\\
& \quad \logic\land h_1\alpha h_2=\alpha h_2=h_1\alpha \logic\land
h_2 \alpha h_1=h_2\alpha\biggr).
\end{align*}
\end{enumerate}

This means that minimal proper subspaces of the endomorphism~$\psi$
($h_1A$ and $h_2A$), containing the groups $f_1A$ and $f_2A$, respectively,
have the same power, i.e.,
$f_1A,f_2A\in \Ker \psi^m\setminus \Ker \psi^{m+1}$
for some natural~$m$, or $f_1A,f_2A\in \varphi^m(\Gamma_1)$.
It is clear that in this case $f_1=f_m^j(g)$ and $f_2=f_m^i(g)$
(or, in more general case,
$f_1=f_m^j(g_1)$ and $f_2=f_m^j(g_2)$).

We shall call an element $f_1\in \mathbf F_l(g)$
a~\emph{$\varphi$-successor} of the element
$f_2\in \mathbf F_l(g)$ ($f_1\sim_\varphi f_2+1$) if
$$
\exists f\losp (f\sim_\varphi f_1 \logic\land
f\varphi f_2=\varphi f_2=f\varphi).
$$

A~similar formula will introduce a~notion of
a~\emph{$\varphi$-greater element}
($f_1>_\varphi f_2$) as an element for which
a~minimal proper subspace
of an endomorphism~$\psi$ containing $f_1A$ has greater power
than the corresponding subspace for~$f_2A$.

Let us fix $g\in \mathbf F_l'$ and
construct an interpretation of the group $A=D\oplus G$ for this~$g$.

Consider the set $\Endom_g$ of all homomorphisms $h\colon gA\to A$
satisfying the following conditions:
\begin{enumerate}
\item
$\forall f\in \mathbf F_l\losp
\Bigl(fg=f\Rightarrow \Bigl(hf=0 \logic\lor
\exists \tilde f\in \mathbf F^D\losp
(hf=\tilde f hf\ne 0) \logic\lor
\Bigl(\,\bigvee\limits_{i=1}^k \exists \tilde f\in \mathbf F_i\losp
(hf=\tilde fhf\ne 0)\Bigr)\Bigr)\Bigr)$;
\item
$\exists f\losp (\Fin(f) \logic\land \Idem(f) \logic\land fh=h)$;
\item
$\bigwedge\limits_{i=1}^k \forall f\in \mathbf F_i\losp
\neg \Bigl(\exists f_1,\dots,\exists f_{p^i}\in \mathbf F\losp
\Bigl(\,\bigwedge\limits_{q\ne s}
f_q\ne f_s \logic\land f_1,\dots,f_{p^i}\in \Gamma_0(g) \logic\land
hf_1=fhf_1\ne 0 \logic\land \dots \logic\land
hf_{p^i}=fhf_{p^i}\ne 0\Bigr) \logic\land
\forall f'\in \mathbf F_l\losp
(hf'=fhf'\ne 0\Rightarrow f'\in \Gamma_0(g))\Bigr)$;
this means that the inverse images of the elements
of the bounded summand can be contained only in the set $\Gamma_0(g)$;
\item
$ \forall f\in \mathbf F^D\losp \exists f'\in \mathbf F_l\losp
(hf'=fhf'\ne 0\Rightarrow f'\notin \Gamma_0)$;
contrary to the previous assertion, this means that the inverse images
of the elements of the divisible summand can not be contained
in~$\Gamma_0(g)$;
\item
$\forall f\in \mathbf F^D\losp
\neg \Bigl(\exists f_1,\dots,f_p\in \mathbf F_l\losp
\Bigl(\Bigl(\,
\bigwedge\limits_{q\ne s} f_q\ne f_s \logic\land f_q\sim_\varphi f_s\Bigr)
\logic\land f_1g=f_1 \logic\land\dots \logic\land
f_pg=f_p \logic\land hf_1=fhf_1\ne 0 \logic\land \dots \logic\land
hf_p=fhf_p\ne 0\Bigr)\Bigr)$;
this means that no element from~$\mathbf F_D$ can have more than
$p-1$ $\varphi$-equivalent inverse images;
\item
$\forall f\in \mathbf F^D\losp
\exists f'\in \mathbf F_l\losp
\neg (\exists f''\in \mathbf F_l\losp
(f''>_\varphi f' \logic\land hf''=fhf''\ne 0))$,
i.e., every element from~$\mathbf F_D$
contains only a~finite number of inverse images in~$\mathbf F_l$.
\end{enumerate}

Two elements $h_1$~and~$h_2$ of the set $\Endom_g$
are said to be equivalent
(notation: $h_1\sim h_2$) if we have the following formula:
\begin{align*}
& \exists f_1\losp \exists f_2\losp
((gf_1g)\cdot (gf_2g)=(gf_2g)\cdot (gf_1g)=g
\\
& \quad \logic\land
\forall f\in \mathbf F_l\losp
(fg=f
\Rightarrow
\forall f'\in \mathbf F^D\cup \mathbf F_1\cup \dots \cup \mathbf F_k\losp
(h_1f=f'h_1f\ne 0\Leftrightarrow (gf_1gh_2)f=f'(gf_1gh_2)f\ne 0))
\\
& \quad \logic\land \forall f\in \mathbf F_l\losp
(gf_1g\cdot f\sim_\varphi f)).
\end{align*}

The obtained set $\Endom_g/{\sim}$ will be denoted by $\tEnd_g$.

Now we shall show how to find a~bijection between the set
$\tEnd_g$ and the group~$A$.
Let us consider some $h\in \tEnd_g$. For every $f\in \mathbf F_i$
let us consider
the intersection of the inverse image of $fA$ with the set $\Gamma_0 A$.
Suppose that this inverse image contains $m_f$~elements.
Thus, we get the set
$$
M_G =\bigcup\limits_{i=1}^k \{ \langle f,m_f\rangle\mid
f\in \mathbf F_i,\ m_f=1,\dots ,p^i-1\}.
$$
For every $f\in \mathbf F^D$ and every natural~$j$
let us consider the intersection
of the inverse image of~$fA$ with the set $\Gamma_j(g)$. Suppose that
this inverse image
contains $l_f^j$~elements, and the maximal nonzero~$j$
is equal to $\gamma_f$. Then we get the set
$$
M_D=\{ \langle f,\gamma_f,l_f^1,\dots,l_f^{\gamma_f}\rangle\mid
f\in \mathbf F^D,\ \gamma_f\in \mathbb N,\
l_f^1,\dots,l_f^{\gamma_f}\in \{ 0,\dots,p-1\}\}.
$$
Now an element~$h$ will be mapped to the following sum:
$$
\sum_{f\in M_G} m_f a(f)+\sum_{f\in M_D} l_f^1c_1(f)+\dots+
l_f^{\gamma_f}c_{\gamma_f}(f)\in A.
$$
It is clear that such a~mapping is a~bijection
between the sets $\tEnd_g$ and~$A$, and this bijection becomes
an isomorphism if and only if we introduce addition on the set $\tEnd_g$
with the help of a~formula similar to the formulas
from Secs.\ \ref{sec5_3}~and~\ref{sec6_2}.
Therefore we have the following proposition.

\begin{proposition}\label{p6.4}
Let $A_1=D_1\oplus G_1$, $A_2=D_2\oplus G_2$,
the groups $D_1$~and~$D_2$ be divisible, the groups $G_1$~and~$G_2$
be infinite and bounded by the number~$p^k$, $|D_1|< |G_1|$, and
$|D_2|< |G_2|$. Then
$\Endom (A_1)\equiv \Endom(A_2)\Rightarrow A_1\equiv_{\mathcal L_2} A_2$.
\end{proposition}

\section{Groups with Unbounded Basic Subgroups}\label{sec7}
\subsection[The Case $A=D\oplus G$, Where
$|D|\ge |G|$, and Other Cases]{The Case $\boldsymbol{A=D\oplus G}$,
Where $\boldsymbol{|D|\ge |G|}$, and Other Cases}\label{sec7_1}
Let us separate our problem into three cases.

1. $A=D\oplus G$, where $|D|\ge |G|$ and $G$~is any unbounded group.
We shall not consider this case in details, because its proof
is similar to the proof from Sec.~\ref{sec6_2}.

This case resembles the case $A=D\oplus G$, where $|D|\ge |G|$,
$D$~is a~divisible group,
and $G$~is a~bounded group (see Sec.~\ref{sec6_3}).
Here we give only a~sketch of the proof.

Since $|G|\le |D|$, we have that a~basic subgroup of the group~$G$
(or the group~$A$) is of power not greater than
the power of~$D$. Hence there exists
an embedding $\varphi_1\colon B\to D_1$, where $D=D_1\oplus D_2\oplus D_3$
and
$|D|=|D_1|=|D_2|=|D_3|$.
This embedding will be fixed,
and after that we can assume that the group~$B$
is a~subgroup of~$D_1$. Further, $|G|\le |D|$ implies $|G/B|\le |D|$,
whence there exists an embedding $\varphi_2\colon G/B\to D_2$
(i.e., a~mapping from~$G$ into~$D_2$ which is equal to zero on~$B$).
Thus the group $G/B$ can also be considered as a~subgroup of~$D_2$.

Now we shall find some definable sets:
\begin{enumerate}
\item
the set~$\mathbf F_1$ consists of $|B|$~independent indecomposable
projections
on quasicyclic direct summands in a~minimal direct summand of~$D_1$,
containing $\varphi_1(B)$ as a~subgroup;
\item
the set~$\mathbf F_2$ consists of $|G/B|$ independent indecomposable
projections on quasicyclic direct summands of $\varphi_2(G/B)$;
\item
the sets $\mathbf F$~and~$\mathbf F_3$ consists of $\mu=|D|$ independent
projections on quasicyclic direct summands of the groups $D$~and~$D_3$,
respectively;
\item
the set $\mathbf F'$ consists of
$\mu$~independent projections on countably generated
direct summands of the group~$D$.
\end{enumerate}

For every $g\in \mathbf F'$ an interpretation of the group~$A$
will be constructed in the following way:
we shall consider homomorphisms $h\colon gA\to A$ such that the images
of the subgroups $fA$ ($f\in \mathbf F$, $fA\subset gA$)
are either zero or are contained in~$f'A$
($f'\in \mathbf F_1\cup \mathbf F_2\cup \mathbf F_3$), and
$h(gA)$ is finite-dimensional.

The inverse images of~$fA$, where $f\in \mathbf F_1$, will interpret
the summands from~$B$ in the decomposition of the element $a\in A$
in a~quasibasis; the inverse images of~$fA$,
where $f\in \mathbf F_2$, are the summands from $G/B$, i.e.,
$c_{i,j}$ for $i\in \omega$, $j\in |G/B|$;
the inverse images of~$fA$, where $f\in \mathbf F_3$,
are the summands from~$D$.
The rest of the procedure is similar to that
from the previous sections.

So we have given a~sketch of the proof of the following proposition.

\begin{proposition}\label{p7.1}
Let $A_1=D_1\oplus G_1$, $A_2=D_2\oplus G_2$,
where the groups $D_1$,~$D_2$ are divisible,
the groups $G_1$,~$G_2$ are reduced, $|D_1|\ge |G_1|$, and $|D_2|\ge |G_2|$.
Then $\Endom(A_1)\equiv \Endom(A_2)\Rightarrow A_1\equiv_{\mathcal L_2} A_2$.
\end{proposition}

In this section, we shall assume that $A=D\oplus G$, where $|D|< |G|$.

\smallskip

2. $A=D\oplus G$, where $|D|< |G|$, $B$ is a~basic subgroup in~$G$,
and $r(B)=r_{\mathrm{fin}}(B)$.

The case $r_{\mathrm{fin}}(B)> \omega$ will be considered
in Sec.~\ref{sec7_4},
and the case $r_{\mathrm{fin}}=\omega$ will be considered
in Sec.~\ref{sec7_5}.

If $r(B)> \omega$, then $|A|=r(B)$ and $|D|< r(B)$.
If $r(B)=\omega$, then $|A|\le |\mathcal P(\omega)|$, therefore
if we do not assume the continuum-hypothesis, then we can meet
the situation, where $\omega< |D| < |A|\le 2^\omega$,
which is not good for us. Thus for the simplicity of arguments,
we shall assume the continuum-hypothesis.

Hence, if $A=D\oplus G$, where $|D|< |G|$, $r(B)=r_{\mathrm{fin}}(B)$, then
we shall interpret the theory $\Th_2^{r(B)}(A)$ in the ring $\Endom(A)$.

\smallskip

3. $A=D\oplus G$, where $|D|< |G|$ and $r(B)\ne r_{\mathrm{fin}}(B)$.
If in this case $r_{\mathrm{fin}}(B)> \omega$, then we can obtain the full
second order theory of the group~$A$. This case is considered in
Sec.~\ref{sec7_4}.

If $r_{\mathrm{fin}}(B)=\omega$, then, assuming the continuum-hypothesis,
we can find in the group~$A$ a~bounded direct summand
of power~$|A|$, and in this case we can define the complete
second order theory of the group~$A$.
This case is considered in Sec.~\ref{sec7_6}.

In Secs.\ \ref{sec7_2}~and~\ref{sec7_3},
we shall find some definable objects, which are important for all cases.

\subsection{Definable Objects}\label{sec7_2}
In this section, we assume that $A=D\oplus G$,
the group~$D$ is divisible (it can be zero),
the group~$G$ is reduced and has an unbounded basic subgroup~$B$,
$$
B=B_1\oplus \dots\oplus B_n\oplus \dots,
$$
where
$$
B_n\cong \bigoplus_{\mu_n} \mathbb Z(p^n),
$$
$r(D)=\mu_D$, $|B|=\bigcup\limits_{n\in \mathbb N} \mu_n=\mu_B$,
$|G|=\mu_G$ (if $\mu_B> \omega$, then $\mu_G=\mu_B$), and
$\mu=|A|=\max(\mu_D,\mu_G)$.

We suppose that projections $\rho_D$~and~$\rho_G$
on the summands $D$~and~$G$ of the group~$A$, respectively, are fixed.

By~$Z$ we shall denote the center of the ring $\Endom(A)$.
As we remember (see Theorem~\ref{t2.9}),
each of its elements multiplies all elements of~$A$
by some fixed $p$-adic number.

For any indecomposable projections $\rho_1$~and~$\rho_2$ from~$\Endom(G)$
we shall write $o(\rho_1)\le o(\rho_2)$ if
$$
\forall c\in Z\losp (c\rho_2=0\Rightarrow c\rho_1=0).
$$
It is clear that this formula holds if and only if the order
of the finite cyclic direct summand~$\rho_1A$ is not greater than the order
of the summand~$\rho_2A$.

Similarly,
\begin{align*}
(o(\rho_1)< o(\rho_2))&\prisv
(o(\rho_1)\le o(\rho_2)) \logic\land \neg (o(\rho_2)\le o(\rho_1)),\\
(o(\rho_1)= o(\rho_2))&\prisv
(o(\rho_1)\le o(\rho_2)) \logic\land (o(\rho_2)\le o(\rho_1)).
\end{align*}

For every indecomposable projection~$\rho$ we shall consider the following
formula sets.

1. The formula
$$
\mathrm{Ord}_{\rho}(f)\prisv \Idem(f) \logic\land
\forall f'\losp (\Idem^*(f') \logic\land f'f=f'\Rightarrow o(f')=o(\rho))
$$
defines the projections~$f$ on direct summands $fA$ in~$A$
which are direct sums of cyclic groups of order~$o(\rho A)$.

2. The formula
$$
\mathrm{MaxOrd}_{\rho}(f)\prisv \Idem(f) \logic\land
\mathrm{Ord}_\rho(f) \logic\land
\forall f'\losp (\mathrm{Ord}_\rho(f') \Rightarrow \neg (ff'=f))
$$
defines the projections~$f$ on maximal direct summands~$fA$ in~$A$
which are direct sums of cyclic groups of order~$o(\rho A)$.

3. The formula
$$
\mathrm{Rest}_{\rho}(f) \prisv \Idem(f) \logic\land
\forall f'\losp (\Idem^*(f') \logic\land f'f=f'\Rightarrow o(f')\le o(\rho))
$$
defines the projections~$f$ on direct summands $fA$ in~$A$
which are direct sums of cyclic groups of order at most $o(\rho A)$.

4. The formula
$$
\mathrm{MaxRest}_{\rho}(f) \prisv
\Idem(f) \logic\land \mathrm{Ord}_\rho(f) \logic\land
\forall f'\losp (\mathrm{Ord}_\rho(f') \Rightarrow \neg (ff'=f))
$$
defines the projections~$f$ on maximal direct summands $fA$ in~$A$
which are direct sums of cyclic groups of order at most $o(\rho A)$.

5. The formula
$$
\overline{\Base}(\varphi) \prisv
\forall \rho\losp \exists f\losp
(\mathrm{MaxRest}_\rho(f) \logic\land
\forall f'\losp (\Idem^*(f')\Rightarrow
(f'f=f'\Leftrightarrow
\forall c\in Z\losp (cf'\ne 0\Rightarrow c(f'\varphi)\ne 0)))
$$
postulates that for every natural~$n$ there exists a~maximal
$p^n$-bounded direct summand of the group~$A$
which is included in~$\varphi A$. Therefore,
the group~$\varphi A$ necessarily
contains some basic subgroup of the group~$A$.

6. The formula
\begin{multline*}
\Base(\varphi) \prisv
\overline{\Base}(\varphi) \logic\land
\forall f^*\losp (\Idem^*(f^*) \logic\land f^*\varphi \ne 0\Rightarrow
\exists \rho\losp \exists f\losp
(\mathrm{MaxRest}_\rho(f)
\\
\logic\land \forall f'\losp (\Idem^*(f')\Rightarrow
(f'f=f'\Leftrightarrow \forall c\in Z\losp
(cf'\ne 0\Rightarrow c(f\varphi)\ne 0))) \logic\land f^*f=f^*))
\end{multline*}
is true for every endomorphism $\varphi\in \Endom(A)$ whose image
is a~basic subgroup in~$A$.

Let us suppose that we have a~fixed endomorphism~$\varphi_B$
such that $\Base(\varphi_B)$.

\subsection{Definable Special Sets}\label{sec7_3}
We shall consider three different cases:
\begin{enumerate}
\item
$\mu_B=\omega$;
\item
$\mu_B> \omega$ and
$\forall k\in \mathbb N\losp \exists n\in \mathbb N\losp
(n> k \logic\land \mu_n=\mu_B)$.
This is always true if $\cf \mu_B >\omega$;
\item
$\mu_B> \omega$, $\cf \mu_B=\omega$,
$\forall n\in \mathbb N\losp (\mu_n<\mu_B)$.
\end{enumerate}

\smallskip

Case 1. $\mu_B=\omega$.

Let us consider the formula
\begin{align*}
& \mathrm{Intr}(f) \prisv
[\forall f' \losp
(\Idem(f') \logic\land f'\varphi_B\ne 0\Rightarrow f'f\ne 0)]\logic\land
[\forall f_1\losp \forall f_2\losp
(\Idem^*(f_1) \logic\land \Idem^*(f_2) \logic\land o(f_1)=o(f_2)
\\
& \quad {}\logic\land \forall c\in Z\losp
(cf_1\ne 0\Rightarrow cf_1f\ne 0) \logic\land \forall c\in Z\losp
(cf_2\ne 0\Rightarrow cf_2 f\ne 0)\Rightarrow
f_1f_2\ne 0 \logic\land f_2f_1\ne 0)]
\\
& \quad {}\logic\land
[\forall \rho'\losp
(\Idem^*(\rho')\Rightarrow \exists f'\losp
(\Idem^*(f') \logic\land o(f') > o(\rho') \logic\land
\forall c\in Z\losp
(cf'\ne 0\Rightarrow cf'f\ne 0)))].
\end{align*}
The first part of this formula, enclosed in square brackets,
postulates $fA\subset \varphi_B A=B$. The second part states that
there is at most one
cyclic direct summand of the same order in the image of~$fA$.
The third part states that
the orders of direct summands in $fA$ are unbounded.

Therefore this formula gives us an endomorphism~$f$ with image~$B'$
being a~cyclic direct summand in~$B$,
$$
B'\cong \bigoplus_{i\in \mathbb N} \mathbb Z(p^{n_i}),
$$
where $(n_i)$ is an increasing sequence.

This endomorphism is supposed to be fixed and is denoted by~$f_B$.

Now we shall consider endomorphisms from~$B'$ into~$A$.
Namely we shall consider only functions~$f$ with the condition
$$
\forall \rho\losp (\Idem(\rho) \logic\land \rho f_B=0\Rightarrow\rho f=0).
$$
Two functions $f_1$~and~$f_2$ satisfying this condition
are said to be equivalent if
$$
\forall \rho\losp (\Idem^*(\rho) \logic\land
\forall c\in Z\losp (c\rho\ne 0\Rightarrow c\rho f\ne 0)
\Rightarrow f_1 \rho=f_2\rho),
$$
i.e., if they coincide on the group~$B'$.

Consequently, if we factorize the set of all described functions
by this equivalence, then we obtain the group $\Hom(B',A)$.

Introduce now the formula
$$
o(\rho_1)\ge o(\rho_2)^2
$$
for two indecomposable idempotents $\rho_1$~and~$\rho_2$ as follows:
$$
\forall c\in Z\losp (c\rho_2\ne 0\Rightarrow c^2 \rho_1\ne0).
$$
This formula means that $|\rho_1A|\ge |\rho_2A|^2$.

Similarly we can introduce the formulas
$$
o(\rho_1)> o(\rho_2)^2\ \ \text{and}\ \ o(\rho_1)=o(\rho_2)^2.
$$
Suppose now that our function~$f_B$ satisfies the additional condition
\begin{align*}
& \forall f'\losp
(\Idem^*(f') \logic\land \forall c\in Z\losp
(cf'\ne 0\Rightarrow cf'f_B\ne 0)\Rightarrow pf'\ne 0)
\\*
& \quad \logic\land
\forall \rho'\losp \forall f'\losp
(\Idem^*(f') \logic\land \forall c\in Z\losp
(cf'\ne 0\Rightarrow cf'f_B\ne 0) \logic\land o(f')=o(\rho')
\\*
& \quad
\Rightarrow \forall f\losp
(\Idem^*(f) \logic\land \forall c\in Z\losp
(cf\ne 0\Rightarrow cff_B\ne 0) \logic\land
o(f)> o(\rho')\Rightarrow o(f)> o(\rho')^2)).
\end{align*}
This condition means that
\begin{enumerate}
\item
every cyclic direct summand in~$B'$ of the smallest order has order
greater than~$p$ (i.e., at least not smaller than~$p^2$);
\item
for every direct cyclic summand in~$B'$ of order~$p^k$
the next cyclic summand of greater order has order
greater than~$p^{2k}$.
\end{enumerate}

Therefore,
$$
B'\cong \bigoplus_{i\in \omega} \mathbf Z(p^{n_i}),
$$
where $n_1\ge 2$, $n_{i+1}> 2n_i$.

Now consider the formula
\begin{align*}
& \mathrm{Ins}(\psi) \prisv
[\exists f\losp (\Idem^*(f) \logic\land \forall c\in Z\losp
(cf\ne 0 \Rightarrow cff_B\ne 0)
\\
& \quad \logic\land \forall f'\losp
(\Idem^*(f') \logic\land \forall c\in Z\losp
(cf'\ne 0\Rightarrow cf'f_B\ne 0) \Rightarrow
o(f)\le o(f')) \logic\land \psi f=pf)]
\\
& \quad \logic\land [\forall f_1\losp \forall c_1\in Z\losp
(\Idem^*(f_1) \logic\land \forall c\in Z\losp
(cf_1\ne 0\Rightarrow cf_1f_B\ne 0)
\\
& \quad {} \logic\land (\psi f_1=c_1f_1\Rightarrow \exists f_2\losp
(\Idem^*(f_2) \logic\land \forall c\in Z\losp
(cf_2\ne 0\Rightarrow cf_2f_B\ne 0) \logic\land o(f_2)> o(f_1)
\\
& \quad \logic\land \forall f'\losp
(\Idem^*(f') \logic\land \forall c\in Z\losp
(cf'\ne 0\Rightarrow cf'f_B\ne 0)
\Rightarrow
o(f')\le o(f_1) \logic\lor o(f')\ge o(f_2)) \logic\land \psi f_2=pc_1 f_2)))].
\end{align*}
The condition from the first square brackets states that there exists
a~cyclic summand of the smallest order such that the action
of the endomorphism~$\psi$ on it is multiplication by~$p$.
The second condition states that
for every natural~$i$ there exists a~direct cyclic summand
$\langle a_i \rangle$ of order~$p^{n_i}$ such that
the action of~$\psi$ on it is multiplication by~$p^i$.

Suppose that for some different cyclic direct summands
$\langle a_i\rangle$ and $\langle b_i\rangle$ from~$B'$ the action of~$\psi$
on them is multiplication by~$p^i$. Let
$b_i=\sum \alpha_k a_k +\sum \beta_l a_l + a_i$,
where $o(a_k)< o(b_i)$, $o(a_l)> o(b_i)$, $k< i$, $l> i$,
$$
\psi(b_i)=p^i b_i=\sum p^k \alpha_k a_k +\sum p^l \beta_l a_l+ p^i a_i=
\sum p^i \alpha_k a_k+\sum p^i \beta_l a_l +p^i a_i.
$$
Let $k< i$. Then we have $p^i\alpha_ka_k=p_k\alpha_ka_k=0$, i.e.,
$\alpha_k$~is divisible by~$p^{n_k-k}$
and $p^{n_k-k}$ is divisible by~$p^k$.
Therefore, we can write
$$
b_i=\sum \alpha_k p^ka_k+\sum \beta_k p^{n_l-n_i}a_l + a_i.
$$
We have $p^i\beta_l p^{n_l-n_i}a_l=p^l \beta_l p^{n_l-n_i}a_l=0$.
Note that every cyclic direct summand $\langle a\rangle$
either uniquely corresponds to an element of the center $c\in Z$
such that $\psi a=ca$, or there is no such element. We can consider
only those summands that correspond to elements of the center.

Let us suppose that we have some homomorphism $f\colon B'\to A$ such that
$o(f(a_i))\le p^i$. Let $\psi(a_i)=p^ia_i$ and $\psi(b_i)=p^i b_i$.
Let us find $f(b_i)$:
$$
f(b_i)=\sum \alpha_k p^k f(a_k)+\sum \beta_l p^{n_l-n_i} f(a_l)+f(a_i).
$$
Since $o(f(a_k))\le p^k$, we have $\sum \alpha_k p^k f(a_k)=0$.
Since $o(f(a_l))\le p^l$, we have $\sum \beta_l p^{n_l-n_i} f(a_l)=0$.
Therefore $f(b_i)= f(a_i)$.
Thus every element of the center~$Z$ of the form $p^n\cdot E$
is mapped (under this homomorphism $f\colon B'\to A$)
to some uniquely defined element $a\in A$ with the condition $o(a)\le p^n$.

\smallskip

Case 2. $\forall k\in \omega\losp \exists n\in\omega\losp
(n > k \logic\land \mu_n =\mu_B)$.

Consider the formula
\begin{multline*}
\mathrm{ECard}(\rho) \prisv \Idem^*(\rho) \logic\land
\exists \psi\losp \forall f\losp
(\Idem^*(f)\logic\land \forall c\in Z\losp
(cf\ne 0\Rightarrow cf\varphi_B\ne 0)
\\
\Rightarrow \exists f'\losp (\Idem^*(f') \logic\land
\forall c\in Z\losp (cf'\ne 0\Rightarrow cf'\varphi_B\ne 0)
\logic\land o(f')=o(\rho) \logic\land f\psi f'\ne 0)).
\end{multline*}
This formula states that the set of independent cyclic summands of order
$o(\rho)$ has the same power as the whole group~$B$, because there exists
a~homomorphism~$\psi$ from a~direct summand of the group~$B$
(which is isomorphic to the sum
of cyclic groups of order~$o(\rho)$)
such that its image intersects with every cyclic summand of~$B$.
Therefore, $\mu_{o(\rho)}=\mu_B$.

Now let us consider the formula
\begin{align*}
& \mathrm{Fine}(f) \prisv
[\forall f'\losp
(\Idem(f') \logic\land f'\varphi_B\ne 0 \Rightarrow f'f\ne 0)]
\\
& \quad {}\logic\land
[\forall f_1\losp \forall f_2\losp
(\Idem^*(f_1) \logic\land \Idem^*(f_2) \logic\land o(f_1)=o(f_2)
\\
& \quad {}\logic\land \forall c\in Z\losp
(cf_1\ne 0\Rightarrow cf_1 f\ne 0) \logic\land
\forall c\in Z\losp (cf_2\ne 0\Rightarrow cf_2 f\ne 0)
\Rightarrow f_1f_2\ne 0 \logic\land f_2f_1\ne 0)]
\\
& \quad \logic\land
[\forall \rho'\losp (\mathrm{ECard}(\rho')\Leftrightarrow
\exists f'\losp
(\Idem^*(f') \logic\land o(f')= o(\rho') \logic\land
\forall c\in Z\losp (cf'\ne 0\Rightarrow cf'f\ne 0)))].
\end{align*}
The first part of the formula, enclosed in square brackets, postulates
$fA\subset \varphi_B A=B$.
The second part, enclosed in the second square brackets,
states that the image of $fA$ does not contain more than one
cyclic direct summand of any order. The third part states that
all direct cyclic summands have order~$p^n$, where $\mu_n=\mu_B$.

To make the further constructions, we need
to recall Sec.~\ref{sec3}.

Formulation of Theorem~\ref{shel_t4_1} need not be changed,
but Lemma~\ref{shel_l4_2} and the proof of the theorem
with the help of this lemma must be changed a~little for our case.

There is a~new formulation of Lemma~\ref{shel_l4_2}:
\textit{there exists a~formula $\varphi(f)$ with one free variable~$f$
such that $\varphi(f)$ holds in $\Endom(B')$ if and only if
there exists an ordinal number $\alpha\in \mu$ such that for all
$\beta \in \mu$ and all $m\in \omega$}
$$
f(a_{\langle 0,\beta\rangle}^m)=a_{\langle \alpha,\beta\rangle}^m.
$$

Now we shall write the proof of the theorem with the help of the lemma.

Let a~function $f_0^*$ map (for every $m\in \omega$)
the set $\{ a_{\langle 0,\alpha\rangle}^m\mid \alpha\in \mu\}$
onto the set $\{ a_t^m\mid t\in I^*\}$,
and $f_0^*(a_{\langle \alpha,\beta\rangle}^m)=
f_0^*(a_{\langle 0,\beta\rangle}^m)$.
Suppose that we have a~set $\{ f_i\}_{i\in \mu}$ and let the function~$f^*$
be such that
$$
f^*(a_{\langle \alpha,\beta\rangle}^m)=
f_\alpha\circ f_0^*(a_{\langle \alpha,\beta\rangle}^m).
$$
Let $f_1^*$ map (for every~$m\in \omega$)
the set $\{ a_t^m\mid t\in I^*\}$ onto the set
$\{ a_{\langle 0,\beta\rangle}^m\mid \beta\in \mu\}$.
Let the formula $\tilde \varphi(f',\dots)$ say that there exists~$f$
such that
\begin{enumerate}
\item
$\varphi(f)$;
\item
$f'\circ f_0^*\circ f_1^*=f^*\circ f\circ f_1^*$.
\end{enumerate}

Then $\Endom(B')\vDash \varphi(f)$ if and only if
there exists $\alpha\in \mu$ such that
$$
\forall \beta\in \mu\losp \forall m\in \omega\losp
f(a_{\langle 0,\beta\rangle}^m)=a_{\langle \alpha,\beta\rangle}^m.
$$
Therefore
\begin{multline*}
f'\circ f_0^*\circ f_1^*(a_t^m)=f^*\circ f\circ f_1^*(a_t^m)
\Leftrightarrow
f'\circ f_0^*(a_{\langle 0,\beta\rangle}^m)=
f^*\circ f(a_{\langle 0,\beta\rangle}^m)
\\
{}\Leftrightarrow
f'\circ f_0^*(a_{\langle 0,\beta\rangle}^m)=
f^*(a_{\langle \alpha,\beta\rangle}^m)
\Leftrightarrow
f'\circ f_0^*(a_{\langle 0,\beta\rangle}^m)=
f_{\alpha}\circ f_0^*(a_{\langle 0,\beta\rangle}^m).
\end{multline*}

Let $f_0^*(a_{\langle 0,\beta\rangle}^m)=a_{t_\beta}^m$. Then
$$
f'(a_{t_\beta}^m)= f_{\alpha}(a_{t_\beta}^m),
$$
what we needed.

Now we need to change the proof of the lemma. The case $\mu=\omega$
is not interesting for us, therefore we shall begin from the second case.

Suppose that the cardinal number~$\mu_B$ is regular.
Represent it in the form of the union of the sets $I_0$, $I$, and~$J$,
where $|I_0|=|J|=|I|=\mu_B$,
$J=\{ \langle \alpha,\beta\rangle\mid \alpha,\beta\in \mu\}\cup \{ 0\}$,
$I=\{ \langle \alpha,\delta,n\rangle\mid \alpha\in \mu,\allowbreak\
{\delta\in \mu},\allowbreak \ {\cf \delta=\omega},\allowbreak\
n\in \omega\}$, and
$a_\alpha^{\beta,n}=a_{\langle \alpha,\beta,n\rangle}$.
As in Sec.~\ref{sec3},
for every limit ordinal $\delta\in \mu$ such that $\cf \delta=\omega$,
we choose an increasing sequence $(\delta_n)_{n\in \omega}$
of ordinal numbers less than~$\delta$ such that their limit is~$\delta$ and
for each $\beta\in \mu$ and $n\in \omega$ the set
$\{ \delta\in \mu\mid \beta=\delta_n\}$ is a~stationary subset in~$\mu$.

Consider an independent set of generators
of cyclic direct summands from~$B$ such that
\begin{enumerate}
\item
the order of every generator from this set is equal to~$p^n$,
where $\mu_n=\mu_B$;
\item
for every $n\in \omega$ such that $\mu_n=\mu_B$,
the set of all elements of order~$p^n$ from this set
has the power~$\mu_B$.
\end{enumerate}
Let us denote this set by
$$
\{ a_t^n\mid t\in J\cup I_0\cup I=I^*,\ o(a_t^n)=p^n,\ \mu_n=\mu_B\}.
$$
Let us define the functions $f_0^*,\dots,f_{14}^*$,
similar to the functions from the case~II from Sec.~\ref{sec3},
but with some addition:
$f_0^*(a_t^m)=a_0^m$,
$f_1^*(a_t^m)=a_t^{m-1}$ for $m> 0$, $f_1^*(a_t^0)=0$,
$f_2^*(a_{\langle \alpha,\beta\rangle}^m)=a_{\langle 0,0\rangle}^m$,
$f_3^*(a_{\langle \alpha,\beta\rangle}^m)=a_{\langle \alpha,0\rangle}^m$,
$f_4^*(a_{\langle \alpha,\beta\rangle}^m)=a_{\langle 0,\beta\rangle}^m$,
$f_5^*(a_{\langle \alpha,\beta\rangle}^m)=a_{\langle \beta,\alpha\rangle}^m$,
$f_6^*(a_{\langle \alpha,\beta\rangle}^m)=
a_{\langle \alpha,\alpha\rangle}^m$,
for $\delta\in \mu_B$, $\cf \delta=\omega$,
$f_7^*(a_{\langle \alpha,\delta\rangle}^m)=
a_{\langle\alpha,\delta,0\rangle}^m$,
for $\delta\in \mu_B$, $\cf \delta\ne \omega$,
$f_7^*(a_{\langle \alpha,\delta\rangle}^m)=
a_{\langle \alpha,\delta\rangle}^m$,
$f_8^*(a_{\langle \alpha,\beta\rangle}^m)=a_{\langle \alpha,\beta\rangle}^m$,
$f_8^*(a_{\langle \alpha,\delta,n\rangle}^m)=
a_{\langle \alpha,\delta, n+1\rangle}^m$,
$f_9^*(a_{\langle \alpha,\beta\rangle}^m)=a_{\langle \alpha,\beta\rangle}^m$,
$f_9^*(a_{\langle \alpha,\delta,n\rangle}^m)=
a_{\langle \alpha,\delta_n\rangle}^m$.

Let, as above,
\begin{align*}
B_0&= \langle \{ a_{\langle 0,0\rangle}^m\mid m\in \omega\}\rangle,\\
B_1&= \langle \{ a_{\langle \alpha,0\rangle}^m\mid
\alpha\in \mu,\ m\in \omega\}\rangle,\\
B_2&= \langle \{ a_{\langle 0,\beta\rangle}^m\mid
\beta\in \mu,\ m\in \omega\}\rangle,\\
B_3&= \langle \{ a_{\langle \alpha,\beta\rangle}^m\mid
\alpha,\beta\in \mu,\ m\in \omega\}\rangle,\\
B_4&= \langle \{ a_{\langle \alpha,\beta\rangle}^m,
a_{\langle \alpha,\beta,n\rangle}^m\mid
\alpha,\beta\in \mu,\ n,m\in \omega\}\rangle,\\
B_5&= \langle\{ a_{\langle 0,\beta\rangle}^m,a_{\langle 0,\beta,n\rangle}^m
\mid \beta\in \mu,\ n,m\in \omega\}\rangle,\\
B_6&= \langle \{ a_{\langle 0,\beta\rangle}^m\mid
\beta\in \mu,\ \cf \beta=\omega,\ m\in \omega \}\rangle,\\
B_7&= \langle \{ a_{\langle \alpha,\beta\rangle}^m\mid
\alpha,\beta\in \mu,\ \cf \beta =\omega,\ m\in \omega\}\rangle.
\end{align*}
It is clear that $f_3^*$, $f_4^*$, and $f_9^*$
are projections onto $B_1$, $B_2$, and~$B_3$, respectively.
Let $f_{10}^*$, $f_{11}^*$, $f_{12}^*$, and $f_{13}^*$ be
projections onto $B_4$, $B_5$, $B_6$, and $B_7$, respectively.

All functions which are considered later satisfy the following
formula $\varphi^0(f,\dots)$:
$$
f_0^* f=f f_0^*=f_0^* \logic\land ff_1^*=f_1^* f.
$$
Let us see what this formula means.
Its first part gives us
$f(a_0^m)=ff_0^*(a_0^m)=f_0^*(a_0^m)=a_0^m$ and
$f_0^* f(a_i^m)=f_0^*(a_i^m)=a_0^m$, therefore
$f(a_i^m)=\alpha_1 a_{i_1}^m+\dots+\alpha_k a_{i_k}^m$.
The second part gives
$$
ff_1^*(a_t^m)=f_1^* f(a_t^m)\Leftrightarrow
f(a_t^{m-1})=f_1^* \alpha_1(m)a_{t_1(m)}^m+\dots+\alpha_k(m) a_{t_k(m)}^m)=
\alpha_1(m) a_{t_1(m)}^{m-1}+\dots+\alpha_k(m)a_{t_k(m)}^{m-1},
$$
therefore for every $t\in I^*$ and for all
$l,m\in \omega$ we have $\alpha_i(m)=\alpha_i(l)$ and $t_i(m)=t_i(l)$.

Now we apply Lemma~\ref{shel_l1_3} with
$J=\{ \langle \alpha,\beta\rangle\mid \alpha,\beta \in \mu\}$,
$J_\beta=\{ \langle \alpha,\beta\rangle\mid \alpha\in \mu\}$, $I=I^*$,
and $f=f_{14}^*$.

The formula
\begin{multline*}
\varphi^1(f,g;f_1^*,\dots,f_{14}^*) \prisv
\varphi^0(f) \logic\land \varphi^0(g) \logic\land
\exists h_1\losp \exists h_2\losp
(h_1 f h_1^{-1} =f_2^* \logic\land h_2 gh_2^{-1}=f_2^*)
\\
{}\logic\land f_9^* f=f \logic\land f_9^* g=g \logic\land
\exists h\losp (f f_{14}^*=f_{14}^* h \logic\land
h f_9^* =f_9^* h f_9^* \logic\land hf=g)
\end{multline*}
says that
\begin{enumerate}
\item
$f$ and $g$ are conjugate to~$f_2^*$;
\item
$\Rng f, \Rng g\subset B_3$;
\item
$\exists h\losp
(h\circ f_{14}^*=f_{14}^*\circ h \logic\land
\Rng h|_{B_3}\subseteq B_3 \logic\land h\circ f=g)$.
\end{enumerate}

We shall write $\varphi^1(f,g;\dots)$ also in the form $f\le g$.

If $f$ and $g$ are conjugate to $f_2^*$,
$f(a_{\langle 0,0\rangle}^m)=a_{\langle \alpha,\beta\rangle}^m$ and
$g(a_{\langle 0,0\rangle}^m)=
\tau(a_{\langle \alpha_1, \beta_1\rangle}^m,\dots,
a_{\langle \alpha_k,\beta_k\rangle}^m)$, then $f\le g$
if and only if
$\beta\le \beta_1,\dots,\allowbreak \beta\le \beta_k$.

The formula
\begin{align*}
& \varphi_2(f,f_1^*,\dots,f_{14}^*) \prisv
\varphi^0(f) \logic\land
(ff_4^*=f_9^* ff_4^*) \logic\land
(ff_{11}^*=f_{10}^* ff_{11}^*)
\logic\land
(ff_2^*=f_3^* ff_2^*) \logic\land
(ff_{12}^*=f_{13}^*ff_{12}^*)
\\*
& \quad \logic\land
\forall g\losp
(\varphi^0(g) \logic\land \exists h\losp
(hgh^{-1}=f^*) \logic\land f_4^* g=g
\Rightarrow
\varphi^1(f,fg;f_1^*,\dots,f_{14}^*))
\\*
& \quad \logic\land
(ff_{11}^* f_7^*=f_7^* ff_{11}^*) \logic\land
(ff_{11}^* f_8^*=f_8^* ff_{11}^*) \logic\land
(ff_{11}^* f_9^*=f_9^* ff_{11}^*) \logic\land
(ff_9^*f_3^*=f_3^*ff_9^*)
\end{align*}
says that
\begin{enumerate}
\item
$\Rng f|_{B_2}\subseteq B_3$; $\Rng f|_{B_5}\subseteq B_4$;
$\Rng f|_{B_0}\subseteq B_1$; $\Rng f|_{B_6}\subseteq B_7$;
\item
for any~$g$ satisfying $\varphi^0(g)$ and
conjugate to~$f_2^*$, if $\Rng g\subseteq B_2$, then
$g\le f\circ g$;
\item
$f|_{B_5}$ commutes with $f_7^*$, $f_8^*$, and $f_9^*$;
\item
$f|_{B_3}$ commutes with~$f_3^*$.
\end{enumerate}

Then, similarly to Statement~\ref{shel_s4_3},
we can prove that the formula $\varphi_2(f,\dots)$ holds in $\Endom(B')$
if and only if
for any $\beta\in \mu$
$$
f(a_{\langle 0,\beta\rangle}^m)=
\tau(a_{\langle \alpha_1,\beta\rangle}^m,\dots,
a_{\langle \alpha_k,\beta\rangle}^m)\ \ \text{and}\ \
f(a_{\langle 0,\beta,n\rangle}^m)=
\tau(a_{\langle \alpha_1,\beta,n\rangle}^m,\dots,
a_{\langle \alpha_k,\beta,n\rangle}^m)
$$
for some linear combination~$\tau$ and ordinal numbers
$\alpha_1,\dots,\alpha_k \in \mu$ (which do not depend on~$\beta$).

The formula
$$
\varphi^3(f,f_1^*,\dots,f_{14}^*)\prisv
(ff_4^*=f_9^*ff_4^*) \logic\land \exists f_1\losp
(f_1f_4^*=ff_4^* \logic\land \varphi^2(f_1,f_1^*,\dots,f_{14}^*))
$$
says that
\begin{enumerate}
\item
$\Rng f|_{B_2}\subseteq B_3$;
\item
$\exists f_1 (f_1|_{B_2}=f|_{B_2} \logic\land \varphi_2(f_1))$.
\end{enumerate}

The formula $\varphi_3(f,\dots)$ holds if and only if
$$
f(a_{\langle 0,\beta\rangle}^m)=
\tau(a_{\langle \gamma_1,\beta\rangle}^m,\dots,
a_{\langle \gamma_k,\beta\rangle}^m)
$$
for every $\beta\in \mu$ and some $\tau; \gamma_1,\dots,\gamma_k$.
This follows immediately from Statement~\ref{shel_s4_3}.

Let the formula $\varphi_4(f)$ say that
\begin{enumerate}
\item
$\Rng f|_{B_2}\subseteq B_3$;
\item
$\varphi_3(f)$;
\item
$\forall g\losp (\varphi_3(g)\Rightarrow
g\circ f_5\circ f|_{B_0}=f_5\circ f\circ f_5\circ g|_{B_0} \logic\land
f_5^*\circ f\circ f_5^*\circ f|_{B_0}=f_6^*\circ f|_{B_0} \logic\land
f_2^*\circ f|_{B_0}=f_2^*|_{B_0}$.
\end{enumerate}

The formula $\varphi_4(f)$ holds if and only if
$$
f(a_{\langle 0,\beta\rangle}^m)=
\tau(a_{\langle \gamma_1,\beta\rangle}^m,\dots,
a_{\langle \gamma_k,\beta\rangle}^m),
$$
where $\tau$ is a~beautiful linear combination.

Now we suppose that $\mu_B$ is a~singular cardinal number.

We let $\mu_1< \mu$, where $\mu_1$ is a~regular cardinal and $\mu_1> \omega$.
Let $I^*\setminus J=I_0\cup
\{ \langle \alpha,\delta,n\rangle\mid \alpha\in \mu,\allowbreak\
{\delta\in \mu_1},\allowbreak\
{\cf \delta=\omega},\allowbreak\ n\in \omega\}$, $|I_0|=\mu$.

For every limit ordinal $\delta\in \mu_1$ such that $\cf \delta=\omega$,
similarly to the previous case we shall choose an increasing sequence
$(\delta_n)_{n \in\omega}$ of ordinal numbers less than~$\delta$,
with limit~$\delta$, such that
for any $\beta\in \mu_1$ and $n\in \omega$ the set
$\{ \delta\in \mu_1\mid \beta=\delta_n\}$
is a~stationary subset on~$\mu_1$.

Let
\begin{align*}
B_1&=\langle \{ a_{\langle \alpha,0\rangle}^m\mid
\alpha\in \mu,\ m\in \omega \}\rangle,\\
B_2&= \langle \{ a_{\langle 0,\beta\rangle}^m\mid \beta\in \mu_1,\
m\in \omega\}\rangle,\\
B_3&=\langle \{ a_{\langle\alpha,\beta\rangle}^m\mid \alpha\in \mu,\
\beta\in \mu_1,\ m\in \omega\}\rangle.
\end{align*}

As in the previous case, we can define the functions~$f_i^*$
in such a~way that for some
$\varphi'({\dots})$ the formula $\varphi'(f;\dots)$ holds
$\Endom(B')$ if and only if there exists an ordinal number $\alpha\in \mu$
such that for every $\beta\in \mu_1$ and every $m\in \omega$
$$
f(a_{\langle 0,\beta\rangle}^m)=a_{\langle \alpha,\beta\rangle}^m.
$$

Let the formula $\varphi^1(f,\dots)$ say that
\begin{enumerate}
\item
$\Rng f|_{B_0}\subseteq B_2$;
\item
for every~$g$ we have
$\varphi^0(g)\Rightarrow (f\circ g)|_{B_0}=(g\circ f)|_{B_0}$.
\end{enumerate}

It is easy to check that the formula $\varphi^1(f,\dots)$ holds
if and only if
there exist a~linear combination~$\sigma$ and distinct ordinal numbers
$\beta_1,\dots,\beta_m\in \mu_1$ such that for every $\alpha\in \mu$
and every $m\in \omega$
$$
f(a_{\langle \alpha,0\rangle}^m)=
\sigma(a_{\langle \alpha,\beta_1\rangle}^m,\dots,
a_{\langle \alpha,\beta_m\rangle}^m).
$$

As the cardinal number~$\mu_1$ is regular, we can use the previous case.
Thus, there is a~formula $\varphi^2(f;\dots)$ such that $\varphi^2(f)$
holds in $\Endom(B')$ if and only if there exists
$\beta\in \mu_1$ such that for every $\alpha\in \mu$ and every $m\in \omega$
$$
f(a_{\langle \alpha,0\rangle}^m)=a_{\langle \alpha,\beta\rangle}^m.
$$

Let $\mu=\bigcup\limits_{i\in \cf \mu} \mu_i$, where $\mu_i\in \mu$ and
the sequence $(\mu_i)$ increases. We have just proved that for every
$\gamma\in \cf \mu$ there exists a~function $\bar f_\gamma^*$ such that
\begin{enumerate}
\item
the formula $\varphi^2[f,\bar f_\gamma^*]$ holds in $\Endom(B')$
if and only if
there exists $\beta\in \mu_\gamma^+$ such that for all
$\alpha\in \mu$ and all $m\in \omega$
$$
f(a_{\langle \alpha,0\rangle}^m)=a_{\langle \alpha,\beta\rangle}^m;
$$
\item
${f_{\gamma,0}^*}$ is a~projection onto
$$
\langle \{ a_{\langle \alpha,\beta\rangle}^m\mid
\alpha\in \mu,\ \beta\in \mu_\gamma^+\}\rangle.
$$
\end{enumerate}

Further, there exists a~formula~$\varphi^3$ and a~vector
of functions~$g^*$ such that the formula
$\varphi^3(\bar f,\bar g^*)$ holds if and only if
$\bar f=\bar f_{\gamma}^*$ for some $\gamma\in \mu$.

Let now the formula $\varphi^4(f,\bar g^*)$ say that there
exists~$\bar f_1$ such that $\varphi^3(\bar f_1,\bar g^*)$ and for
every~$\bar f_2$ satisfying the formulas $\varphi^3(\bar f_2,\bar g^*)$
and $\Rng (\bar f_1)_0\subseteq \Rng (\bar f_2)_0$
we have also $\varphi^2(f,\bar f_2)$.
If the formula $\varphi^4(f,g^*)$ is true, then there exists
$\bar f_1=\bar f_{\gamma}^*$ for some $\gamma\in \mu$
and for every $\bar f_2=f_{\lambda}^*$ (where $\lambda\ge \gamma$)
we have the formula
$$
f(a_{\langle \alpha,0\rangle}^m)=a_{\langle \alpha,\beta\rangle}^m,
$$
where $\beta< \mu_{\lambda}^+$.

Let $f$ be such that
$$
f(a_{\langle \alpha,0\rangle}^m)=a_{\langle \alpha,\beta\rangle}^m,\quad
\beta \in \mu.
$$
Then $\beta\in \mu_\gamma^+$ for $\gamma\in \cf \mu$ and so
the formula $\varphi^4(f,g^*)$ is true for some~$g^*$.

Now we only need to consider the formula
$\varphi^4(f_5^*\circ f\circ f_5^*)$, which is the required formula.

Therefore the case~2 is completely studied, in this case we have
(similarly to Sec.~\ref{sec5}) a~formula (with parameters)
which holds for a~set of $\mu_B$
independent projections~$f$
satisfying the formula $\mathrm{Fine}(f)$.
Thus we can suppose that we have a~set of projections onto~$\mu_B$
from independent direct summands of the group~$B$ isomorphic to the group
$$
\bigoplus_{i\in \omega} \mathbb Z(p^{n_i}).
$$
This set will be denoted by~$\mathbf F$.

\smallskip

Case 3. $\forall n\in \omega\losp (\mu_n< \mu_B)$ and $\mu_B> \omega$.
Naturally, in this case the cardinal number~$\mu_B$ is singular
and $\cf \mu_B=\omega$.

Choose in the sequence $(\mu_i)_{i\in \omega}$
an increasing subsequence $(\mu_{n_i})_{i\in \omega}$
with limit~$\mu_B$. For every natural~$i$,
if the number $\mu_{n_i}$ is regular, then by~$\mu^i$
we shall denote $\mu_{n_i}$, and if it is not regular, then by $\mu^i$
we shall denote some regular cardinal number smaller than $\mu_{n_i}$
and greater than~$\omega$ in such a~way that in the result the limit
of the sequence $(\mu^i)_{i\in\omega}$ is also equal to~$\mu_B$.
For every natural~$i$ suppose that we have a~set~$I_i^*$
of power~$\mu^i$ which is the union of the following
sets:
$J_i=\{ \langle \alpha,\beta\rangle\mid \alpha,\beta\in \mu^i\}\cup \{0\}$,
$|I^0_i|=\mu^i$, and
$I_i=\{ \langle \alpha,\delta,n\rangle \mid \alpha\in \mu^i,\
\delta\in \mu^i,\ \cf \delta=\omega,\ n\in \omega\}$.
Let us for every $i\in \omega$ have $\mu^i$~independent
generating direct cyclic summands of order~$p^{n_i}$,
denoted by~$a_t^i$, where $t\in I^*$.

Let $\langle \{ a_t^i\mid t\in I^*,\ i\in \omega\}\rangle =B'$.

Again we need to change the formulation of Theorem~\ref{shel_t4_1},
and Lemma~\ref{shel_l4_2} will be corrected again:
{\itshape there exists a~formula $\varphi(f)$ with one free variable~$f$
such that $\varphi(f)$ holds in $\Endom(B')$ if and only if
there exists a~sequence of ordinal numbers $(\alpha_i)_{i\in \omega}$,
where $\alpha_i\in \mu^i$, such that for every $m\in \omega$ and every
$\beta_m \in \mu^m$}
$$
f(a_{\langle 0,\beta_m\rangle}^m)=a_{\langle \alpha_m,\beta_m\rangle}^m.
$$

All changes in the proof of the theorem with the help of the lemma
are clear, so we shall not write them here.
In the proof of the lemma we need only the third case.
As above, for every limit ordinal $\delta\in \mu^i$ such that
$\cf \delta=\omega$,
we shall choose an increasing subsequence $(\delta_n)_{n\in \omega}$
of ordinal numbers less than~$\delta$, with limit~$\delta$,
such that for any $\beta\in \mu$ and $n\in \omega$
the set $\{ \delta\in \mu\mid \beta=\delta_n\}$
is a~stationary subset in~$\mu$.

We shall again introduce the functions $f_0^*,\dots,f_{14}^*$,
which will differ a~little from the similar functions
from the previous case.
Namely, let $f_1^*(a_t^m)=a_0^m$,
$f_2^*(a_{\langle \alpha,\beta\rangle}^m)=a_{\langle 0,0\rangle}^m$,
$f_3^*(a_{\langle \alpha,\beta\rangle}^m)=a_{\langle \alpha,0\rangle}^m$,
$f_4^*(a_{\langle \alpha,\beta\rangle}^m)=a_{\langle 0,\beta\rangle}^m$,
$f_5^*(a_{\langle \alpha,\beta\rangle}^m)=a_{\langle \beta,\alpha\rangle}^m$,
$f_6^*(a_{\langle \alpha,\beta\rangle}^m)=a_{\langle \alpha,\alpha\rangle}^m$,
for $\delta\in \mu_B$, $\cf \delta=\omega$,
$f_7^*(a_{\langle \alpha,\delta\rangle}^m)=
a_{\langle\alpha,\delta,0\rangle}^m$,
for $\delta\in \mu_B$, $\cf \delta\ne \omega$,
$f_7^*(a_{\langle \alpha,\delta\rangle}^m)=
a_{\langle \alpha,\delta\rangle}^m$,
$f_8^*(a_{\langle \alpha,\beta\rangle}^m)=a_{\langle \alpha,\beta\rangle}^m$,
$f_8^*(a_{\langle \alpha,\delta,n\rangle}^m)=
a_{\langle \alpha,\delta, n+1\rangle}^m$,
$f_9^*(a_{\langle \alpha,\beta\rangle}^m)=a_{\langle \alpha,\beta\rangle}^m$,
$f_9^*(a_{\langle \alpha,\delta,n\rangle}^m)=
a_{\langle \alpha,\delta_n\rangle}^m$.

Let, as above,
\begin{align*}
B_0&= \langle \{ a_{\langle 0,0\rangle}^m\mid m\in \omega\}\rangle,\\
B_1&= \langle \{ a_{\langle \alpha,0\rangle}^m\mid
\alpha\in \mu^m,\ m\in \omega\}\rangle,\\
B_2&= \langle \{ a_{\langle 0,\beta\rangle}^m\mid
\beta\in \mu^m,\ m\in \omega\}\rangle,\\
B_3&= \langle \{ a_{\langle \alpha,\beta\rangle}^m\mid
\alpha,\beta\in \mu^m,\ m\in \omega\}\rangle,\\
B_4&= \langle \{ a_{\langle \alpha,\beta\rangle}^m,
a_{\langle \alpha,\beta,n\rangle}^m\mid \alpha,\beta\in \mu^m,\
n,m\in \omega\}\rangle,\\
B_5&= \langle\{ a_{\langle 0,\beta\rangle}^m,a_{\langle 0,\beta,n\rangle}^m
\mid \beta\in \mu^m,\ n,m\in \omega\}\rangle,\\
B_6&= \langle \{ a_{\langle 0,\beta\rangle}^m\mid
\beta\in \mu^m,\ \cf \beta=\omega,\ m\in \omega \}\rangle,\\
B_7&= \langle \{ a_{\langle \alpha,\beta\rangle}^m\mid
\alpha,\beta\in \mu^m,\ \cf \beta =\omega,\ m\in \omega\}\rangle.
\end{align*}
Let $f_{10}^*$, $f_{11}^*$, $f_{12}^*$, and $f_{13}^*$ be projections onto
$B_4$, $B_5$, $B_6$, and $B_7$, respectively.

All functions which will be considered later satisfy the
following formula $\varphi^0(f,\dots)$:
$$
f_1^* f=f f_1^*=f_1^*.
$$
This formula implies $f(a_0^m)=ff_0^*(a_0^m)=f_0^*(a_0^m)=a_0^m$ and
$f_0^* f(a_i^m)=f_0^*(a_i^m)=a_0^m$, and therefore
$f(a_i^m)=\alpha_1 a_{i_1}^m+\dots+\alpha_k a_{i_k}^m$.

For every $m\in \omega$ we apply Lemma~\ref{shel_l1_3} with
$J^m=\{ \langle \alpha,\beta\rangle\mid \alpha,\beta \in \mu^m\}$,
$J_\beta^m=\{ \langle \alpha,\beta\rangle\mid \alpha\in \mu^m\}$, and
$I^m=I^*_m$. Let us for every $m\in \omega$ have the corresponding
function~$f^m$ on the group
$B^m=\langle \{ a_t^m\mid t\in I^*_m\}\rangle$. Construct
with its help the function~$f_{14}^*$, which coincides with~$f^m$ on every
subgroup~$B^m$.

As above, the formula
\begin{multline*}
f\le g \prisv \varphi^1(f,g;f_1^*,\dots,f_{14}^*) \prisv
\varphi^0(f) \logic\land \varphi^0(g) \logic\land
\exists h_1\losp \exists h_2\losp
(h_1 f h_1^{-1} =f_2^* \logic\land h_2 gh_2^{-1}=f_2^*)
\\
\logic\land f_9^* f=f \logic\land f_9^* g=g \logic\land
\exists h\losp
(f f_{14}^*=f_{14}^* h \logic\land h f_9^* =f_9^* h f_9^* \logic\land hf=g)
\end{multline*}
says that
\begin{enumerate}
\item
$f$ and $g$ are conjugate to~$f_2^*$;
\item
$\Rng f, \Rng g\subset B_3$;
\item
$\exists h\losp
(h\circ f_{14}^*=f_{14}^*\circ h \logic\land
\Rng h|_{B_3}\subseteq B_3 \logic\land h\circ f=g)$.
\end{enumerate}

If $f$ and $g$ are conjugate to~$f_2^*$,
$f(a_{\langle 0,0\rangle}^m)=a_{\langle \alpha^m,\beta^m\rangle}^m$, and
$g(a_{\langle 0,0\rangle}^m)=
\tau^m(a_{\langle \alpha_1^m, \beta_1^m\rangle}^m,\dots,
a_{\langle \alpha_{k_m}^m,\beta_{k_m}^m\rangle}^m)$,
then $f\le g$ if and only if $\beta^m\le \beta_1^m,\dots,\allowbreak
\beta^m \le \beta_{k_m}^m$.

The formula $\varphi_2(f,f_1^*,\dots,f_{14}^*)$ says that
\begin{enumerate}
\item
$\Rng f|_{B_2}\subseteq B_3$; $\Rng f|_{B_5}\subseteq B_4$;
$\Rng f|_{B_0}\subseteq B_1$; $\Rng f|_{B_6}\subseteq B_7$;
\item
for every $g$ satisfying the formula $\varphi^0(g)$ and
conjugate to~$f_2^*$ from $\Rng g\subseteq B_2$ it follows that
$g\le f\circ g$;
\item
$f|_{B_5}$ commutes with $f_7^*$, $f_8^*$, and $f_9^*$;
\item
$f|_{B_3}$ commutes with $f_3^*$.
\end{enumerate}

The formula $\varphi_2(f,\dots)$ holds in $\Endom(B')$
if and only if
for every $m\in \omega$ and every $\beta\in \mu^m$
$$
f(a_{\langle 0,\beta\rangle}^m)=
\tau^m(a_{\langle \alpha_1^m,\beta\rangle}^m,\dots,
a_{\langle \alpha_{k_m}^m,\beta\rangle}^m)\ \ \text{and}\ \
f(a_{\langle 0,\beta,n\rangle}^m)=
\tau^m(a_{\langle \alpha_1^m,\beta,n\rangle}^m,\dots,
a_{\langle \alpha_{k_m}^m,\beta,n\rangle}^m)
$$
for some linear combination~$\tau^m$ and ordinal numbers
$\alpha_1^m,\dots,\alpha_{k_m}^m \in \mu^m$,
which do not depends on~$\beta$.

The formula $\varphi^3(f,f_1^*,\dots,f_{14}^*)$ says that
\begin{enumerate}
\item
$\Rng f|_{B_2}\subseteq B_3$;
\item
$\exists f_1\losp (f_1|_{B_2}=f|_{B_2} \logic\land \varphi_2(f_1))$.
\end{enumerate}

The formula $\varphi_3(f,\dots)$ holds if and only if
for every $m\in \omega$
$$
f(a_{\langle 0,\beta\rangle}^m)=
\tau^m(a_{\langle \gamma_1^m,\beta\rangle}^m,\dots,
a_{\langle \gamma_{k_m}^m,\beta\rangle}^m)
$$
for every $\beta\in \mu^m$ and some
$\tau^m; \gamma_1^m,\dots,\gamma_{k_m}^m$.

The formula $\varphi_4(f)$ says that
\begin{enumerate}
\item
$\Rng f|_{B_2}\subseteq B_3$;
\item
$\varphi_3(f)$;
\item
$\forall g\losp (\varphi_3(g)\Rightarrow
g\circ f_5\circ f|_{B_0}=f_5\circ f\circ f_5\circ g|_{B_0} \logic\land
f_5^*\circ f\circ f_5^*\circ f|_{B_0}=f_6^*\circ f|_{B_0} \logic\land
f_2^*\circ f|_{B_0}=f_2^*|_{B_0}$.
\end{enumerate}

The formula $\varphi_4(f)$ holds if and only if
for every $m\in \omega$
$$
f(a_{\langle 0,\beta\rangle}^m)=
\tau^m(a_{\langle \gamma_1^m,\beta\rangle}^m,\dots,
a_{\langle \gamma_{k_m}^m,\beta\rangle}^m),
$$
where $\tau$ is a~beautiful linear combination, i.e., there exists a~set
$\gamma^1,\dots, \gamma^n,\dots$, where $\gamma^i\in \mu^i$ is such that
$$
f(a_{\langle 0,\beta\rangle}^m)=a_{\langle \gamma^m,\beta\rangle}^m
$$
for all $m\in \omega$ and $\beta\in \mu^m$.

Therefore we suppose that a~set of $\mu_B$ independent projections
satisfying the formula $\mathrm{Fine}(f)$ is fixed.
This set will also be denoted by~$\mathbf F$.
In what follows we shall not distinguish the second and the third cases.

\subsection{Final Rank of the Basic Subgroup Is Uncountable}\label{sec7_4}
Let us change the formula $\mathrm{Fine}(f)$ a~little:
\begin{align*}
& \mathrm{Fine}(f) \prisv
[\forall f'\losp (\Idem(f') \logic\land f'\varphi_B\ne 0\Rightarrow
f'f\ne 0)]
\\
& \quad \logic\land
[\forall f_1\losp \forall f_2\losp
(\Idem^*(f_1) \logic\land \Idem^*(f_2) \logic\land o(f_1)=o(f_2)
\\
& \quad \logic\land \forall c\in Z\losp
(cf_1\ne 0\Rightarrow cf_1f\ne 0) \logic\land \forall c\in Z\losp
(cf_2\ne 0\Rightarrow cff_2\ne 0)\Rightarrow
f_1f_2\ne 0 \logic\land f_2f_1\ne 0)]
\\
& \quad \logic\land
[\forall \rho'\losp (\Idem^*(\rho')\Rightarrow
\exists f'\losp (\Idem^*(f') \logic\land
o(f')> o(\rho') \logic\land \forall c\in Z\losp
(cf'\ne 0\Rightarrow cf'f\ne 0)))]
\\
& \quad \logic\land [\forall f'\losp (\Idem^*(f') \logic\land
\forall c\in Z\losp
(cf'\ne 0\Rightarrow cf'f_B\ne 0)\Rightarrow pf'\ne 0)]
\\
& \quad \logic\land
[\forall \rho'\losp \forall f'\losp
(\Idem^*(f') \logic\land \forall c\in Z\losp
(cf'\ne 0\Rightarrow cf'f_B\ne 0) \logic\land o(f')=o(\rho')
\\
& \quad \Rightarrow \forall f\losp
(\Idem^*(f) \logic\land \forall c\in Z\losp
(cf\ne 0\Rightarrow cff_B\ne 0) \logic\land
o(f)> o(\rho')\Rightarrow o(f)> o(\rho')^2))].
\end{align*}
The first part of the formula, which is enclosed in square brackets,
postulates $fA\subset \varphi_B A=B$. The second part, which is enclosed
in square brackets, states that the image of~$fA$
contains at most one cyclic direct summand of one order.
The third part states that the orders of direct summands of~$fA$
are unbounded.
The fourth part states that
the cyclic summand of~$fA$ of the smallest order has order
greater than~$p$ (i.e., at least not smaller than~$p^2$).
Finally, the fifth part states that
for every direct cyclic summand in~$fA$ of order~$p^k$
the next cyclic summand of greater order has order
greater than~$p^{2k}$.

We shall again write the formula $\mathrm{Ins}(\psi)$
from the case~1 from Sec.~\ref{sec7_4}, which says
that
\begin{enumerate}
\item
for every group $fA$, where $f\in \mathbf F$,
there exists a~direct cyclic summand
of the smallest order such that
the action of~$\psi$ on it is multiplication by~$p$;
\item
for every natural~$i$ and every $f\in \mathbf F$
there exists a~direct cyclic summand
$\langle a_i \rangle\subset fA$ of order~$p^{n_i}$ such that
the action of~$\psi$ on it is multiplication by~$p^i$.
\end{enumerate}

Let us fix some endomorphism~$\Psi$ that satisfies the formula
$\mathrm{Intr}(\Psi)$.

Further, fix an endomorphism $\Gamma\colon B'\to B'$
that for every $f\in \mathbf F$
satisfies the following conditions:
\begin{enumerate}
\item
$f \Gamma f=f \Gamma =\Gamma f$, i.e.,
the endomorphism $\Gamma$ maps $fA$ into $fA$;
\item
$\forall \rho\losp (\Idem^*(\rho) \logic\land \rho f=\rho \logic\land
\forall \rho'\losp
(\Idem^*(\rho') \logic\land \rho' f=\rho'\Rightarrow o(\rho')\ge o(\rho))
\Rightarrow \Gamma \rho =0)$,
i.e., the endomorphism $\Gamma$ maps a~cyclic summand
of the smallest order of~$fA$ into zero;
\item
$\forall \rho_1\losp
(\Idem^*(\rho_1) \logic\land \rho_1 f=\rho_1\Rightarrow
\exists \rho_2\losp
(\Idem^*(\rho_2) \logic\land \rho_2 f=\rho_2 \logic\land
o(\rho_1 )< o(\rho_2) \logic\land
\forall \rho\losp (\Idem^*(\rho) \logic\land {\rho f=\rho}\Rightarrow
\neg (o(\rho)>o(\rho_1) \logic\land o(\rho)< o(\rho_2))) \logic\land
\rho_1 \Gamma \rho_2=\Gamma \rho_2 \logic\land \forall c\in Z\losp
(c\rho_1\ne 0\Rightarrow c\rho_1 \Gamma \rho_2\ne 0)))$;
this means that $\Gamma$ maps every generator~$a_i$ of a~cyclic direct
summand of the group~$fA$ (isomorphic to $\mathbb Z(p^{n_i})$) into
the generator~$a_{i-1}$ of a~cyclic direct summand of~$fA$
(isomorphic to $\mathbb Z(p^{n_{i-1}})$).
\end{enumerate}

This endomorphism gives us a~correspondence between generators of cyclic
summands in the group $fA$, for every $f\in \mathbf F$.
We shall assume that it is fixed.

At first, suppose for simplicity that the final rank of a~basic subgroup
of~$A$ coincides with its rank, and that it is uncountable.
Then $|A|=|B|=\mu$.
We suppose that the set~$\mathbf F$ of $\mu$~independent projections
on direct summands of the group~$B$, isomorphic to
$$
\bigoplus_{i\in \omega} \mathbb Z(p^{n_i}),
$$
where the sequence $(n_i)$ is such that $n_1\ge 2$, $n_{i+1}> 2n_i$,
is fixed. Let us fix $f\in \mathbf F$ and interpret the group~$A$ on~$f$.

Let us consider the set $\Endom_f$ of all homomorphisms
$h\colon fA\to A$ that satisfy the following condition:
\begin{align*}
& \exists \rho (\Idem^*(\rho) \logic\land
\exists c\in Z\losp (\Psi \rho=c\rho) \logic\land \rho f=\rho
\\
& \quad \logic\land
\forall \rho'\losp (\Idem^*(\rho') \logic\land \exists c\in Z\losp
(\Psi \rho'=c\rho') \logic\land \rho' f=\rho' \logic\land
(o(\rho')< o(\rho) \logic\lor o(\rho')> o(\rho))\Rightarrow h\rho'=0)
\\
& \quad \logic\land
\forall c\in Z\losp (\Psi \rho=c\rho\Rightarrow pc h \rho=0)).
\end{align*}
This condition means that
\begin{enumerate}
\item
there exists such $i\in \omega$ that $h(a_k)=0$ for every $k\ne i$;
\item
$o(h(a_i))\le p^i$.
\end{enumerate}

Naturally, two such homomorphisms $h_1$~and~$h_2$ are said to be
equivalent if $h_1(a_i)=h_2(a_j)\ne 0$ for some $i,j\in \omega$.

Therefore two homomorphisms $h_1$~and~$h_2$ from $\Endom_f$
are said to be equivalent if there exists
a~homomorphism~$h$ that satisfies the following two conditions:
\begin{enumerate}
\item
$\forall \rho\losp (\Idem^*(\rho) \logic\land \rho f=\rho \logic\land
\exists c\in Z\losp
(\Psi \rho =c\rho) \logic\land h_1 \rho\ne 0\Rightarrow h\rho =h_1\rho)
\logic\land
\forall \rho\losp (\Idem^*(\rho) \logic\land \rho f=\rho \logic\land
\exists c\in Z\losp
(\Psi \rho =c\rho) \logic\land h_2 \rho\ne 0\Rightarrow h\rho =h_2\rho)$;
this means that the homomorphism~$h$ coincides with~$h_1$ on
the element~$a_i$
that satisfies $h_1(a_i)\ne 0$, and it coincides with~$h_2$ on
the element~$a_j$
that satisfies $h_2(a_j)\ne 0$;
\item
$\forall \rho\losp (\Idem^*(\rho) \logic\land \rho f=\rho \logic\land
\exists c\in Z\losp
(\Psi \rho=c\rho) \logic\land \exists \rho_1\losp \exists \rho_2\losp
(\Idem^*(\rho_1) \logic\land \Idem^*(\rho_2) \logic\land
\rho_1 f=\rho_1 \logic\land \rho_2 f=\rho_2 \logic\land
\exists c_1\in Z\losp
(\Psi \rho_1=c_1\rho_1) \logic\land \exists c_2\in Z\losp
(\Psi \rho_2=c_2\rho_2) \logic\land
h_1\rho_1\ne 0 \logic\land h_2 \rho_2\ne 0 \logic\land
{o(\rho)> o(\rho_1)} \logic\land o(\rho)\le o(\rho_2))\Rightarrow
h\Gamma \rho=h\rho)$;
this means that $h(a_j)=h(a_{j-1})=\dots=h(a_{i+1})=h(a_i)$.
\end{enumerate}

Thus $h_1(a_i)=h(a_i)=h(a_j)=h_2(a_j)$, which is what we need.

If we factorize the set $\Endom_f$ by this equivalence, we shall obtain
the set $\tEnd_f$. A~bijection between this set and the group~$A$
is easy to establish. We only need to introduce addition. Namely,
\begin{multline*}
(h_3=h_1\oplus h_2) \prisv \exists h_1'\losp \exists h_2'\losp
(h_1'\sim h_1 \logic\land h_2'\sim h_2
\\
\logic\land h_3=h_1+h_2 \logic\land
\forall \rho\losp (\Idem^*(\rho) \logic\land \rho f=\rho \logic\land
\forall c\in Z\losp (\Psi \rho =c\rho)\Rightarrow
(h_1 \rho =0\Leftrightarrow h_2\rho=0))).
\end{multline*}

We have interpreted the group~$A$ for every $f\in \mathbf F$,
and so we can prove the main theorem for this case.

\begin{proposition}\label{p7.2}
Suppose that $p$-groups $A_1$~and~$A_2$
are direct sums $D_1\oplus G_1$ and $D_1\oplus G_2$, where
the groups $D_1$~and~$D_2$ are divisible, the groups $G_1$~and~$G_2$
are reduced and unbounded,
$|D_1|\le |G_1|$, $|D_2|\le |G_2|$, $B_1$~and~$B_2$ are basic subgroups
of the groups $A_1$~and~$A_2$, respectively,
and the final ranks of the groups
$B_1$~and~$B_2$ coincide with their ranks and are uncountable.
Then elementary equivalence of the rings $\Endom(A_1)$ and
$\Endom(A_2)$ implies equivalence of the groups $A_1$~and~$A_2$
in the language~$\mathcal L_2$.
\end{proposition}

\begin{proof}
As usual, we shall consider an arbitrary sentence~$\psi$
of the second order group
language and show an algorithm which translates this sentence~$\psi$
to the sentence~$\tilde \psi$ of the first order ring language
such that $\tilde \psi$ holds in $\Endom(A)$ if and only if $\psi$
holds in~$A$.

Consider the formula
\begin{multline*}
\mathrm{Min} (f) \prisv
f\in \mathbf F \logic\land \forall f'\losp
(f'\in \mathbf F \Rightarrow \forall c\losp \forall \rho'\losp
\forall \rho\losp (c\in Z \logic\land \Idem^*(\rho') \logic\land
\Idem^*(\rho)
\\
\logic\land \rho f=\rho \logic\land \rho' f'=\rho' \logic\land
\Psi \rho'=c\rho' \logic\land
\Psi \rho =c\rho\Rightarrow o(\rho)\le o( \rho'))).
\end{multline*}
This formula gives us a~summand $f(A)$ in~$B$
such that for every $i\in \omega$
the number $n_i(f)$ is minimal among all
$n_i(f')$ for $f'\in \mathbf F$.

Consider the formula
\begin{multline*}
\mathrm{Basic}(\Lambda) \prisv
\exists f\losp (f\in \mathbf F \logic\land \mathrm{Min}(f)
\logic\land \forall f'\losp \forall c\losp \forall \rho'\losp
(f'\in \mathbf F \logic\land c\in Z \logic\land \Idem^*(\rho') \logic\land
\rho' f'=\rho'
\\
\logic\land
\Psi \rho'=c\rho'\Rightarrow \exists \rho\losp
(\Idem^*(\rho) \logic\land \rho f=\rho \logic\land
\Psi \rho=c\rho
\logic\land \rho \Lambda \rho' =\Lambda \rho' \logic\land
\forall c'\in Z\losp (c \rho \ne 0\Rightarrow c\rho \Lambda \rho'\ne 0)))).
\end{multline*}
This formula defines an endomorphism~$\Lambda$
that maps $a_t^i$ to~$a_0^i$ for every $i\in \omega$ and every $t\in \mu$.
The corresponding~$f_0$ will be denoted by $f_{\mathrm{min}}$.

Translate the sentence~$\psi$ to the sentence
$$
\tilde \psi \prisv \exists \bar g\losp \exists \Gamma\losp
\exists \Psi\losp \exists \Lambda\losp
\exists f_{\mathrm{min}}\in
\mathbf F\psi'(\bar g,\Gamma,\Psi,\Lambda, f_{\mathrm{min}}),
$$
where the formula $\psi'(\dots)$ is obtained from the sentence~$\psi$
with the help of the following translations of subformulas from~$\psi$:
\begin{enumerate}
\item
the subformula $\forall x$ is translated to the subformula
$\forall x\in \tEnd_{f_{\mathrm{min}}}$;
\item
the subformula $\exists x$ is translated to the subformula
$\exists x\in \tEnd_{f_{\mathrm{min}}}$;
\item
the subformula $\forall P_m(v_1,\dots,v_m)({\dots})$
is translated to the subformula
$$
\forall f_1^P\dots \forall f_m^P\losp
\biggl(\forall g\in \mathbf F\losp
\biggl(\,\bigwedge_{i=1}^m
(f_i^Pg\in \Endom_g)\biggr)\Rightarrow \ldots\biggr);
$$
\item
the subformula $\exists P_m(v_1,\dots,v_m)({\dots})$
is translated to the subformula
$$
\exists f_1^P\dots \exists f_m^P\losp
\biggl(\forall g\in \mathbf F\losp \biggl(\,\bigwedge_{i=1}^m
(f_i^Pg\in \Endom_g)\biggr) \logic\land \ldots\biggr);
$$
\item
the subformula $x_1=x_2$
is translated to the subformula $x_1\sim x_2$;
\item
the subformula $x_1=x_2+x_3$ is translated to the subformula
$x_1\sim x_2\oplus x_3$;
\item
the subformula $P_m(x_1,\dots,x_m)$
is translated to the subformula
$$
\exists g\in \mathbf F\losp
\biggl(\,\bigwedge_{i=1}^m (f_i^Pg)=x_i \Lambda g\biggr).
$$
\end{enumerate}

The rest of the proof is similar to the previous cases.
\end{proof}

Now we shall consider the case where
the final rank of~$B$ is greater
than~$\omega$ and does not coincide with the rank of~$B$.

In this case, $A=G\oplus G'$, where the group~$G$
satisfies the conditions of the
previous proposition, the group~$G'$ is bounded, and its power is
greater than~$|G|$. Let $|G|=\mu$ and $|A|=|G'|=\mu'$.

Let us define, for the group~$G$, the set~$\mathbf F$
from Proposition~\ref{p7.2},
and the set~$\mathbf F'$ of $\mu'$~independent projections on
countably generated subgroups of the group~$G'$, from
Sec.~\ref{sec5}.

The formula
$$
\mathrm{Add}(\varphi) \prisv \forall f'\in \mathbf F'\losp
\exists f\in \mathbf F\losp \forall \rho'\losp
(\Idem^*(\rho') \logic\land \rho' f'=\rho'
\Rightarrow \exists \rho\losp
(\Idem^*(\rho) \logic\land \rho f=\rho \logic\land
\varphi \rho'=\rho \varphi \rho'\ne 0))
$$
defines a~functions from the set~$\mathbf F'$ onto the set~$\mathbf F$.

\begin{proposition}\label{p7.3}
Suppose that $p$-groups $A_1$~and~$A_2$
are direct sums $D_1\oplus G_1$ and $D_1\oplus G_2$, where
the groups $D_1$~and~$D_2$ are divisible, the groups
$G_1$~and~$G_2$ are reduced and unbounded,
$|D_1|< |G_1|$, $|D_2|< |G_2|$, $B_1$~and~$B_2$ are basic subgroups
of the groups $A_1$~and~$A_2$, respectively,
and the final ranks of the groups
$B_1$~and~$B_2$ do not coincide with their ranks and are uncountable.
Then elementary equivalence of the rings $\Endom(A_1)$ and
$\Endom(A_2)$ implies equivalence of the groups $A_1$~and~$A_2$
in the language~$\mathcal L_2$.
\end{proposition}

\begin{proof}
We only write an algorithm of translation.
Let us translate the sentence~$\psi$ to the sentence
$$
\tilde \psi \prisv \exists \bar g\losp \exists \Gamma\losp
\exists \Psi\losp \exists \Lambda\losp
\exists f_{\mathrm{min}} \in \mathbf F\losp
\exists \bar g_1 \dots \exists \bar g_k\losp
\exists \bar g'\losp \exists
\tilde g\in \mathbf F'\psi'(\bar g,\Gamma,\Psi,\Lambda, f_{\mathrm{min}}),
$$
where the formula $\psi'({\dots})$ is obtained from the sentence~$\psi$
with the help of the following translations of subformulas from~$\psi$:
\begin{enumerate}
\item
the subformula $\forall x$ is translated to the subformula
$\forall x\in \tEnd_{f_{\mathrm{min}}}\losp \forall x'\in \tEnd_{\tilde g}$;
\item
the subformula $\exists x$ is translated to the subformula
$\exists x\in \tEnd_{f_{\mathrm{min}}}\losp \exists x'\in \tEnd_{\tilde g}$;
\item
the subformula $\forall P_m(v_1,\dots,v_m)({\dots})$
is translated to the subformula
\begin{multline*}
\forall f_1^P\dots \forall f_m^P\losp
\forall {f_1^P}'\dots \forall {f_m^P}'\losp
\forall \varphi_1^P\dots \forall \varphi_m^P\losp
\biggl(\,\bigwedge_{i=1}^m \mathrm{Add}(\varphi_i^P)
\\
\logic\land \forall g\in \mathbf F\losp
\biggl(\,\bigwedge_{i=1}^m f_i^Pg\in \Endom_g\biggr) \logic\land
\forall g\in \mathbf F'\losp
\biggl(\,\bigwedge_{i=1}^m {f_i^P}'g\in \Endom_g\biggr)\Rightarrow
\ldots\biggr);
\end{multline*}
\item
the subformula $\exists P_m(v_1,\dots,v_m)({\dots})$
is translated to the subformula
\begin{multline*}
\exists f_1^P\dots \exists f_m^P\losp
\exists {f_1^P}'\dots \exists {f_m^P}'\losp
\exists \varphi_1^P\dots \exists \varphi_m^P\losp
\biggl(\,\bigwedge_{i=1}^m \mathrm{Add}(\varphi_i^P)
\\
\logic\land \forall g\in \mathbf F\losp
\biggl(\,\bigwedge_{i=1}^m f_i^Pg\in \Endom_g\biggr)\logic\land
\forall g\in \mathbf F'\losp
\biggl(\,\bigwedge_{i=1}^m {f_i^P}'g\in \Endom_g\biggr) \logic\land
\ldots\biggr);
\end{multline*}
\item
the subformula $x_1=x_2$ is translated to the subformula
$x_1\sim x_2 \logic\land x_1'\sim x_2'$;
\item
the subformula $x_1=x_2+x_3$
is translated to the subformula
$x_1\sim x_2\oplus x_3 \logic\land x_1'\sim x_2'\oplus x_3'$;
\item
the subformula $P_m(x_1,\dots,x_m)$
is translated to the subformula
$$
\exists g\in \mathbf F\losp \exists g'\in \mathbf F'\losp
\biggl(\,\bigwedge_{i=1}^m (\varphi_i^P(g')=g \logic\land
f_i^Pg=x_i \Lambda g \logic\land {f_i^P}'g'=x_i'\tilde g h)\biggr).\qed
$$
\end{enumerate}
\renewcommand{\qed}{}
\end{proof}

\subsection{The Countable Restriction of the Second Order Theory
of the Group in the Case Where the Rank of the Basic Subgroup Is Countable}
\label{sec7_5}
Let the group~$A$ have a~countable basic subgroup~$B$.
As above, we suppose that an endomorphism~$\varphi_B$
with the image~$B$ is fixed.

Further, we suppose that an endomorphism~$f_B$ considered in the
case~1 in Sec.~\ref{sec7_3} and an endomorphism~$\Psi$
satisfying the formula $\mathrm{Ins}(\Psi)$ from the same section
are fixed. As we remember,
$$
B=\varphi_B (A)\supset B'=f_B (A)\cong
\bigoplus_{i\in \omega}\mathbb Z (p^{n_i}),
$$
where $n_0\ge 2$, $n_{i+1}> 2n_i$. Generators of cyclic summands
of $f_B(A)$, where $\Psi(a_i)=p^i a_i$,
will be denoted by~$a_i$ ($i\in \omega$).

Further, as in Sec.~\ref{sec7_4}, we fix a~homomorphism
$\Gamma\colon B'\to B'$ satisfying the conditions
\begin{enumerate}
\item
$\forall f\losp (\Idem^*(f) \logic\land \forall c\in Z\losp
(cf\ne 0\Rightarrow cff_B\ne 0)\Rightarrow f\Gamma f=f\Gamma=\Gamma f)$,
i.e., $f\in \Endom(B')$;
\item
$\forall \rho\losp
(\Idem^*(\rho) \logic\land \forall c\in Z\losp
(c\rho\ne 0\Rightarrow c\rho f_B\ne 0) \logic\land
\forall \rho'\losp (\Idem^*(\rho') \logic\land
\forall c\in Z\losp (c\rho'\ne 0\Rightarrow c\rho' f_B\ne 0)\Rightarrow
o(\rho')\ge o(\rho))\Rightarrow \Gamma \rho=0)$,
i.e., the endomorphism~$\Gamma$ maps a~cyclic summand of the smallest order
(in the group) into~$0$;
\item
$\forall \rho_1\losp
(\Idem^*(\rho_1) \logic\land \forall c\in Z\losp
(\rho_1\ne 0\Rightarrow c\rho_1 f_B\ne 0)\Rightarrow
\exists \rho_2\losp (\Idem^*(\rho_2) \logic\land
\forall c\in Z\losp
({c\rho_2\ne 0} \Rightarrow {c\rho_2 f_B \ne 0}) \logic\land
o(\rho_1)< o(\rho_2) \logic\land
\forall \rho\losp (\Idem^*(\rho) \logic\land
\forall c\in Z\losp (c\rho\ne 0\Rightarrow c\rho f_B \ne 0)\Rightarrow
\neg ({o(\rho)> o(\rho_2)} \logic\land o(\rho)< o(\rho_2))) \logic\land
\rho_1\Gamma=\Gamma \rho_2 =\rho_1\Gamma \rho_2 \logic\land
\forall c\in Z\losp (c\rho_1\ne 0\Rightarrow c\rho_1 \Gamma \rho_2 \ne 0)))$,
i.e., $\Gamma$ maps every generator~$a_i$ of a~cyclic
direct summand of the group~$B'$ (isomorphic to $\mathbb Z(p^{n_i})$)
into the generator~$a_{i-1}$ of a~cyclic direct summand (isomorphic to
$\mathbb Z(p^{n_{i-1}})$).
\end{enumerate}

It is clear that for interpretation of the group~$A$ for our function~$f_B$
we can use the sets $\Endom_{f_B}$ and $\tEnd_{f_B}$
from the previous subsection. Therefore, every element $a\in A$
is mapped to a~class of homomorphisms $g\colon B'\to A$ satisfying
the conditions $g(a_i)=a$ and $g(a_j)=0$
if $j\ne i$, where $i$ is such that $p^i\ge o(a)$.

Now suppose that we want to interpret
some set $X=\{ x_i\}_{i\in I}\subset A$, where $|X|\le \omega$,
in the ring $\Endom(A)$.
It is clear that there exists a~sequence $(k_i)$, $i\in I$,
and a~homomorphism $h\colon B'\to A$ such that
$$
h(a_{k_i})x_i,\quad i\in I,
$$
and $o(x_i)\le p^{k_i}$. The set $\{ x_i\}$ is uniquely defined by
the homomorphism~$h$ and the sequence~$(k_i)$.
Therefore, every
set $\{ x_i\}$ can be mapped to a~pair of endomorphisms
consisting of a~projection onto the subgroup
$\langle \{ a_{k_i}\mid i\in I\}\rangle$
and a~homomorphism~$h$.

Similarly, every $n$-place relation in~$A$
is mapped to a~projection
on $\langle \{ a_{k_i}\mid i\in I\}\rangle$ and
an $n$-tuple of homomorphisms $h_1,\dots,h_n$.

Introduce the formulas
\begin{align*}
& \mathrm{Proj}(\rho) \prisv \forall \rho'\losp
(\Idem^*(\rho') \logic\land \forall c\in Z\losp
(c\rho'\ne 0\Rightarrow c\rho'\rho\ne 0)
\\
& \quad {}\Rightarrow \forall c\in Z\losp
(c\rho'\ne 0\Rightarrow c\rho'f_B\ne 0) \logic\land
\exists \rho''\losp (\Idem^*(\rho'') \logic\land \forall c\in Z\losp
(c\rho''\ne 0\Rightarrow c\rho''\rho\ne0)
\\
& \quad \logic\land o(\rho'')=o(\rho') \logic\land
\exists c\in Z\losp (\Psi \rho''=c\rho'')) \logic\land
\exists c\in Z\losp (\Psi \rho'=c\rho'))
\end{align*}
(a~projection on a~direct summand in~$B'$,
generated by $\{ a_{k_i}\mid i\in I\}$) and
$$
\mathrm{Hom}(h) \prisv \forall \rho'\losp
(\Idem^*(\rho') \logic\land \forall c\in Z\losp
(c\rho'\ne 0\Rightarrow c\rho' f_B\ne 0)\Rightarrow
\exists c\in Z\losp (\Psi \rho'=c\rho'\Rightarrow pch\rho'=0)).
$$

Now we are ready to prove the following proposition.

\begin{proposition}\label{p7.4}
Let $p$-groups $A_1$~and~$A_2$ be unbounded and have countable
basic subgroups. Then $\Endom(A_1)\equiv \Endom(A_2)$
implies $\Th_2^\omega(A_1)=\Th_2^\omega(A_2)$.
\end{proposition}

\begin{proof}
Suppose that we have a~sentence $\psi\in \Th_2^\omega(A_1)$. Then for
every predicate variable $P_n(v_1,\dots,v_n)$ included in~$\psi$,
the set $\{ \langle a_1,\dots,a_n\rangle\in A^n\mid P(a_1,\dots,a_n)\}$
is at most countable.

We shall show a~translation of the sentence~$\psi$ into the first order
sentence $\tilde \psi \in \Th_1(\Endom(A_1))$.

Let us translate the sentence~$\psi$ to the sentence
$$
\tilde \psi=\exists \varphi_B\losp \exists f_B\losp \exists \Psi\losp
\exists \Gamma\losp \psi'(\varphi_B,f_B,\Psi,\Gamma),
$$
where the formula $\psi'({\dots})$ is obtained from the sentence~$\psi$
with the help of the following translations of subformulas from~$\psi$:
\begin{enumerate}
\item
the subformula $\forall x$ is translated to the subformula
$\forall x\in \tEnd_{f_B}$;
\item
the subformula $\exists x$ is translated to the subformula
$\exists x\in \tEnd_{f_B}$;
\item
the subformula $\forall P_m(v_1,\dots,v_m)({\dots})$
is translated to the subformula
$$
\forall \rho^P\losp \forall h_1^P \dots \forall h_m^P\losp
(\mathrm{Proj}(\rho^P) \logic\land \mathrm{Hom}(h_1^P)\land \dots \land
\mathrm{Hom}(h_m^P)\Rightarrow \ldots);
$$
\item
the subformula $\exists P_m(v_1,\dots,v_m)({\dots})$
is translated to the subformula
$$
\exists \rho^P\losp \exists h_1^P \dots \exists h_m^P\losp
(\mathrm{Proj}(\rho^P) \logic\land \mathrm{Hom}(h_1^P)
\land\dots \land \mathrm{Hom}(h_m^P) \logic\land \ldots);
$$
\item
the subformula $x_1=x_2$ is translated to the subformula $x_1\sim x_2$;
\item
the subformula $x_1=x_2+x_3$ is translated to the subformula
$x_1\sim x_2\oplus x_3$;
\item
the subformula $ P_m(x_1,\dots,x_m)$
is translated to the subformula
$$
\exists \rho\losp (\Idem^*(\rho) \logic\land
\rho \rho^P=\rho \logic\land \exists c\in Z\losp
(\Psi\rho =c\rho) \logic\land h_1^P \rho=x_1 \land \dots \land
h_m^P \rho =x_m),
$$
i.e., there exists a~cyclic summand
$\langle a_{k_i}\rangle$ in~$B'$ such that
$i\in I$ and $h_1(a_{k_i})=x_1,\dots,\allowbreak h_m(a_{k_i})=x_m$.
\end{enumerate}
\end{proof}

\subsection[The Final Rank of the Basic Subgroup
Is Equal to $\omega$
and Does Not Coincide\\ with Its Rank]{The Final Rank of the Basic Subgroup
Is Equal to $\boldsymbol{\omega}$
and Does Not Coincide with Its Rank}\label{sec7_6}
As above, we suppose that $A=D\oplus G$, where the group~$D$
is divisible, the group~$G$ is reduced and unbounded,
and $|D|< |G|$. Let a~basic subgroup~$\bar B$ of~$A$ (and of~$G$)
have the form $B\oplus B'$, where $|B'|=|\bar B|$,
$|B|=\omega$, and $B'$ is bounded.

Note that from $|D|< |G|$ it follows that $|D|\le \omega$, because
we assume the continuum-hypothesis, which implies
$|B'|=|\bar B|=\omega_1=2^\omega=c$.

The condition $|D|\le \omega$ means that if the groups $A_1=D_1\oplus G_1$
and $A_2=D_2\oplus G_2$ have the described type and
$\Endom(A_1)\equiv \Endom(A_2)$, then $D_1\cong D_2$, and
in order to prove $A_1\equiv_{\mathcal L_2} A_2$
we only need to prove $G_1\equiv_{\mathcal L_2} G_2$
(with the condition that an endomorphism between $D$ and~$G$ is fixed,
see Proposition~\ref{p7.3}).

Therefore for simplicity of arguments we suppose that the group~$A$
is reduced, i.e., $A=G$ and $D=0$.

We fix an endomorphism~$\varphi_B$ with the image~$B$.
Further, let us fix an endomorphism~$f_B$
with the image~$B'$ which is a~direct summand in~$B$ isomorphic to
$$
\bigoplus_{i\in \omega} \mathbb Z(p^{n_i}),
$$
where $n_0\ge 2$, $n_{i+1}> 2n_i$.

Naturally, we also suppose that endomorphisms $\Psi$~and~$\Gamma$,
described in Secs.\ \ref{sec7_3}~and~\ref{sec7_5}, are fixed.

Let $\rho_1$~and~$\rho_2$ be indecomposable projections onto
cyclic direct summands of the group~$B'$ satisfying the formulas
$\exists c\in Z\losp (\Psi \rho_1=c\rho_1) \logic\land
\exists c\in Z\losp (\Psi \rho_2=c\rho_2)$ and $o(\rho_1)> o(\rho_2)$.
Then we shall write $\gamma\in \langle\Gamma \rangle_{\rho_1,\rho_2}$
if an endomorphism $\gamma$~satisfies the formula
\begin{multline*}
\exists \gamma'\losp
(\gamma' \rho_2=\rho_2 \logic\land
\forall \rho\losp (\Idem^*(\rho) \logic\land
\forall c\in Z\losp (c\rho\ne 0\Rightarrow c\rho f_B\ne 0)
\\
{}\logic\land \exists c\in Z\losp (\Psi \rho =c\rho) \logic\land
o(\rho)\le o(\rho_1 ) \logic\land o(\rho)>o(\rho_2)\Rightarrow
\gamma' \rho=\gamma'\Gamma \rho)
\logic\land \gamma \rho_1 =\gamma' \rho_1).
\end{multline*}
This formula means that there exists an
endomorphism $\gamma'\colon B'\to B'$ such that,
if $a_i$ is a~generator in~$\rho_2 A$ and $a_{i+k}$ is a~generator
in~$\rho_1 A$, then we have
\begin{enumerate}
\item
$\gamma'(a_i)=a_i$;
\item
$\forall i\in \{ 1,\dots,k\}\losp
\gamma'(a_{i+k})=\gamma'(a_{i+k-1})=\dots=\gamma'(a_{i+1})=\gamma'(a_i)=a_i$.
\end{enumerate}

Further, $\gamma(a_{i+k})=\gamma'(a_{i+k})=a_i$. Thus, we have
$\gamma (a_{i+k})=a_i$.

Now consider the formula
\begin{align*}
& \mathrm{Onto}(\Lambda)\prisv
\biggl[\forall \bar \rho\losp \forall \bar c\losp
\biggl(\Idem^*(\bar \rho) \logic\land
\forall c\in Z\losp
(c\bar \rho \ne 0\Rightarrow c\bar \rho \varphi_B \ne 0)
\\
& \quad \logic\land \bar c \in Z \logic\land \bar c \bar \rho \ne 0
\Rightarrow \exists \rho_1 \dots \exists \rho_{p-1}\losp
\biggl(\biggl(\,\bigwedge_{i,j=1}^{p-1} \rho_i\rho_j=\rho_j\rho_i=0\biggr)
\logic\land \Idem^*(\rho_1)\land \dots \land \Idem^*(\rho_{p-1})
\\
& \quad \logic\land
\biggl(\,\bigwedge_{i=1}^{p-1} \forall c\in Z\losp
(c\rho_i \ne 0\Rightarrow c\rho_i f_B\ne 0)\biggr) \logic\land
\biggl(\,\bigwedge_{i=1}^{p-1} \exists c\in Z\losp
(\Psi \rho_i=c\rho_i) \biggr)
\\
& \quad \logic\land
\biggl(\,\bigwedge_{i=1}^{p-1} \bar c\bar \rho \Lambda \rho_i\ne 0\biggr)
\logic\land
\biggl(\,\bigwedge_{i=1}^{p-1} \forall c\in Z\losp
(c\bar c\bar \rho\ne 0\Rightarrow
c\bar c \bar \rho\Lambda \rho_i\ne 0)\biggr)
\\
& \quad \logic\land
\biggl(\,\bigwedge_{i,j=1;\ i\ne j}^{p-1}
o(\rho_i)> o(\rho_j)\Rightarrow
\forall \gamma\in \langle \Gamma\rangle_{\rho_i,\rho_j}
(\Lambda \gamma \rho_i\ne \Lambda \rho_i)
\biggl)\!\biggl)\!\biggl)\!\biggr]
\\
& \quad \logic\land
[\Hom(\Lambda)] \logic\land
[\forall \rho_1\losp \forall \rho_2\losp
(
\Idem^*(\rho_1) \logic\land \Idem^*(\rho_2)
\\
& \quad \logic\land
\forall c\in Z\losp (c\rho_1 \ne 0\Rightarrow c\rho_1f_B\ne 0) \logic\land
\forall c\in Z\losp (c\rho_2\ne 0\Rightarrow c\rho_2 f_B\ne 0)
\\
& \quad \logic\land \exists c\in Z\losp
(\Psi \rho_1=c\rho_1) \logic\land
\exists c\in Z\losp (\Psi \rho_2=c\rho_2) \logic\land
o(\rho_1)\ge o(\rho_2)
\\
& \quad \Rightarrow \exists \rho_3 (\Idem^*(\rho_3)\logic\land
\forall c\in Z\losp (c\rho_3\ne 0 \Rightarrow c\rho_3 f_B\ne 0) \logic\land
\exists c\in Z\losp (\Psi \rho_3=c\rho_3)
\\
& \quad \logic\land
((o(\rho_3)> o(\rho_1) \logic\land o(\rho_3)> o(\rho_2) \logic\land
\exists \gamma_1\in \langle \Gamma\rangle_{\rho_3,\rho_1}\losp
\exists \gamma_2 \in \langle \Gamma\rangle_{\rho_3,\rho_2}\losp
(\Lambda \rho_3=\Lambda \gamma_1 \rho_3+\Lambda \gamma_2 \rho_3))
\\
& \quad \logic\lor (o(\rho_3)< o(\rho_1) \logic\land
o(\rho_3)< o(\rho_2) \logic\land
\exists \gamma_2\in \langle \Gamma\rangle_{\rho_1,\rho_2}\losp
\exists \gamma_3 \in \langle \Gamma\rangle_{\rho_1,\rho_3}\losp
(\Lambda \gamma_3\rho_1=\Lambda \rho_1+\Lambda \gamma_2 \rho_1))
\\
& \quad \logic\lor (o(\rho_3)< o(\rho_1) \logic\land
o(\rho_3)> o(\rho_2) \logic\land
\exists \gamma_2\in \langle \Gamma\rangle_{\rho_1,\rho_2}\losp
\exists \gamma_3 \in \langle \Gamma\rangle_{\rho_1,\rho_3}\losp
(\Lambda \gamma_3\rho_1=\Lambda \gamma_1 +\Lambda \gamma_2 \rho_1)))))]
\\
& \quad \logic\land
[\forall \rho_1\rho_2 (\Idem^*(\rho_1) \logic\land \Idem^*(\rho_2)
\logic\land
\forall c\in Z\losp (c\rho_1\ne 0\Rightarrow c\rho_1f_B\ne 0)
\\
& \quad \logic\land
\forall c\in Z\losp (c\rho_2\ne 0\Rightarrow c\rho_2 f_B\ne 0) \logic\land
\exists c\in Z\losp (\Psi \rho_1=c\rho_1) \logic\land
\exists c\in Z\losp (\Psi \rho_2=c\rho_2)
\\
& \quad \logic\land o(\rho_1)> o(\rho_2)\Rightarrow
\exists \rho\losp (\Idem^*(\rho) \logic\land
\forall c\in Z\losp (c\rho\ne 0\Rightarrow c\rho \varphi_B\ne 0)
\logic\land \rho\Lambda \rho_2\ne 0)
\\
& \quad \logic\land
\forall \gamma\in \langle \Gamma \rangle_{\rho_1,\rho_2}\losp
(\Lambda \rho_1\ne \Lambda \gamma \rho_1))].
\end{align*}

The first condition in brackets means that for every generator~$b$
of the cyclic direct summand $\langle b\rangle$ of the group~$B$
and for every integer $p$-adic number~$c$
there exist at least $p-1$ numbers $i_1,\dots,i_{p-1}\in \omega$
such that $\Lambda (a_{i_k})=\xi_k c\cdot b$,
where $\xi_k$ are distinct and not divisible by~$p$,
$\xi_k c\cdot b\ne \xi_l c\cdot b$ for all $k\ne l$.

The following conditions mean that for any $a_i$~and~$a_j$
there exists~$a_k$ such that $\Lambda(a_k)=\Lambda(a_)+\Lambda(a_j)$.
Therefore the endomorphism~$\Lambda$ is an epimorphism $B'\to B$.

The last condition in brackets means that $\Lambda$ induces a~bijection
between the set $\{ a_i\mid i\in \omega\}$ and the group~$B$,
with the condition $o(\Lambda(a_i))\le p^i$.

Let us fix an endomorphism~$\Lambda$.

Now recall that we have a~bounded group~$B'$ of some uncountable
power~$\mu$ that has a~definable set $\mathbf F=\mathbf F(\bar g)$
consisting of $\mu$~independent indecomposable projections onto direct
summands of the group~$B'$, and a~set $\mathbf F'=\mathbf F'(\bar g')$
consisting of $\mu$~independent projections onto countably generated direct
summands of the group~$B'$, and for every $f\in \mathbf F'$ the set
consisting of $f_t\in \mathbf F$ such that $f_tA$ is a~direct summand
in~$fA$ is countable.
Denote the subset $\{ f_t\in \mathbf F\mid f_tA\subset fA\}$ of~$\mathbf F$
by~$\mathbf F_f$.
It is clear that the set $\mathbf F_f$ is definable.

Let us fix endomorphisms $\Pi_1$ and~$\Pi_2$
introducing the order on the set
$f(A)$, where $f\in \mathbf F$:
\begin{align*}
& \mathrm{Order}(\Pi_1,\Pi_2) \prisv
\forall f\in \mathbf F\losp
(
\exists ! f_0\in \mathbf F_f\losp
(
\Pi_2 f_0=0 \logic\land{}
\\
& \quad \logic\land
\forall f_1\losp
(
f_1\in \mathbf F_f \logic\land f_1\ne f_0 \Rightarrow
\exists f_2\in \mathbf F_f\losp
(f_1\ne f_2 \logic\land{}
\\
& \quad \logic\land f_2 \Pi_2=\Pi_2f_1=f_2\Pi_2 f_1 \logic\land
\forall c\in Z\losp (cf_1\ne 0\Rightarrow c\Pi_2 f_1\ne 0) \logic\land{}
\\
& \quad \logic\land f_1\Pi_1 =\Pi_1 f_2=f_1\Pi_1f_2 \logic\land
\forall c\in Z\losp
(cf_2\ne 0\Rightarrow c\Pi_1 f_2\ne 0)) \logic\land{}
\\
& \quad \logic\land
\forall f_1\losp \forall f_2\in \mathbf F_f\losp
\biggl(
f_1\ne f_2\Rightarrow \bigwedge_{i=1}^2 \forall f_3,f_4\in \mathbf F_f
\\
& \quad
(f_3\Pi_i f_1=f_3\Pi_i=\Pi f_1\ne 0 \logic\land
f_4\Pi_i f_2=f_4\Pi_i f_2=f_4\Pi_i=
\Pi_i f_2\ne 0\Rightarrow f_3\ne f_4)\biggr) \logic\land
\\
& \quad \logic\land \forall f_1\in \mathbf F_f\losp
\exists f_2\in \mathbf F_f\losp
(\Pi_2 f_2=f_1\Pi_2=f_1\Pi_2 f_2 \ne 0) \logic\land
\forall f_1 \in \mathbf F_f\losp
(\Pi_2\Pi_1 f=f \logic\land{}
\\
& \quad \logic\land \forall f'\losp
(\Idem^*(f') \logic\land f'f=f' \logic\land \Fin(f') \logic\land
f'\ne 0 \Rightarrow{}
\\
& \quad \Rightarrow \exists f_1,f_2\in \mathbf F_f\losp
(f_1 f=f_1 \logic\land f_2 f=f_2 \logic\land
f_2\Pi_2=f_2\Pi_2 f_1=\Pi_2 f_1\ne 0)))
))).
\end{align*}
Up to the order~$\Pi_1$ we can suppose that for every $f_t\in \mathbf F$
($t\in \mu$) we have a~basis in $f_t(A)$ consisting of $f_t^0,f_t^1,\dots$,
where $\Pi_1(f_t^iA)=f_t^{i+1}A$.

Now let us write the conditions for the homomorphism~$\Delta$.

1.
For every function $f\in \mathbf F$ we shall introduce (by a~formula)
a~function
$$
\rho_{\mathrm{on},f}\colon B'\to f_tA
$$
such that
$$
\rho_{\mathrm{on},f}(a_i)=f_t^i
$$
for all natural~$i$.

The condition for the homomorphism~$\Delta$ now has the form
$$
\forall f\in \mathbf F\losp \exists \beta\losp
(\Hom(\beta) \logic\land (1) \logic\land (2) \logic\land (3) \logic\land (4)),
$$
where

(1) $\forall \rho\losp
(\Idem^* (\rho) \logic\land \forall c\in Z\losp
(c\rho\ne 0\Rightarrow c\rho f_B\ne 0) \logic\land
\exists c\in Z\losp (\Psi \rho = c \rho)\Rightarrow
\forall c\in Z\losp (\Psi \rho=c\rho \Rightarrow
c\beta \rho\ne 0 \logic\land pc\beta \rho \ne 0))$,
and this means that $o(\beta(a_i))=p^i$;

(2) $\forall \rho\losp \forall \rho'\losp
(\Idem^*(\rho) \logic\land \Idem^*(\rho') \logic\land
\forall c\in Z\losp (c\rho\ne 0\Rightarrow c\rho f_B\ne 0) \logic\land
\forall c\in Z\losp (c\rho'\ne 0\Rightarrow c\rho' \varphi_B\ne 0)\allowbreak
\exists c\in Z\losp (\Psi \rho =c\rho)\Rightarrow
\neg (\rho' \beta \rho=\beta \rho))$, and this means that
$\beta (a_i)\notin B$;

(3) $\forall \rho\losp (\Idem^*(\rho) \logic\land
\forall c\in Z\losp (c\rho \ne 0\Rightarrow c\rho f_B\ne 0) \logic\land
\exists c\in Z\losp (\Psi \rho =c\rho)\Rightarrow p\beta \rho =
p\beta (\Gamma \rho)+\Lambda \Delta \rho_{\mathrm{on},f} \rho)$,
and this means that $p\beta(a_i)=\rho (a_{i-1})+\Lambda \Delta (f_t^i)$;

(4) $\forall c\in Z\losp
(c\Lambda \Delta \rho_{\mathrm{on},f} \rho \ne 0\Rightarrow
c\beta \rho \ne 0)$, and this means that $o(\Lambda \Delta (f_t^i))\le p^i$.

This condition means that for every $f_t$, where $t\in \mu$, $\Delta$ is
a~mapping
from $\{ f_t^i\mid i\in \omega\}$ into $\{ a_i\mid i\in \omega\}$
such that there is a~sequence $c_{1,t},\dots,c_{m,t},\dots\in A$
such that $c_{i,t}\notin B$,
$o(c_{i,t})=p^i$, $pc_{i,t}=c_{i-1,t}+b_{i,t}$,
where $b_{i,t}\in B$, $b_{i,t}=\beta(\Delta(f_t^i))$,
and $o(b_{i,t})\le p^i$.

As we know, such a~sequence can be considered as a~sequence
of the elements of a~quasibasis of~$A$, and it can be assumed to be
uniquely defined with the help of the sequence
$$
(b_{i,t}\mid i\in \omega, b_{i,t}=\beta(\Delta(f_t^iA))).
$$
An endomorphism~$\beta$ is uniquely defined by every $f\in \mathbf F$,
for abbreviation we shall write $\mathrm{Quasi}_f(\beta)$.

2. Every set of elements from~$\Gamma$
gives us a~set of sequences of elements of the quasibasis
$$
\{ f_t\mid t\in J\}\leftrightarrow
\{ (c_{1,t},\dots,c_{i,t},\dots)\mid t\in J\}=C_J,
$$
and the set $C_J$ gives us a~linear space
$\bar C_J=\langle C_J\rangle$
such that we can determine whether some sequence belongs
to this linear space or not.
We do not want to write the formulas, because they are too complicated,
but we shall write only conditions:
(a)~every existing sequence belongs to~$\bar C_I$;
(b)~linear spaces $\langle C_{J_1}\rangle$ and $\langle C_{J_2}\rangle$
($J_1\cap J_2=\varnothing$) intersect by~$B$.

Therefore, our homomorphism $\Delta$ is a~bijection between
the set $\{f_t^i\mid t\in I\}$ and the quasibasis
$\{ c_{it}\mid i\in \omega,\ t\in \mu\}$
of the group~$A$. We assume that $\Delta$ is fixed.

Now we can interpret the second order theory of the group~$A$.
We shall do the following.

In addition to the set~$\mathbf F$, consisting of
$\mu$~independent projections onto countably generated
direct summands of~$B$,
we shall also consider a~similar set~$\mathbf G$ with the only
additional condition that for every $g\in \mathbf G$
we shall select one fixed projection~$g_0$ onto a~countably generated
direct summand $g_0A$ of~$gA$ such that if $gA=g_0A\oplus A'$,
then $A'$ is also countably generated.

Let us fix some $g\in \mathbf G$ and consider a~homomorphism~$h$
such that $h(g_0A)\in \langle a_i\rangle$ for some $i\in \omega$,
and if $g_t\in \mathbf G_g\setminus \mathbf G_g^0$, then either
$h(g_tA)\subset fA$ for some $f\in \mathbf F$ or $h(g)=0$.

Let the image~$h$ under~$gA$ be finitely generated and the inverse image
of every $f\in \mathbf F$ contain at most $p-1$ elements
of~$\mathbf G_g$.

Then every such~$h$ is mapped to an element from~$A$ as follows:
if $h$ on $\mathbf G_g\setminus \mathbf G_g^0$ is a~finite subset
$\mathcal F$ of $\{ g_t^i\mid t\in \mu,\ i\in \omega\}$,
and $h(g_t^0)=a_k$, then we obtain an element
$$
\sum_{f_t^i\in \mathcal F} \alpha_{ti}c_{it}+b,
$$
where $\alpha_{it}$ is the multiplicity of the inverse image of~$f_t^i$
and $b=\Lambda(a_k)$.

It is clear that, as above, for such a~mapping~$h$ we can
write $h\in \Endom_g$. Two elements
$h_1,h_2\in \Endom_g$
are said to be equivalent if there exists an automorphism~$\alpha$
of the group~$gA$ substituting elements of~$\mathbf G_g$ and replacing
$g_0\in \mathbf G_g^0$ such that $h_1\alpha$ and $h_2$ coincide.
As usual, the set
$\Endom_g$ factorized by such an equivalence
is denoted by $\tEnd_g$.
Addition on the set $\tEnd_g$ is obvious:
the images of an element $g_0$ are added, the number of inverse images
of every~$f_{it}$ is added, and, if it is greater than~$p-1$,
then we have an excess of~$p$
of inverse images of~$f_{it}$, and we add one more inverse image
of $f_{i-1,t}$ and we add
the known $b_{it}=\beta(\Delta (f_{it}))$ to the image of~$g_0$.

The rest of the proof is similar to the previous cases, because
we have $\mu$~independent elements $g_t \in \mathbf G$
and for each of them we can interpret the theory $\Th(A)$.

\section{The Main Theorem}\label{sec8}

We recall (see Sec.~\ref{sec4_2}) that if
$A=D\oplus G$, where $D$ is divisible and $G$
is reduced, then the \emph{expressible rank of the group~$A$} is
the cardinal number
$$
r_{\mathrm{exp}}=\mu=\max(\mu_D,\mu_G),
$$
where $\mu_D$ is the rank of~$D$, and $\mu_G$ is the rank of a~basic
subgroup of~$G$.

{\bf Theorem 1.}
\emph{For any infinite $p$-groups $A_1$~and~$A_2$
elementary equivalence of endomorphism rings $\Endom(A_1)$
and $\Endom(A_2)$ implies coincidence of the second order theories
$\Th_2^{r_{\mathrm{exp}}(A_1)}(A_1)$ and
$\Th_2^{r_{\mathrm{exp}}(A_2)}(A_2)$ of the groups $A_1$~and~$A_2$,
bounded by the cardinal numbers
$r_{\mathrm{exp}}(A_1)$ and $r_{\mathrm{exp}}(A_2)$, respectively.
}
\smallskip
\begin{proof}
Since the rings $\Endom(A_1)$ and $\Endom(A_2)$ are elementary equivalent,
they satisfy the same first order sentences.
If in the ring $\Endom(A_1)$ for some natural~$k$
the sentence
$$
\forall x\losp (p^k x=0) \logic\land \exists x\losp (p^{k-1} x\ne 0)
$$
holds, then the group~$A_1$ is bounded and the maximum
of the orders of its elements
is equal to~$p^k$. It is clear that in this case the same holds also
for the group~$A_2$, and the theorem follows from
Proposition~\ref{p5.1} (Sec.~\ref{sec5_4}).

Now suppose that neither the group~$A_1$, nor the group~$A_2$ is bounded.
Let for some natural~$k$ the sentence $\psi_{p^k}$
(from Sec.~\ref{sec4_4}) hold in the ring $\Endom(A_1)$.

Then this sentence holds also in the ring $\Endom(A_2)$,
and therefore the groups $A_1$ and $A_2$ are direct sums
$D_1\oplus G_1$ and $D_2\oplus G_2$, respectively, where the groups
$D_1$~and~$D_2$ are divisible and the groups $G_1$~and~$G_2$
are bounded by the number~$p^k$.
Further, the sentence $\psi_{p^k}$ fixes the projections
$\rho_D$~and~$\rho_G$ on the groups $D$~and~$G$, respectively.
If $\rho_G=0$, then the groups $A_1$~and~$A_2$ are divisible
and in this case the theorem follows from Proposition~\ref{p6.2}.

Let the rings $\Endom(A_1)$ and $\Endom(A_2)$ satisfy the sentence
\begin{align*}
& \tilde \psi_{p^k}^2 \prisv
\exists \rho_D\losp \exists \rho_G\losp
(
\psi_{p^k}(\rho_D,\rho_G) \logic\land
\exists h\losp (\rho_D h\rho_G=h\rho_G
\\
& \quad \logic\land \forall \rho_1\losp\forall \rho_2\losp
(\Idem^*(\rho_1) \logic\land \Idem^*(\rho_2) \logic\land
\rho_1\rho_G=\rho_1 \logic\land \rho_2\rho_G=\rho_2 \logic\land
\rho_1\rho_2=\rho_2\rho_1=0
\\
& \quad \Rightarrow \exists \rho_1'\losp \exists \rho_2'\losp
(\Idem^*(\rho_1') \logic\land \Idem^*(\rho_2') \logic\land
\rho_1'\rho_D =\rho_1' \logic\land \rho_2'\rho_D=\rho_2'
\\
& \quad \logic\land \rho_1'\rho_2'=\rho_2'\rho_1'=0 \logic\land
\rho_1'h\rho_1=h\rho_1\ne 0 \logic\land \rho_2'h\rho_2=h\rho_2\ne 0)))).
\end{align*}
This sentence (in addition to the conditions of the sentence $\psi_{p^k}$)
says that there exists an endomorphism~$h$ of the group~$A$
such that $h$ maps~$G$ into~$D$ and any two independent cyclic summands
$\rho_1A$ and $\rho_2A$ of the group~$G$ are mapped into independent
quasicyclic summands $\rho_1'A$ and $\rho_2'A$ of~$D$, i.e.,
there exists an embedding of~$G$ into the group~$D$.
This implies that $|G|\le |D|$, i.e., if the sentence $\psi_{p^k}^2$ holds
in $\Endom(A_1)$ and $\Endom(A_2)$, then the groups $A_1$~and~$A_2$
are isomorphic to direct sums $D_1\oplus G_1$ and $D_2\oplus G_2$,
where $|D_1|\ge |G_1|$ and $|D_2|\ge |G_2|$. In this case, the theorem
follows from Proposition~\ref{p6.3}.

If the sentence $\psi_{p^k}$ holds in $\Endom(A_1)$ and $\Endom(A_2)$,
but the sentence $\psi_{p^k}^2$ is false, then the groups $A_1$~and~$A_2$
are direct sums $D_1\oplus G_1$ and $D_2\oplus G_2$, where $|D_1|< |G_1|$
and
$|D_2|< |G_2|$. In this case, the theorem follows from
Proposition~\ref{p6.4}.

If for no natural~$k$ the sentence $\psi_{p^k}$ belongs to
the theory $\Th(\Endom(A_1))$, then for no natural~$k$ the sentence
$\psi_{p^k}$ belongs to the theory $\Th(\Endom(A_2))$,
and therefore both groups $A_1$~and~$A_2$
have unbounded basic subgroups.

Let us consider the formula
\begingroup
\setlength{\multlinegap}{0pt}
\begin{multline*}
\psi(\rho_D,\rho_G) \prisv \Idem(\rho_D) \logic\land
\Idem(\rho_G) \logic\land (\rho_D\rho_G=\rho_G\rho_D=0) \logic\land
(\rho_D+\rho_G=1)
\\*
\logic\land \forall x\losp (\rho_D x\rho_D =0 \logic\lor
p(\rho_D x \rho_D)\ne 0) \logic\land
\forall \rho'\losp (\Idem^*(\rho') \logic\land
\rho'\rho_G =\rho'\Rightarrow \neg
(\forall x\losp (\rho'x\rho'=0 \logic\lor p(\rho' x\rho')\ne 0))).
\end{multline*}
\endgroup%
This formula says that $A$ is the direct sum of its subgroups
$\rho_DA$ and $\rho_GA$,
the group $\rho_DA$ is divisible, and the group $\rho_GA$ is reduced.

Let us consider the sentence
$$
\psi^2 \prisv \exists \rho_D\losp \exists \rho_G\losp \exists h\losp
(\psi(\rho_D,\rho_G) \logic\land \rho_D h\rho_G=h\rho_G
\logic\land \forall \rho\losp
(\Idem^*(\rho) \logic\land \rho \rho_G=\rho\Rightarrow
\forall c\in Z\losp (c\rho\ne 0\Rightarrow ch\rho\ne 0))).
$$
This sentence says that~$A$ is the direct sum of
a~divisible subgroup $D=\rho_DA$ and a~reduced subgroup $G=\rho_GA$, and
there exists an embedding $h\colon G\to D$. Therefore $|G|\le |D|$.
If the rings $\Endom(A_1)$ and $\Endom(A_2)$ satisfy the sentence~$\psi^2$,
then for the groups $A_1$~and~$A_2$ we have $A_1=D_1\oplus G_1$,
$A_2=D_2\oplus G_2$, $|D_1|\ge |G_1|$, $|D_2|\ge |G_2|$,
and the theorem follows from Proposition~\ref{p7.1}.

Now suppose that the rings $\Endom(A_1)$ and $\Endom(A_2)$
do not satisfy the sentence~$\psi^2$.
In this case, $A_1=D_1\oplus G_1$, $A_2=D_2\oplus G_2$, $|D_1|< |G_1|$, and
$|D_2|< |G_2|$.

Recall the formulas from Sec.~\ref{sec7_2}.

Let us consider the sentence
\begin{align*}
& \psi^3\prisv \exists \rho_D\losp \exists \rho_G\losp
(\psi(\rho_D,\rho_G) \logic\land \neg \psi^2
\logic\land \exists \varphi_B\losp
(\Base(\varphi_B) \logic\land \forall \rho\losp
(\Idem^*(\rho)\Rightarrow \exists \rho'\losp
(\Idem^*(\rho') \logic\land o(\rho')> o(\rho)
\\
& \quad \logic\land \exists f\losp
(\mathrm{Ord}_{\rho'}(f) \logic\land \forall f'\losp
(\Idem^*(f') \logic\land f'f=f'\Rightarrow \forall c\in Z\losp
(cf'\ne 0\Rightarrow cf'\varphi_B\ne 0))
\\
& \quad \logic\land \exists h\losp
(\forall f_1\losp (\Idem^*(f_1) \logic\land \forall c\in Z\losp
(cf_1\ne 0\Rightarrow cf_1\varphi_B\ne 0)
\\
& \quad \Rightarrow \exists f_2\losp
(\Idem^*(f_2) \logic\land f_2f=f_2 \logic\land
f_1h=hf_2=f_1hf_2\ne 0)))))))).
\end{align*}
This sentence says that
\begin{enumerate}
\item
$A=\rho_DA\oplus \rho_GA=D\oplus G$, where $D$ is divisible,
$G$~is reduced, and $|D|< |G|$;
\item
$\varphi_B$ is an endomorphism with the image $\varphi_B(A)$
coinciding with some basic subgroup~$B$;
\item
for every natural~$k$ there exists a~natural~$n$ such that
in the group~$B$ there exists a~direct summand (which is a~sum
of cyclic groups of order~$p^n$) having the same power
as the group~$B$.
\end{enumerate}

Therefore the sentence~$\psi^3$ says that the final rank of a~basic
subgroup of the group~$G$ coincide with its rank.

The sentence
\begin{multline*}
\psi^4 \prisv \exists \varphi_B\losp
(\Base(\varphi_B) \logic\land \forall \rho\losp
(\Idem^*(\rho)
\Rightarrow \forall f\losp
(\mathrm{Ord}_\rho(f)\Rightarrow \Fin(f) \logic\lor
\exists h\losp (\forall f_1\losp (\Idem^*(f_1)
\\
\logic\land \forall c\in Z\losp
(cf_1\ne 0\Rightarrow cf_1 \varphi_B\ne 0)\Rightarrow
\exists f_2\losp (\Idem^*(f_2)
\logic\land f_2f=f_2 \logic\land f_1h=hf_2=f_1hf_2\ne 0))))))
\end{multline*}
means that in a~basic subgroup~$B$ for every natural~$n$ every direct
summand which is a~direct sum of cyclic groups of order~$p^n$ either
is finite or has the same power as the group~$B$, so
the group~$B$ is countable.

Therefore if the rings $\Endom(A_1)$ and $\Endom(A_2)$ satisfy
the sentence $\psi^3 \logic\land \neg \psi^4$,
then the groups $A_1$~and~$A_2$
are direct sums $D_1\oplus G_1$ and $D_2\oplus G_2$, where $|D_1|< |G_1|$,
$|D_1|< |G_2|$, and the final ranks of basic subgroups of $A_1$~and~$A_2$
coincide with their ranks and are uncountable. In this case
the theorem follows from Proposition~\ref{p7.1}. If the rings
$\Endom(A_1)$ and $\Endom(A_2)$ satisfy the sentence
$\psi^3 \logic\land \psi^4$,
then their basic subgroups
are countable and in this case the theorem follows
from Proposition~\ref{p7.4}.
Now we have only two cases, and to separate them
we shall write the sentence
\begin{align*}
& \psi^5 \prisv \exists \varphi_B\losp \exists \bar \rho\losp
(\Base(\varphi_B) \logic\land \Idem^*(\bar \rho)
\logic\land \forall \rho\losp
(\Idem^*(\rho) \logic\land o(\rho)> o(\bar \rho)
\\
& \quad \Rightarrow
\forall f\losp (\mathrm{Ord}_\rho(f)\Rightarrow \Fin(f)
\logic\lor \exists h\losp
(\forall f_1(\Idem^*(f_1) \logic\land o(f_1)> o(\rho)
\\
& \quad \logic\land
\forall c\in Z\losp (cf_1\ne 0\Rightarrow cf_1\varphi_B\ne 0)
\Rightarrow \exists f_2\losp
(\Idem^*(f_2) \logic\land f_2f=f_2 \logic\land f_1h=hf_2=f_1hf_2\ne 0)))))),
\end{align*}
which means that there exists a~number~$k$ such that in a~basic
subgroup~$B$ for every natural $n>k$,
every direct summand which is a~sum of cyclic groups of order~$p^n$
either is finite or has the same power
as the direct summand of the group~$B$ generated by all generators
of order greater than~$p^k$.

Naturally, this means that the final rank of the group~$B$ is countable.

Now if the rings $\Endom(A_1)$ and $\Endom(A_2)$
satisfy the sentence $\neg \psi^3 \logic\land \neg \psi^5$,
then the final ranks of basic subgroups of $A_1$~and~$A_2$
are uncountable and do not coincide with their ranks.
In this case, the theorem follows from Proposition~\ref{p7.3}.

If the rings $\Endom(A_1)$ and $\Endom(A_2)$ satisfy the sentence
$\neg \psi^3 \logic\land \psi^5$, then the final ranks of basic subgroups
of $A_1$~and~$A_2$ are countable and do not coincide with their ranks,
and in this case the theorem follows from Sec.~\ref{sec7_6}.
\end{proof}

\end{document}